\documentstyle[amssymb,twoside,graphicx,epsfig,epstopdf]{asme_submission}

\newtheorem{remark}{Remark}
\newtheorem{theorem}{Theorem}
\newtheorem{lemma}[theorem]{Lemma}
\newtheorem{corollary}[theorem]{Corollary}

\input{tcilatex}
\begin{document}

\author{ Andreas M\"uller 
\affiliation{University of Michigan - Shanghai Jiao Tong University Joint Institute,\\
	Shanghai, China,
	andreas.mueller@ieee.org} \and Zdravko Terze 
\affiliation{Faculty of Mechanical Engineering and Naval Architecture   \\
		University of Zagreb, Croatia,
		zdravko.terze@fsb.hr} \vspace{-10mm} }
\title{ \vspace{-9ex} The significance of the configuration space Lie group
for the constraint satisfaction in numerical time integration of multibody
systems }
\maketitle

\begin{abstract}
The dynamics simulation of multibody systems (MBS) using spatial velocities
(non-holonomic velocities) requires time integration of the dynamics
equations together with the kinematic reconstruction equations (relating
time derivatives of configuration variables to rigid body velocities). The
latter are specific to the geometry of the rigid body motion underlying a
particular formulation, and thus to the used configuration space (c-space).
The proper c-space of a rigid body is the Lie group $SE\left( 3\right) $,
and the geometry is that of the screw motions. The rigid bodies within a MBS
are further subjected to geometric constraints, often due to lower kinematic
pairs that define $SE\left( 3\right) $ subgroups. Traditionally, however, in
MBS dynamics the translations and rotations are parameterized independently,
which implies the use of the direct product group $SO\left( 3\right) \times 
{\Bbb R}^{3}$ as rigid body c-space, although this does not account for
rigid body motions. Hence, its appropriateness was recently put into
perspective.

In this paper the significance of the c-space for the constraint
satisfaction in numerical time stepping schemes is analyzed for holonomicaly
constrained MBS modeled with the 'absolute coordinate' approach, i.e. using
the Newton-Euler equations for the individual bodies subjected to geometric
constraints. The numerical problem is considered from the kinematic
perspective. It is shown that the geometric constraints a body is subjected
to are exactly satisfied if they constrain the motion to a subgroup of its
c-space. Since only the $SE\left( 3\right) $ subgroups have a practical
significance it is regarded as the appropriate c-space for the constrained
rigid body. Consequently the constraints imposed by lower pair joints are
exactly satisfied if the joint connects a body to the ground. For a general
MBS, where the motions are not constrained to a subgroup, the $SE\left(
3\right) $ and $SO\left( 3\right) \times {\Bbb R}^{3}$ yield the same order
of accuracy. Hence an appropriate configuration update can be selected for
each individual body of a particular MBS, which gives rise to tailored
update schemes. Several numerical examples are reported illustrating this
statement.

The practical consequence of using $SE\left( 3\right) $ is the use of screw
coordinates as generalized coordinates. To account for the inevitable
singularities of 3-parametric descriptions of rotations, the kinematic
reconstruction is additionally formulated in terms of (dependent) dual
quaternions as well as a coordinate-free ODE on the c-space Lie group. The
latter can be solved numerically with Lie group integrators like the
Munthe-Kaas integration method, which is recalled in this paper.

\end{abstract}

\vspace{-6ex}

\section{Introduction}

The seemingly simple problem addressed in this paper is how to numerically
reconstruct the finite motion of a constrained rigid body within a MBS from
its velocity field so that the overall system of geometric constraints is
satisfied. When a rigid body moves it performs a translation together with a
rotation since a general rigid body motion is a screw motion, with coupled
rotation and translation. Even though standard numerical integration schemes
for MBS neglect the geometry of Euclidean motion in the sense that, within
the integration schemes, the position and orientation updates are performed
independently. Whether or not their dependence is respected has to do with
the geometric model used to represent rigid body motions, i.e. with the
configuration space (c-space) Lie group. It is know that rigid body motions
form the Lie group $SE\left( 3\right) $ \cite{murray,SeligBook}.

With the 'absolute coordinate' formalism (i.e. representing the spatial
configuration of each body by a set of six generalized coordinates) the
equations governing the dynamics of a constrained MBS comprising $n$ rigid
bodies are commonly written in the form 
\begin{eqnarray}
{\bf M}\left( {\bf q}\right) \dot{{\bf V}}+{\bf J}^{T}{\bf \lambda } &=&{\bf %
Q}\left( {\bf q},{\bf V},t\right) \ \ \ \ a)  \label{BH} \\
{\bf V} &=&{\bf A}\left( {\bf q}\right) \dot{{\bf q}}\ \ \ \ b)%
\addtocounter{equation}{-1} \\
h\left( {\bf q}\right) &=&{\bf 0}.\ \ \ \ c)\addtocounter{equation}{-1}
\end{eqnarray}%
The $N=6n$ dimensional coordinates vector ${\bf q}=\left( {\bf \theta }_{i},%
{\bf r}_{i}\right) \in {\Bbb V}^{N}$ comprises the position vector ${\bf r}%
_{i}$ and the vector ${\bf \theta }_{i}$ consisting of 3 (or 4 dependent)
rotation parameters for body $i=1,\ldots ,n$, and ${\bf V}=\left( {\bf %
\omega }_{i},{\bf v}_{i}\right) \in {\Bbb R}^{N}$ is composed of the angular
and linear velocity vectors ${\bf \omega }_{i}$ and ${\bf v}_{i}$,
respectively. The matrix ${\bf J}$ is the constraint Jacobian corresponding
to the system (\ref{BH}c) of geometric constraints.

The equations (\ref{BH}) constitute a DAE system on the coordinate manifold $%
{\Bbb V}^{N}$ considered as vector space. From a kinematic point of view
this formulation raises two issues regarding their numerical solution:

\begin{enumerate}
\item The motion of the MBS is deduced from the velocity ${\bf V}$ by the 
{\it kinematic reconstruction equations} (\ref{BH}b). The accuracy of their
numerical solution depends directly on the underlying geometry of rigid body
motions, which is encoded in the mapping ${\bf A}$. In the standard MBS
formulation the rotations and positions are reconstructed separately
according to \vspace{-3ex}

\begin{equation}
\left( 
\begin{array}{c}
{\bf \omega }_{i} \\ 
{\bf v}_{i}^{\text{s}}%
\end{array}%
\right) =\left( 
\begin{array}{cc}
{\bf B}_{i}\left( {\bf \theta }_{i}\right) & {\bf 0} \\ 
{\bf 0} & {\bf I}%
\end{array}%
\right) \left( 
\begin{array}{c}
\dot{{\bf \theta }}_{i} \\ 
\dot{{\bf r}}_{i}%
\end{array}%
\right) ,i=1,\ldots ,n.  \label{VelRigidBody}
\end{equation}%
\vspace{-3ex}

The underlying geometry is that of $SO\left( 3\right) \times {\Bbb R}^{3}$,
which does not account for the coupling of rotations and translation
inherent to screw motions. Nevertheless, the kinematic equations (\ref%
{VelRigidBody}) correspond to a valid parameterization of rigid body {\em %
configurations}. The interdependence of ${\bf \omega }_{i}$ and ${\bf v}%
_{i}^{\text{s}}$ is ensured by solving (\ref{BH}a) and (\ref{BH}c), and an
analytic solution of (\ref{VelRigidBody}) correctly reflects the bodies'
screw motions. However, when (\ref{BH}) is solved numerically with a finite
step size, and (\ref{BH}b) is used to predict finite (screw) {\em motion}
increments, also the kinematic reconstruction equations (\ref{BH}b) must
properly reflect the geometry of screw motions. Moreover, (\ref{VelRigidBody}%
) can only predict the finite motion if ${\bf r}_{i}$ are the coordinates of
a point on the rotation axis, as for instance in the case of an
unconstrained body with its body-fixed reference frame located at the COM. A
generic motion of a constrained body, as part of a MBS, will not comply with
the decoupling assumption encoded in (\ref{VelRigidBody}). The matrix ${\bf B%
}$ is specific to the rotation parameterization. If Euler angles are used,
for instance, ${\bf B}$ corresponds to the kinematic Euler equations \cite%
{McCauley}.\newline
Consequently, the kinematic reconstruction equations (\ref{BH}b) shall be
amended in order to respect the interrelation of rotations and translations,
which boils down to the appropriate choice of the rigid body configuration
space being a Lie group. The implications of using the Lie group $SE\left(
3\right) $ as well as $SO\left( 3\right) \times {\Bbb R}^{3}$ are studied in
this paper.

\item The second issue regards the violations of the constraints (\ref{BH}c)
that occur when numerically solving (\ref{BH}). This has been a central
problem in numerical MBS dynamics. However, the investigations have
exclusively been focussed on reducing or correcting constraint violations by
means of stabilization and projection methods \cite%
{BauchauLaulusa2008,Blajer2011} rather then aiming to avoid such violations.
It is immediately clear that the constraint satisfaction is affected by the
accuracy with which the finite motions are reconstructed from the velocity
field ${\bf V}$ solving (\ref{BH}a), which indeed depends on the feasibility
of the relation (\ref{BH}b). Even more, besides the accuracy with which the
system dynamics is captured by the numerical integrator, it is crucial to
ensure the kinematic consistency of the MBS, thus the constraint
satisfaction is imperative. This is the focus of this paper. In this respect
it is important to observe that the majority of mechanisms is built with
lower kinematic pairs (Reuleaux pairs). The latter are characterized by
their isotropy groups, i.e. subgroups of $SE\left( 3\right) $ leaving the
contact surface invariant. It is clear that, if a numerical update step does
not respect these motion groups, the lower pair constraints will be violated.
\end{enumerate}

The reconstruction equations (\ref{BH}b) represent a first-order relation,
and from a computational perspective the question arises whether the
decoupling significantly affects the accuracy of the numerical solution of (%
\ref{BH}). The goal of this paper is to study the extend to which different
forms of this first-order relation affect the reconstruction of finite
motions of a constrained MBS with numerical time integration methods. To
this end, the kinematic reconstruction equations (\ref{BH}b) will be
reformulated as an ordinary differential equation (ODE) on the c-space Lie
group of the form \vspace{-5ex}

\begin{equation}
\widehat{{\bf V}}_{i}={\bf dexp}_{-{\bf \Phi }_{i}}\dot{{\bf \Phi }}_{i}
\label{Videxp}
\end{equation}%
with ${\bf \Phi }_{i},\widehat{{\bf V}}_{i}\in {\frak g}$ and ${\frak g}$
being the Lie algebra of the respective c-space Lie group used to describe
the motion of body $i$. The two Lie groups considered here are $SE\left(
3\right) $ and $SO\left( 3\right) \times {\Bbb R}^{3}$. In case of $SE\left(
3\right) $ the exponential coordinates ${\bf \Phi }_{i}$ are the screw
coordinates used to describe the rigid body motions, while in case of $%
SO\left( 3\right) \times {\Bbb R}^{3}$ they consist of the scaled rotation
vector and the Cartesian position coordinates, and (\ref{Videxp}) reduces to
(\ref{VelRigidBody}). The operator ${\bf dexp}$ is the right-translated
differential of the exponential mapping on the c-space Lie group. Hence, the
use of a specific c-space is reflected by the corresponding exponential
mapping, which clearly leads to different numerical properties of (\ref%
{Videxp}).

Spatial kinematics rests on the theory of screws \cite%
{BottemaRoth1979,Dimentberg}. Instantaneous screw motions form the Lie
algebra $se\left( 3\right) $ that generates the Lie group $SE\left( 3\right) 
$ of finite rigid body motions. Both are canonically related via the
exponential mapping. This fact is the cornerstone of modern kinematics. In
context of computational MBS dynamics this fact has not yet been
sufficiently exploited, however. Chevalier \cite{chevallier84,chevallier94}
introduced a coordinate free formulation for the dynamics of tree-topology
MBS, which has unfortunately not been recognized by the MBS community. In
robotics the Lie group formulation has been given due attention after
Brocket \cite{Brockett1984} introduced the product-of-exponential (POE)
formula. This gave rise to Lie group formulations for MBS such as \cite%
{Mladenova1999,Mladenova2006,mueller2003,PloenPark1997,ParkKim2000}. Liu in 
\cite{Liu1988} already presented a MBS modeling approach using screw theory
and Lie group methods. Now excellent introductions can be found in the text
books by Angeles \cite{Angeles2003}, Murray et. al \cite{murray}, and Selig 
\cite{SeligBook}. On one hand, the Lie group setting provides a framework
for systematic matrix formulations of MBS dynamics as discussed in the
recent book by Featherstone \cite{Featherstone2008} and Uicker et al. \cite%
{Uicker2013} that are tailored for code implementation. On the other hand,
it also provides the basis for analytical investigations as reported in \cite%
{EmamiMMT}.

These contributions regard the dynamics modeling of MBS. Even less
publications deal with the numerical solution of the MBS motion equations
using geometric concepts or Lie group time integration schemes such as the
Munthe-Kaas (MK) scheme. The time integration of the rotational dynamics of
rigid bodies on the Lie group $SO\left( 3\right) $ has been a standard
example from the outset of the MK method in \cite%
{Marthinsen1997,muntekaas1998,muntekaas1999,MuntheKaasOwren1999,OwrenMarthinsen1999}%
. But the MK scheme is only applied to MBS dynamics in few publications such
as \cite{ParkChung2005}. In that publication the rigid body motion is
considered on $SE\left( 3\right) $, and the dexp mapping in (\ref{Videxp})
is the one for screw motions.

Recently there is an increased interest in applying Lie group integrators to
MBS, mainly motivated by the fact that no global parameterization is needed,
thus circumventing the singularity problem of spatial rotations, and that
the geometric setting allows preservation of essential invariants. In \cite%
{BruelsCardona2010,BruelsCardonaArnold2012,Krysl} numerical time stepping
schemes were introduced adopting the Newmark/generalized-alpha schemes \cite%
{Erlicher2002,Newmark1959}. A MK integrator is used in \cite{TerzeIMSD2012}
to solve the index 1 DAE MBS formulation. All these approaches use $SO\left(
3\right) \times {\Bbb R}^{3}$ as rigid body c-space. The latter does not
represent proper rigid body motions. The consequences of this fact for the
constraint satisfaction in MBS models have been briefly analyzed in \cite%
{MuellerTerze2013}. Along this line numerical results were reported in \cite%
{BruelsCardonaArnold2011}, but without further analysis. As it will be shown
in this paper, the crucial point is the preservation of invariants of rigid
body motions (only preserved by $SE\left( 3\right) $). The preservation of
invariants is also a central issue in the dynamics of flexible bodies.
Geometrically exact formulations rely on the proper representation of large
motions, and numerical time stepping schemes should respect the geometric
invariants. This has been discussed in \cite{Bottasso2002} where it was
shown that a numerical scheme is only invariant w.r.t. to the chosen
reference frame if the motion of a frame (rigid body, nodal frame) is
represented by $SE\left( 3\right) $ (although the Lie group is not
mentioned). This is the basis for momentum-preserving integration schemes
for finite element models reported in \cite%
{BorriBottasso1994,Borri2001a,Borri2001b,Borri2002}, and for an integration
scheme in \cite{Borri2003} that also respects Newton's third law
(actio-reactio principle). Unfortunately the implications of these results
for time integration schemes have not been sufficiently recognized.
Moreover, occasionally. $SE\left( 3\right) $ is wrongly identified with $%
SO\left( 3\right) \times {\Bbb R}^{3}$ in the literature, e.g. \cite%
{RobHandbook}.

The aim of this paper is to clarify that any Lie group integrator operating
on $SO\left( 3\right) \times {\Bbb R}^{3}$ as rigid body c-space violates
the joint constraints in MBS models according to its order of accuracy,
whereas using $SE\left( 3\right) $ can achieve perfect constraint
satisfaction. The latter does not increase the order of accuracy of the
applied numerical integration scheme, but it ensures consistency of the MBS
model and thus improves the overall accuracy. The problem is addressed for
the absolute coordinates (also called inertial coordinates) formulation.
where (\ref{BH}a) consists of the Newton-Euler equations of the $n$
individual bodies and (\ref{BH}c) represents geometric constraints due to
lower pair joints. The problem is approached from a kinematic perspective.
Numerical aspects of the integration schemes are not a topic of this paper.

The paper is organized as follows. In section 2 the configuration space of a
constrained MBS is introduced using $SE\left( 3\right) $ and $SO\left(
3\right) \times {\Bbb R}^{3}$ as Lie group, respectively. The MBS motion
equations are given in section 3. The classical vector space formulation is
presented in terms of independent canonical parameters (screw coordinates).
To account for parameterization singularities, a formulation in terms of
dual quaternions is also given. Further, a Lie group formulation is reported
that is inherently coordinate-free and not prone to singularities. As a
numerical integration scheme for this formulation the Munthe-Kaas method is
recalled in appendix B. The main result is presented in section \ref%
{secDiscussion} where the consequences of using either c-space Lie group for
the constraint satisfaction are discussed. It is concluded that $SE\left(
3\right) $ yields the best constraint satisfaction at the expense of
slightly more complex kinematic relations. Several examples are reported in
section \ref{secExamples}, and the paper closes with a summary and
conclusion in section \ref{secConclusion}. Regarding the geometric
background the reader is referred to the text books \cite%
{Angeles2003,BottemaRoth1979,Dimentberg,murray,SeligBook}. In order to make
the paper accessible to reader not familiar with rigid body kinematics and
screw theory, the fundamental geometric background is summarized in appendix
A. In particular, the semidirect and direct product representation of rigid
body motions are discussed. Their parameterization in terms of dual
quaternions and ordinary quaternions are recalled. This is important since
the parameterization of $SE\left( 3\right) $ in terms of canonical (screw)
coordinates suffers from the known singularities, and would thus not allow
for using the $SE\left( 3\right) $ configuration update in the classical
vector space MBS formulation.

It is important to emphasize that all considerations in MBS kinematic are
related to certain reference frames and that all Lie groups appearing in
this context are matrix Lie groups.

\section{The Configuration Space Lie Group of a Constrained MBS}

\subsection{MBS C-Space in Left-Invariant $SE\left( 3\right) $ Representation%
}

The Lie group of proper rigid body motions is $SE\left( 3\right) =SO\left(
3\right) \ltimes {\Bbb R}^{3}$ --the semidirect product of the rotation
group and the translation group identified with ${\Bbb R}^{3}$ (see appendix
A). The {\em ambient configuration space }of an MBS comprising $n$ rigid
bodies is the $6n$-dimensional Lie group (the superscript $\ltimes $
indicates the semidirect product) 
\begin{equation}
G^{\ltimes }:=SE\left( 3\right) ^{n}.
\end{equation}%
An element $g=({\bf C}_{1},\ldots ,{\bf C}_{n})\in G^{\ltimes }$ represents
the configuration of $n$ decoupled bodies in a coordinate-free way.
Multiplication is componentwise and inherited from $SE\left( 3\right) $. The
inverse is $g^{-1}=({\bf C}_{1}^{-1},\ldots ,{\bf C}_{n}^{-1})$.

The Lie algebra of $G^{\ltimes }$ is ${\frak g}^{\ltimes }:=se\left(
3\right) ^{n}$ equipped with the Lie bracket inherited componentwise from $%
se\left( 3\right) $ in (\ref{se3bracket}). ${\frak g}^{\ltimes }$ is
isomorphic to $\left( {\Bbb R}^{6}\right) ^{n}$ equipped with the
componentwise screw product. The notations ${\bf q}=\left( {\bf X}%
_{1},\ldots ,{\bf X}_{n}\right) \in \left( {\Bbb R}^{6}\right) ^{n}$ and $%
\widehat{{\bf q}}=(\widehat{{\bf X}}_{1},\ldots ,\widehat{{\bf X}}_{n})\in 
{\frak g}^{\ltimes }$ are used to denote the overall vector of screw
coordinates of the MBS. With the exponential mapping on $SE\left( 3\right) $
the mapping $\exp :{\frak g}^{\ltimes }\rightarrow G^{\ltimes }$ is 
\begin{equation}
\exp \widehat{{\bf q}}=(\exp \widehat{{\bf X}}_{1},\ldots ,\exp \widehat{%
{\bf X}}_{n})\in G^{\ltimes }.
\end{equation}%
Accordingly its right-translated differential is ${\rm dexp}_{\widehat{{\bf q%
}}}\left( \widehat{{\bf q}}^{\prime }\right) =({\rm dexp}_{\widehat{{\bf X}}%
_{1}}(\widehat{{\bf Y}}_{1}^{\prime }),\ldots ,{\rm dexp}_{\widehat{{\bf X}}%
_{n}}(\widehat{{\bf Y}}_{n}^{\prime }))$. The vector ${\bf q}$ represents
the global coordinates of the MBS.

MBS velocities are denoted as{\bf \ }${\bf V}=({\bf V}_{1},\ldots ,{\bf V}%
_{n})\in \left( {\Bbb R}^{6}\right) ^{n}$. The body-fixed MBS twists are
introduced as $\widehat{{\bf V}}=g^{-1}\dot{g}\in {\frak g}^{\ltimes }$. In
terms of $\dot{{\bf x}}$ this is $\widehat{{\bf V}}={\rm dexp}_{-\widehat{%
{\bf q}}}(\dot{\widehat{{\bf q}}})$, which can be written in vector form as $%
{\bf V}={\bf dexp}_{-{\bf q}}\dot{{\bf q}}$.

The ambient c-space Lie group allows representing configurations of $n$
decoupled bodies. It remains to incorporate constraints the MBS is subjected
to. It is assumed that the MBS is subjected to a system of $m$ scleronomic
geometric constraints 
\begin{equation}
h\left( g\right) ={\bf 0}  \label{geomConstr}
\end{equation}%
with the constraint mapping $h:G^{\ltimes }\rightarrow {\Bbb R}^{m}$. Time
differentiation of the geometric constraints gives rise to the velocity
constraints 
\begin{equation}
{\bf J}\left( g\right) \cdot {\bf V}={\bf 0}  \label{VelConstr}
\end{equation}%
where ${\bf J}\left( g\right) :\left( {\Bbb R}^{6}\right) ^{n}\rightarrow 
{\Bbb R}^{m}$ is the Jacobian mapping of $h$ in vector space representation
of ${\frak g}^{\ltimes }$. The geometric constraints define the {\it MBS
configuration space} (an analytic subvariety of the manifold $G^{\ltimes }$) 
\begin{equation}
{\cal V}^{\ltimes }:=\{g\in G^{\ltimes }|h\left( g\right) ={\bf 0}\}.
\label{VSE3}
\end{equation}%
(\ref{geomConstr}), (\ref{VelConstr}), and (\ref{VSE3}) can be expressed in
terms of the (local) coordinates ${\bf q}\in {\frak g}^{\ltimes }$ on $%
G^{\ltimes }$ (screw coordinates). Alternatively, dual quaternions can be
used as dependent (global) coordinates, replacing $SE\left( 3\right) $ by $%
{\Bbb H}_{\varepsilon }$ (appendix \ref{secDualQuat}). The couple $\left( g,%
{\bf V}\right) $, or alternatively $\left( {\bf q},{\bf V}\right) $,
represents the MBS state. The Lie group formulation is merely a formal
description of what is commonly pursued in MBS dynamics.

\begin{remark}
In this paper the 'absolute coordinate' approach for MBS modeling is used.
In coordinate-free terms the 'absolute' configuration of each body w.r.t. a
global reference frame is represented by ${\bf C}_{i}$, and these are
subject to certain joint constraints. Rotation and translation parameters
can be introduced for each rigid body representing six 'absolute
coordinates'.
\end{remark}

\subsection{MBS C-Space in Left-Invariant $SO\left( 3\right) \times {\Bbb R}%
^{3}$ Representation}

The direct product group $SO\left( 3\right) \times {\Bbb R}^{3}$ allows for
representing rigid body {\em configurations} but not rigid body {\em motions}%
. The direct product representation corresponds to the mixed velocity
representation. The corresponding {\it ambient configuration} space of the
MBS is the $6n$-dimensional Lie group (the superscript $\times $ indicates
the direct product) 
\begin{equation}
G^{\times }:=\left( SO\left( 3\right) \times {\Bbb R}^{3}\right) ^{n}
\end{equation}%
with elements $g=({\bf C},\ldots ,{\bf C}_{n})\in G^{\times }$ and inverse $%
g^{-1}=({\bf C}_{1}^{-1},\ldots ,{\bf C}_{n}^{-1}))$. The multiplication is
inherited from $SO\left( 3\right) \times {\Bbb R}^{3}$. The corresponding
Lie algebra is ${\frak g}^{\times }:=\left( so\left( 3\right) \oplus {\Bbb R}%
^{3}\right) ^{n}$, which is isomorphic to $\left( {\Bbb R}^{6}\right) ^{n}$
equipped with the componentwise Lie bracket (\ref{bracketSO3xR3}). Adopting
the above notation, ${\bf q}=\left( {\bf X}_{1},\ldots ,{\bf X}_{n}\right)
\in \left( {\Bbb R}^{6}\right) ^{n}$ and $\widehat{{\bf q}}=(\widehat{{\bf X}%
}_{1},\ldots ,\widehat{{\bf X}}_{n})\in {\frak g}^{\ltimes }$, the
exponential mapping on $G^{\times }$ is introduced as 
\begin{equation}
\exp {\bf q}=(\exp {\bf X}_{1},\ldots ,\exp {\bf X}_{n})\in G^{\times }
\end{equation}%
with exp in (\ref{expSO3xR3}), and its differential in an obvious way. The
mixed velocity of the MBS is defined with (\ref{hybridvel}) as the
left-invariant vector field $\widehat{{\bf V}}^{\text{m}}=g^{-1}\dot{g}=(%
\widehat{{\bf V}}_{1}^{\text{m}},\ldots ,\widehat{{\bf V}}_{n}^{\text{m}%
})\in {\frak g}^{\times }$. It is determined in vector notation as ${\bf V}^{%
\text{m}}={\bf dexp}_{-{\bf q}}\dot{{\bf q}}$ with (\ref{Vhdexp}).

Imposing the geometric constraints (\ref{geomConstr}) yields the {\it MBS
configuration space in direct product representation} 
\begin{equation}
{\cal V}^{\times }:=\{g\in G^{\times }|h\left( g\right) ={\bf 0}\}
\end{equation}%
The velocity constraints attain the form (\ref{VelConstr}), but now with $%
{\bf V}^{\text{m}}$. The ${\bf q}\in {\frak g}^{\times }$ serve as (local)
coordinates on $G^{\times }$. The MBS state is represented by $\left( g,{\bf %
V}^{\text{m}}\right) $ or $\left( {\bf q},{\bf V}^{\text{m}}\right) $. Using
quaternions for the rotation part gives rise to an alternative (global)
parameterization when $SO\left( 3\right) \times {\Bbb R}^{3}$ is replaced by 
${\Bbb H}\times {\Bbb R}^{3}$ (appendix \ref{secQuat}).

\section{Motion Equations of Constrained MBS}

\subsection{Motion equations on the parameter space in canonical coordinates
(axis-angle, position vector)}

The constrained motion equations (\ref{BH}) in terms of the non-holonomic
quasi-velocities ${\bf V}$ are used as governing equations in descriptor
form for the dynamics simulation of constrained MBS in terms of absolute
coordinates ${\bf q}$. To account for a general c-space Lie group, the
kinematic reconstruction equations (\ref{BH}b) are expressed in terms of the
dexp mapping, which gives rise to the formulation 
\begin{eqnarray}
{\bf M}\left( {\bf q}\right) \dot{{\bf V}}+{\bf J}^{T}{\bf \lambda } &=&{\bf %
Q}\left( {\bf q},{\bf V},t\right) \ \ \ \ \ \ \ \ \ \ \ \ \ \ \ a)
\label{BH2} \\
{\bf V} &=&{\bf dexp}_{-{\bf q}}\dot{{\bf q}}\ \ \ \ \ \ \ \ \ \ \ \ \ b)%
\addtocounter{equation}{-1} \\
h\left( {\bf q}\right) &=&{\bf 0}.\ \ \ \ \ \ \ \ c)%
\addtocounter{equation}{-1}
\end{eqnarray}%
The choice of 'absolute coordinates' depends on the c-space Lie group. Using 
$G^{\ltimes }$ as c-space, ${\bf q}\in {\frak g}^{\times }$ comprises the $n$
screw coordinate vectors ${\bf X}_{i}=\left( {\bf \xi }_{i},{\bf \eta }%
_{i}\right) $, and when using $G^{\times }$, ${\bf q}\in {\frak g}^{\times }$
comprises the $n$ vectors ${\bf X}_{i}=\left( {\bf \xi }_{i},{\bf r}%
_{i}\right) $. In either case, ${\bf q}$ are local coordinates and ${\frak g}%
\simeq {\Bbb R}^{6n}$ is the parameter space. Accordingly, ${\bf V}$
consists of body-fixed twists or mixed velocities. With the 'absolute
coordinate approach' the system (\ref{BH2}a) summarizes the Newton-Euler
equations governing the dynamics of the $n$ bodies subject to the geometric
constraints (\ref{BH2}c). This an index 3 DAE system. The dynamics of the
MBS is represented as the dynamics of the representing point ${\bf q}$ on
the subvariety $h^{-1}\left( {\bf 0}\right) \subset G$. It is commonly
treated as an ODE on a vector space by considering the parameter space as
vector space. To this end (\ref{BH2}a) and the acceleration constraints $%
{\bf J}\left( {\bf q}\right) \cdot \dot{{\bf V}}={\bf \eta }\left( {\bf q},%
{\bf V}\right) $ (second time derivative of (\ref{BH2}c)) are combined to
the index 1 system 
\begin{equation}
\left( 
\begin{array}{cc}
{\bf M} & {\bf J}^{T} \\ 
{\bf J} & {\bf 0}%
\end{array}%
\right) \left( 
\begin{array}{c}
\dot{{\bf V}} \\ 
{\bf \lambda }%
\end{array}%
\right) =\left( 
\begin{array}{c}
{\bf Q} \\ 
{\bf \eta }%
\end{array}%
\right) .  \label{index1}
\end{equation}%
For a given state $\left( {\bf q},{\bf V}\right) $ of the MBS, (\ref{index1}%
) can be solved for $\dot{{\bf V}}$ and ${\bf \lambda }$ giving rise to an
explicit ODE $\dot{{\bf V}}=F\left( {\bf q},{\bf V}\right) $ and thus
replacing (\ref{BH2}) by 
\begin{eqnarray}
\dot{{\bf V}} &=&F\left( {\bf q},{\bf V}\right) \ \ \ \ \ \ \ \ \ \ \ \ \ \
\ \ a)  \label{ODEvec1} \\
\dot{{\bf q}} &=&{\bf dexp}_{-{\bf q}}^{-1}{\bf V}\ \ \ \ \ \ \ \ \ \ \ \ \
\ \ \ b)\addtocounter{equation}{-1}.
\end{eqnarray}%
This general formulation of MBS motion equations is applicable to any
c-space Lie group. The velocities ${\bf V}$ are either body-fixed twists or
mixed velocities. A specific choice for the rigid body c-space, with
corresponding dexp mapping, leads to a specific ODE system (\ref{ODEvec1}b).
If $G^{\times }$ and a 3-parametric description of rotations are used, (\ref%
{ODEvec1}b) attains the particular form (\ref{VelRigidBody}) with (\ref%
{3param}), which is commonly applied in MBS formulations. In any case the
motion equations are treated as a system on ${\Bbb V}^{n}$ considered as
vector space.

\begin{remark}
Analytically any kinematic relation in (\ref{BH2}b) is admissible that
corresponds to a valid parameterization of rigid body configurations (not
necessarily motions). Even if the coupling of angular and linear velocities
are not respected by (\ref{BH2}b) and (\ref{ODEvec1}b), as with (\ref{Vhdexp}%
), any kinematic reconstruction equation will yield the same \underline{%
analytic} solution since this dependence is respected by the solution of (%
\ref{BH2}a). But this is not true when numerically solving (\ref{BH2}) as
shown in the remainder of this paper since then (\ref{BH2}b) is used to 
\underline{estimate} finite motions.
\end{remark}

{\bf Assumption:} It is assumed in the following that the velocity
constraints are satisfied.

This assumption is made in order to focus on the configuration update. The
consistency of the velocity with the linear velocity constraints can be
easily achieved by a one-step projection if necessary. In all examples
below, and for the majority of complex MBS, the numerical solution for ${\bf %
V}$ at time step $i$ satisfies the velocity constraints when starting with
one at time step $i-1$ that satisfies them.

\subsection{Motion equations in dependent coordinates (quaternions)}

Any parameterization of rigid bodies in terms of six independent parameters
leads to parameterization singularities, caused by 3-parameter description
of rotations. To circumvent this problem, dual quaternions and Euler
parameters can be used giving rise to a singularity-free parameterization of 
$SE\left( 3\right) $ and $SO\left( 3\right) \times {\Bbb R}^{3}$,
respectively (appendix \ref{secHybrid} and \ref{secDualQuat}). In terms of
these parameters, summarized in the generalized coordinate vector ${\bf q}$,
the motion equations are 
\begin{eqnarray}
\dot{{\bf V}} &=&F\left( {\bf q},{\bf V}\right) \ \ \ \ \ \ \ \ \ \ \ \ \ \
\ \ a)  \label{ODEvec2} \\
\dot{{\bf q}} &=&{\bf H}^{-1}({\bf q}){\bf V}\ \ \ \ \ \ \ \ \ \ \ \ \ \ \ \
b)\addtocounter{equation}{-1}
\end{eqnarray}%
where the coefficient matrix in (\ref{ODEvec2}b) is either (\ref{Hdual}) for 
$SE\left( 3\right) $, or (\ref{Hh}) for $SO\left( 3\right) \times {\Bbb R}%
^{3}$, or (\ref{Hhdual}) if mixed velocities are used in the $SE\left(
3\right) $ formulation. The use of dependent parameters introduces further
constraints due to the normalization and Pl\"{u}cker condition. These
constraints are included in (\ref{BH2}c), and are accounted for by the
constraint Jacobian. The motion equations are again treated as system on a
vector space.

It is straightforward to evaluate $F$ in (\ref{ODEvec2}a) when the rotation
matrices of the bodies are replaced by ${\bf R}={\bf DE}^{T}$ in (\ref{DE})
and the position vectors by ${\bf r}$ in (\ref{Qpos}). It is just remarked
here that the motion equations can also be consistently formulated in terms
of dual quaternions \cite{shoam,DooleyMcCarthy1991}.

\subsection{Motion equations on the c-space Lie group%
\label{secEOMLieGroup}%
}

The best way to avoid parameterization singularities is to not use local
parameters at all. Instead of introducing local coordinates ${\bf q}$, rigid
body configurations are represented coordinate-free by $g\in G^{\ltimes }$
or $g\in G^{\times }$. The velocities ${\bf V}$ belong to the respective Lie
algebra. Then the constrained motion equations attain the form 
\begin{eqnarray}
{\bf M}\left( g\right) \dot{{\bf V}}+{\bf J}^{T}{\bf \lambda } &=&{\bf Q}%
\left( g,{\bf V},t\right) {\ \ \ \ \ \ \ \ \ \ \ \ \ \ \ \ a)}  \label{BH3}
\\
\dot{g} &=&g\widehat{{\bf V}}\ \ \ \ \ \ \ \ \ \ \ \ \ \ \ b)%
\addtocounter{equation}{-1} \\
h\left( g\right) &=&{\bf 0}.\ \ \ \ \ \ \ \ \ \ \ \ \ \ \ c)%
\addtocounter{equation}{-1}
\end{eqnarray}%
where, instead of using a parameterization, the reconstruction equations (%
\ref{BH3}b) are merely the (coordinate-free) definition of twists, either (%
\ref{Vhat}) or (\ref{hybridvel}). With a solution $\dot{{\bf V}}=F\left( g,%
{\bf V}\right) $ of (\ref{index1}) this leads to the motion equations 
\begin{eqnarray}
\dot{{\bf V}} &=&F\left( g,{\bf V}\right) {\ \ \ \ \ a)}  \label{ODE2} \\
\dot{g} &=&g\widehat{{\bf V}}.{\ \ \ \ \ \ \ \ \ \ \ \ \ \ \ \ b)}%
\addtocounter{equation}{-1}
\end{eqnarray}%
Overall (\ref{ODE2}) is an ODE system on the MBS state space ${\frak g}%
\times G$. The system (\ref{ODE2}b) is an ODE on the c-space Lie group $G$,
and equations (\ref{ODE2}a) form an ODE on the vector space ${\frak g}\cong 
{\Bbb R}^{6n}$, as in (\ref{ODEvec1}). The latter can be solved with any
established numerical (vector space) integration scheme, whereas solution of
(\ref{ODE2}b) requires application of a Lie group integration scheme, where
local parameters will have to be used eventually. In this way the
parameterization problem is shifted to the integration scheme \cite%
{BruelsCardona2010,BruelsCardonaArnold2012,EngoMarthinsen2001,Krysl,TerzeIMSD2012}%
. A widely used integration scheme for ODE on Lie groups is the explicit MK
method, which is briefly recalled in appendix B. These schemes implicitly
define a local parameter chart around the current configuration.

\section{The significance of the configuration space for numerical
reconstruction of rigid body motions constrained to subgroups%
\label{secDiscussion}%
\label{secDiscussion}}

The fact that most joints in MBS models constrain the relative motion of
adjacent bodies to $SE\left( 3\right) $ subgroups (in particular Reuleaux
lower pairs) seems to put the use of the direct product group into
perspective. Since in the absolute coordinate formulation the update step
adjusts the absolute configuration bodies, the semidirect representation
does not necessarily preserve the lower pair constraints, however.

Given the velocity of a body, its finite motion is the solution of the
kinematic reconstruction equations (\ref{VdexpSE3}) (respectively (\ref%
{Vhdexp})). Since either formulation is consistent with the respective
definition of velocity, both yield the same {\it analytical} solution for
the MBS motion (but in terms of different coordinates). The difference of
the formulations is critical, however, when the reconstruction equations are
solved numerically and used to estimate finite motions.

When numerically solving (\ref{ODEvec1}b) with a (explicit or implicit)
numerical vector space integration scheme, a constant screw coordinate
vector ${\bf \Phi }^{\left( i\right) }\in {\frak g}$ is constructed and used
to estimate the finite configuration at integration step $i$ as $g\left(
t_{i}\right) =\exp \widehat{{\bf q}}^{\left( i\right) }$ with ${\bf q}%
^{\left( i\right) }={\bf q}^{\left( i-1\right) }+{\bf \Phi }^{\left(
i\right) }$. When integrating (\ref{ODE2}b) with a Lie group integrator such
as a Munthe-Kaas (MK) integration scheme (\ref{MK1}), a ${\bf \Phi }^{\left(
i\right) }$ is constructed that determines the configuration increment $\exp 
\widehat{{\bf \Phi }}^{\left( i\right) }$. Hence, regardless of the
particular integration scheme, a finite estimate is found by evaluating the
dexp mapping. The latter locally reflects the geometry of the configuration
space, and clearly using either $SO\left( 3\right) \times {\Bbb R}^{3}$ or $%
SE\left( 3\right) $, the respective dexp mapping yields different changes of
the local coordinates ${\bf \Phi }$ for the same twist, which is uniquely
determined by the dynamic equations (\ref{ODEvec1}a). Nevertheless the
finite motion estimates are in accordance with the respective geometric
model, but not necessarily with the MBS kinematics. Thus the difference
boils down to the geometry of the c-space. Consequently, only the $SE\left(
3\right) $ update can possibly lead to proper estimates of constrained rigid
body motions.

A numerical integration scheme estimates ${\bf \Phi }^{\left( i\right) }$ as
a (usually) linear combination of ${\bf dexp}^{-1}$ in (\ref{ODEvec1}b) when
evaluated at intermediate time steps $j=1,\ldots ,s$ of the (possibly
implicit) numerical integrator. Using (\ref{dexpInvSeries}), which applies
to any Lie group, a general integration step can be written as%
\begin{equation}
{\bf \Phi }^{\left( i\right) }=\sum_{j=1}^{s}\sum_{r=0}^{\infty }\alpha _{ij}%
{\bf ad}_{{\bf V}^{\left( j\right) }}^{r}{\bf V}^{\left( i\right) }
\end{equation}%
denoting ${\bf V}^{\left( j\right) }:={\bf V}(t_{i-1}+c_{j}\Delta t,{\bf q}%
_{j})$ with some real coefficients $\alpha _{rl}$ and $c_{j}$ specific to
the integration scheme. That is, ${\bf \Phi }^{\left( i\right) }$ is given
in terms of nested Lie brackets of body velocities at intermediate time
steps $t_{i-1}+c_{j}\Delta t$. Hence, the so determined ${\bf \Phi }^{\left(
i\right) }$ belongs to the smallest Lie subalgebra of ${\frak g}$ containing
the absolute body twists ${\bf V}$. This is a screw algebra \cite%
{GibsonHunt19901,GibsonHunt19902,Hunt1978,SeligBook} if $SE\left( 3\right) $
is used as c-space. The configuration increment $g\left( t_{i}\right) $
belongs to the corresponding subgroup of the c-space Lie group $G$,
sometimes called the {\it completion group}.

For the finite update step to respect the geometric constraints it is thus
necessary that the constrained absolute motion belongs to a subgroup of the
c-space Lie group. Since a configuration update step yields ${\bf \Phi }%
^{\left( i\right) }$ in the subalgebra generated by the twists, this
condition is also sufficient.

\begin{corollary}
\label{corollary1}The kinematic motion constraints of a body in a MBS are
satisfied by a configuration update step in terms of linear combinations of
velocity samples ${\bf V}^{\left( j\right) }$ if its motion is constrained
to a subgroup of its c-space Lie group.
\end{corollary}

This corollary applies to a general c-space Lie group. Understanding the
implications for a specific choice requires inspection of the possible
subgroups.

The direct product group $SO\left( 3\right) \times {\Bbb R}^{3}$ possesses
the following subgroups: The groups with elements $C=\left( {\bf I},{\bf r}%
\right) $, representing pure 1,2, or 3 dimensional translations; the groups
with elements $C=\left( {\bf R},{\bf 0}\right) $, representing pure 1 or 3
dimensional rotations such that the origin of body-fixed and spatial frame
coincide; and the groups with elements $C=\left( {\bf R},{\bf r}\right) $,
representing decoupled 1 or 3 dimensional rotation and 1,2, or 3 dimensional
translation. It is crucial to notice that the elements $C=\left( {\bf R},%
{\bf 0}\right) $ cannot represent rotations about arbitrary axes/points
since this would lead to non-zero ${\bf r}$. Clearly these groups have
limited practical relevance.

The 11 subgroups of $SE\left( 3\right) $, on the other hand, do have a
practical significance including modeling of lower pair joints that
correspond to six of the $SE\left( 3\right) $ subgroups \cite{SeligBook} and
are used in MBS modeling. Moreover, also the other $SE\left( 3\right) $
subgroups are used in MBS modeling (occasionally referred to as
'macro-joints'). The immediate consequence is that the constraints imposed
by a lower pair joint connecting a rigid body to the ground are exactly
satisfied by an $SE\left( 3\right) $ update step.

\begin{figure}[bh]
\centerline{a)%
\psfig{file=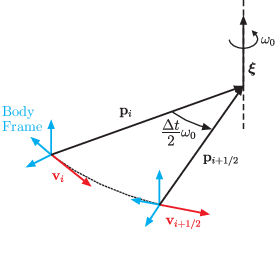,width=5.5cm}~~~~~~~~~~~~~b)%
\psfig{file=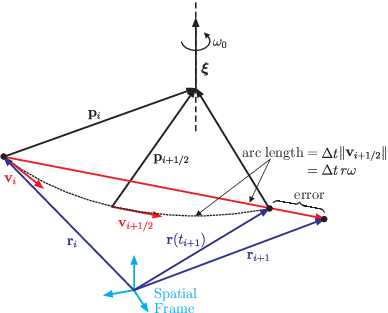,width=7.9cm}}
\caption{a) Frame rotating about a fixed axis with constant angular
velocity. b) Position error when using the $SO\left( 3\right) \times {\Bbb R}%
^{3}$ update.}
\label{figRotConst}
\end{figure}

To illustrate the effect of the two formulations of update steps, consider
the kinematic reconstruction equations $\dot{{\bf X}}={\bf dexp}_{-{\bf X}%
}^{-1}{\bf V}$ in the time interval $[t_{i},t_{i+1}]$, with ${\bf X}\left(
t_{i}\right) ={\bf 0}$ applied to a body rotating with constant angular
velocity $\omega _{0}$ about a constant rotation axis ${\bf \xi }$ at a
position ${\bf p}$ in figure \ref{figRotConst}a). The reference frame
performs a coupled rotation and translation. Using $SE\left( 3\right) $ as
c-space, the motion belongs to its subgroup $SO\left( 2\right) $. If $%
SO\left( 3\right) \times {\Bbb R}^{3}$ is used as c-space, the motion only
forms a submanifold. 

The body-fixed twist of the reference frame is ${\bf V}=\left( \omega {\bf %
\xi },\omega {\bf p}\times {\bf \xi }\right) $. For illustration purpose,
the simple explicit trapezoidal integration scheme is applied. In standard
vector space form, applied to solve $\dot{x}=f\left( t,x\right) $, the
latter reads 
\begin{eqnarray*}
k_{1} &=&f\left( t_{i},x_{i}\right) ,\ k_{2}=f\left( t_{i}+\frac{\Delta t}{2}%
,x_{i}+\frac{\Delta t}{2}k_{1}\right) \\
x_{i+1} &=&x_{i}+\Delta tk_{2}
\end{eqnarray*}%
where $\Delta t$ is the time step size. Applied to the kinematic equations
on the c-space Lie group, this is 
\begin{eqnarray*}
{\bf k}_{1} &=&{\bf dexp}_{-{\bf \xi }_{i}}^{-1}{\bf V}_{i}={\bf V}_{i},\ \
\ {\bf k}_{2}={\bf dexp}_{-\frac{\Delta t}{2}k_{1}}^{-1}{\bf V}_{i+1/2}={\bf %
dexp}_{-\frac{\Delta t}{2}{\bf V}_{i}}^{-1}{\bf V}_{i+1/2} \\
{\bf \xi }_{i+1} &=&{\bf \xi }_{i}+\Delta t{\bf k}_{2}=\Delta t{\bf dexp}_{-%
\frac{\Delta t}{2}{\bf V}_{i}}^{-1}{\bf V}_{i+1/2}\approx \Delta t{\bf V}%
_{i+1/2}-\frac{\Delta t^{2}}{4}\left[ {\bf V}_{i},{\bf V}_{i+1/2}\right]
+\ldots \\
{\bf C}_{i+1} &=&{\bf C}_{i}\exp {\bf \xi }_{i+1}
\end{eqnarray*}%
with ${\bf V}_{i}={\bf V}\left( t_{i}\right) $, ${\bf V}_{i+1/2}={\bf V}%
\left( t_{i}+\Delta t/2\right) $, and noting that ${\bf \xi }_{i}=0$. Denote
with ${\bf p}_{i}$ and ${\bf v}_{i}=\omega {\bf p}_{i}\times {\bf \xi }_{i}$
the position and linear velocity vector, respectively, at time $t_{i}$
(figure \ref{figRotConst}a). Since $\omega _{0}$ is constant, ${\bf p}_{i}=%
{\bf p}$ and ${\bf v}_{i}={\bf v}$ expressed in the body-fixed frame, and
thus ${\bf V}_{i}$, are constant. Hence the trapezoidal rule yields 
\begin{eqnarray}
{\bf \xi }_{i+1} &=&\Delta t{\bf V}_{i+1/2}=\Delta t{\bf V}_{i}=\Delta
t\left( \omega _{0}{\bf \xi },\omega _{0}{\bf p}\times {\bf \xi }\right) 
\nonumber \\
{\bf C}_{i+1} &=&{\bf C}_{i}\exp \left( \Delta t\left( \omega _{0}{\bf \xi }%
,\omega _{0}{\bf p}\times {\bf \xi }\right) \right) .  \label{TrapSol}
\end{eqnarray}%
Apparently, $\left( {\bf \xi },{\bf p}\times {\bf \xi }\right) $ are the Pl%
\"{u}cker coordinates of the screw motion of the body, when rotating about
the axis ${\bf \xi }$ at a distance ${\bf r}$, so that (\ref{TrapSol}) is in
fact the exact solution.\newline
In contrast, if the mixed velocity (\ref{hybridvel}) and the dexp mapping (%
\ref{dexpSO3xR3}) on $SO\left( 3\right) \times {\Bbb R}^{3}$ is used, then
the orientation update is exact but the position (vector space) update ${\bf %
r}_{i+1}={\bf r}_{i}+\Delta t{\bf v}_{i+1/2}$ is not correct as apparent
from figure \ref{figRotConst}b). This is due to the independent treatment of
orientations and translations. Hence $SE\left( 3\right) $ is the proper
c-space for this example.

This example confirms corollary \ref{corollary1}. The perfect constraint
satisfaction is possible (for any integration scheme) because the body is
connected to the ground by a lower pair joint. Moreover, the corollary
regards the motion of constrained bodies rather than their relative motions,
i.e. the joint motions. Since the joint constraints are the actual
constraints included in the motion equations (\ref{BH}), a statement
regarding the joint constraints is in order.

\begin{lemma}
\label{lemma1}Consider two bodies constrained to a c-space subgroup $G_{1}$
and $G_{2}$, respectively. If their relative motions form a subgroup $H$,
i.e. $G_{2}=G_{1}\cdot H$, then the relative motion constraints as well as
the constraints on the motion of the two bodies are satisfied by any
numerical configuration update step in terms of linear combinations of
velocity samples.
\end{lemma}

\paragraph{Proof.}

The proof follows immediately from the assumption that the group $G_{2}$ is
the product of the groups $G_{1}$ and $H$, and thus $H=G_{2}/G_{1}$. That
is, $g_{1}g_{2}\in H$ for any $g_{1}\in G_{1},g_{2}\in G_{2}$. Hence, since
any update of the bodies' motion belongs to the respective subgroup $G_{1}$
and $G_{2}$, the relative motion will be in $H$. \mbox{\ \rule{.1in}{.1in}}

The situation of this lemma is a rather special case since the product of
two groups is not necessarily a group itself. A special case is the rigid
body connected to the ground by a lower pair joint. In this case, the
corollary \ref{corollary1} and lemma \ref{lemma1} are automatically
satisfied with $G_{1}=\{I\}$ and $G_{2}=H\subset SE\left( 3\right) $, if $%
SE\left( 3\right) $ is used as c-space. This is not so if $SO\left( 3\right)
\times {\Bbb R}^{3}$ is used.

As outlined above, lower kinematic pairs form $SE\left( 3\right) $
subgroups, while $SO\left( 3\right) \times {\Bbb R}^{3}$ subgroups have
almost no practical relevance. The only practically relevant exception where
the rigid body motions form a subgroup of $SO\left( 3\right) \times {\Bbb R}%
^{3}$ is the case of an unconstrained free body with its reference frame
attached to the center of mass (COM) (section \ref{secFreeBody}).

When lemma \ref{lemma1} does not hold, the joint constraints will not be
satisfied exactly. In this case, since both formulations use a first-order
relation consistent with the respective kinematic model to estimate the
update, the order of accuracy obtained with the $SE\left( 3\right) $ and
with the $SO\left( 3\right) \times {\Bbb R}^{3}$ update is the same, and the
constraint violations are determined by the accuracy of the integration
scheme.

\begin{corollary}
\label{corollary3}%
For general MBS, the $SE\left( 3\right) $ and the $SO\left( 3\right) \times 
{\Bbb R}^{3}$ update achieve the same order of accuracy. The $SE\left(
3\right) $ update achieves exact satisfaction of geometric constraints for
those joints for which lemma \ref{lemma1} holds.
\end{corollary}

\begin{remark}
When the condition of corollary \ref{corollary1} holds, the motion
constraints of bodies are perfectly satisfied. The joint constraints are
satisfied if lemma \ref{lemma1} holds. It must be emphasized that this does
not imply any increase of the order of accuracy of the overall numerical
solution. But satisfaction of geometric constraints is vital for the
numerical integration. Even more, numerical integration of constrained MBS
(using the index 1 formulation (\ref{index1})) consists of two parts: 1) the
numerical integration of the ODE (\ref{ODEvec1}) and 2) the stabilization of
constraints. The latter is computationally expensive and any means to avoid
constraint violations is beneficial. Note that the conclusion applies to the
'standard' vector space formulation (\ref{ODEvec1}) in 'local' coordinates
as well as to the Lie group integration scheme for the formulation (\ref{BH3}%
).
\end{remark}

\section{Examples%
\label{secExamples}%
\label{secExamples}}

All examples consist of rigid bodies either unconstrained or subjected to
geometric constraints imposed by lower pair joints, and the dynamics
equations (\ref{BH2}a), respectively (\ref{BH3}a), are the Newton-Euler
equations. The explicit form of the latter depends on the velocity
representation (body-fixed twists, mixed velocities), thus the c-space Lie
group, and the used reference frame.

Results are reported when the Lie group formulation (\ref{ODE2}) is
numerically integrated using time step sizes $\Delta t=10^{-2}$s$,10^{-3}$s$%
,10^{-4}$s. The dynamic equations (\ref{ODE2}) are integrated with the
Runge-Kutta 4 (RK4) scheme, and the kinematic equations are solved with the
MK method (\ref{MK1}) based on the RK4 scheme, using the $SE\left( 3\right) $
and $SO\left( 3\right) \times {\Bbb R}^{3}$ update, respectively.

It should be recalled that the primary interest here is the accuracy of
constraint satisfaction rather than that of the solution ${\bf q}\left(
t\right) $. The latter depends indeed on the numerical integration scheme.

\subsection{Unconstrained Rigid Body -- Reference Frame at COM%
\label{secFreeBody}%
}

An unconstrained rigid body is considered assuming that there are no applied
or gravitational forces. Although no kinematic constraints are present, the
momentum conservation imposes motion invariants. The body-fixed reference
frame is located at its COM. The configuration of the body, i.e. of the
reference frame, is represented by $C=\left( {\bf R},{\bf r}^{\text{s}%
}\right) $ with ${\bf r}^{\text{s}}$ being the position vector of the COM
expressed in the space-fixed inertial frame and ${\bf R}$ the rotation
matrix transforming coordinates from body-fixed reference frame to inertia
frame.

In this example the rigid body is a homogenous box with side lengths $%
0.8\times 0.4\times 0.1$~m made of aluminium. Its mass is $m=86.4$~kg and
its inertia tensor w.r.t. the COM is ${\bf \Theta }_{0}={\rm diag~}\left(
1.224,4.68{,}5.76\right) \,$kg\thinspace m$^{2}$.

\subsubsection{Body-Fixed Newton-Euler Equations on $SE(3)$}

The dynamics of the body is governed by the Newton-Euler equations.
Considering $C\in SE\left( 3\right) $, the body velocity represented in the
body-fixed COM frame is a proper twist ${\bf V}=({\bf \omega },{\bf v})\in
se\left( 3\right) $, and the Newton-Euler equations attain the consistent
form 
\begin{equation}
{\bf J}_{0}\dot{{\bf V}}-{\bf ad}_{{\bf V}}^{\ast }{\bf J}_{0}{\bf V}={\bf W}
\label{NE1}
\end{equation}%
where ${\bf W}=({\bf M},{\bf F})\in se^{\ast }\left( 3\right) $ is the
applied wrench expressed in the COM frame, and 
\begin{equation}
{\bf J}_{0}=\left( 
\begin{array}{cc}
{\bf \Theta }_{0} & {\bf 0} \\ 
{\bf 0} & m{\bf I}%
\end{array}%
\right)
\end{equation}%
is the inertia matrix w.r.t. to the COM frame. The matrix ${\bf ad}_{{\bf V}%
}^{\ast }={\bf ad}_{{\bf V}}^{T}$ is given by (\ref{adse3}). The system (\ref%
{NE1}) are the Euler-Poincare equations on $SE\left( 3\right) $.

\subsubsection{Newton-Euler Equations on $SO(3)\times {\Bbb R}^{3}$}

In the mixed velocity representation, the angular velocity ${\bf \omega }$
is expressed in the body-fixed COM frame whereas the linear velocity of the
COM is expressed in the space-fixed inertial frame, denoted ${\bf v}^{\text{s%
}}=\dot{{\bf r}}$. The corresponding Newton-Euler equations, written w.r.t.
the COM frame, are

\begin{equation}
\left( 
\begin{array}{cc}
{\bf \Theta }_{O} & {\bf 0} \\ 
{\bf 0} & m{\bf I}%
\end{array}%
\right) \left( 
\begin{array}{c}
\dot{{\bf \omega }} \\ 
\dot{{\bf v}}^{\text{s}}%
\end{array}%
\right) +\left( 
\begin{array}{c}
\widehat{{\bf \omega }}{\bf \Theta }_{O}{\bf \omega } \\ 
{\bf 0}%
\end{array}%
\right) =\left( 
\begin{array}{c}
{\bf M} \\ 
{\bf F}^{\text{s}}%
\end{array}%
\right)  \label{NE2}
\end{equation}%
where ${\bf F}^{\text{s}}$ is the applied force expressed in the inertial
frame. The equations (\ref{NE2}) can be expressed as the are the
Euler-Poincare equations on $SO\left( 3\right) \times {\Bbb R}^{3}$: ${\bf J}%
_{0}\dot{{\bf V}}-{\bf ad}_{{\bf V}}^{\ast }{\bf J}_{0}{\bf V}={\bf W}$ with 
${\bf W}^{\text{m}}=({\bf M},{\bf F}^{\text{s}})$ and (\ref{adso3xr3}).

\subsubsection{Numerical Results}

The two formulations are integrated with the MK method using the time step
sizes $\Delta t=10^{-2}$s$,10^{-3}$s$,10^{-4}$s. The initial configuration
of the body is set to $C_{0}=\left( {\bf I},{\bf 0}\right) $.

\subparagraph{Spatial Rotation}

Setting the initial angular velocity to ${\bf \omega }_{0}=(10\pi ,2\pi
,0)\, $rad/s and ${\bf v}_{0}=0$, the unconstrained unforced body should
rotate about its COM. That is, the vector ${\bf r}^{\text{s}}$ should remain
zero. The motion equations are integrated for 10\ s.

\begin{figure}[h]
\centerline{\ $
\begin{array}{c@{\hspace{6ex}}c}
a)\psfig{file=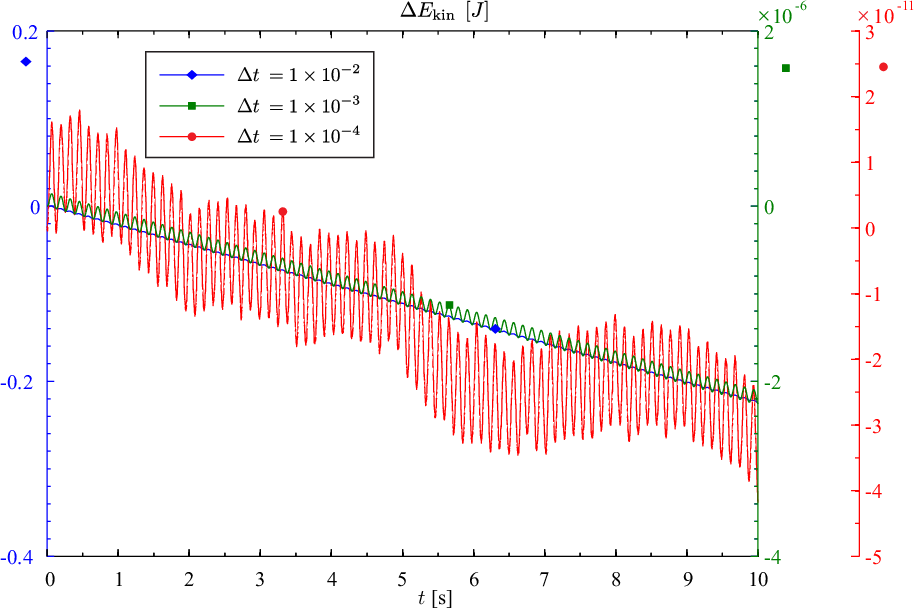,width=8cm} & b)\psfig{file=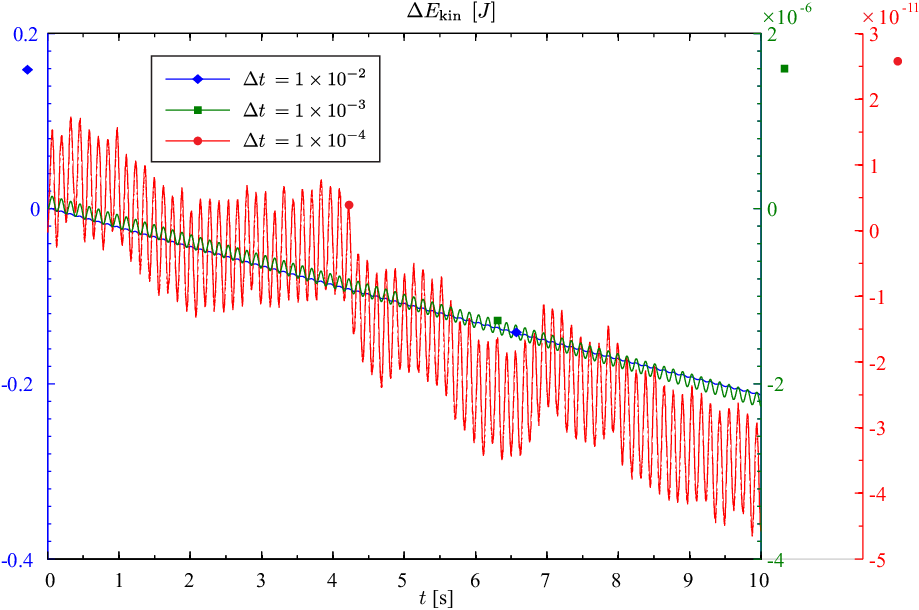,width=8cm}
\end{array}
$ }
\caption{Drift of kinetic energy when integrating a) the $SE(3)$, and b) the 
$SO(3)\times {\Bbb R}^{3}$ formulation.}
\label{EKinFreeBody_COM_RotSkewAxis}
\end{figure}

Both formulations exactly preserve the COM position ${\bf r}^{\text{s}}={\bf %
0}$. It is clear that both formulations perform identical since for a pure
rotation about the COM the $SE\left( 3\right) $ reduces to the $SO\left(
3\right) \times {\Bbb R}^{3}$ formulation. This is visible from the drift of
the kinetic energy, shown in figure \ref{EKinFreeBody_COM_RotSkewAxis},
which should be preserved. 

\subparagraph{Rotation about fixed axis plus linear translation%
\label{secRot+LinVel}%
\label{secRot+LinVel}}

The initial velocities are set to ${\bf \omega }_{0}=(0,0,2\,\pi )$%
\thinspace rad/s and ${\bf v}_{0}=\left( 10,0,0\right) $\thinspace m/s.
Since the body is unconstrained and no gravity forces are present, it should
perform a translation of its COM along the global $x$-axis together with a
rotation about the global $z$-axis. The analytic solution for the position
vector is ${\bf r}^{\text{s}}\left( t\right) =\left( t,0,0\right) $~m/s and
for the rotation angle $\varphi \left( t\right) =2\pi t$. Figure \ref%
{ErrorFreeBody_COM_Rot3AxisLinVel} show the position errors for the
numerical solutions. For the translation the $SO\left( 3\right) \times {\Bbb %
R}^{3}$ update yields very good reconstruction for any step size while the
accuracy achieved with the $SE\left( 3\right) $ update depends on the step
size. Clearly visible is the 4th order convergence of the MK method. The
rotation update performs equally for both variants, as shown in figure \ref%
{AngleErrorFreeBody_COM_Rot3AxisLinVel}, where the relative rotation angle $%
\varepsilon _{r}:=||\log ({\bf R}(t)^{T}{\bf R}_{{\rm num}})||$ of the
analytic solution ${\bf R}(t)$ and the respective numerical solution ${\bf R}%
_{{\rm num}}$ are shown. Notice that the accuracy is within the computation
precision. The figures show essentially the error amplified by the log
mapping.

Due to the momentum conservation, the body performs translation along x-axis
plus rotation about z-axis. This is kinematically equivalent to a rolling
disk of radius $r=\left\Vert {\bf v}_{0}\right\Vert /\left\Vert {\bf \omega }%
_{0}\right\Vert $. This means that the body can be considered being
subjected to a kinematic rolling constraint, or in other words being
connected to the ground by a higher kinematic pair whose motions do not form
a $SE\left( 3\right) $ subgroup. It is a subgroup of $SO\left( 3\right)
\times {\Bbb R}^{3}$, however. In fact, the decoupled rotation and
translation is the ideal situation for the $SO\left( 3\right) \times {\Bbb R}%
^{3}$ update, and according to the corollary 1 it should satisfy the motion
constraints. 
\begin{figure}[h]
\centerline{\ $
\begin{array}{c@{\hspace{6ex}}c}
a)\psfig{file=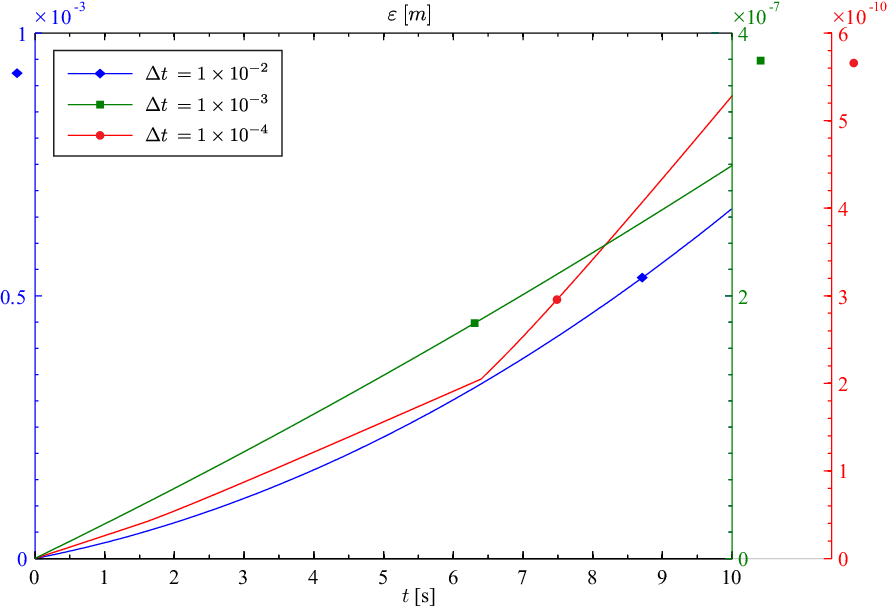,width=8cm}
& b)\psfig{file=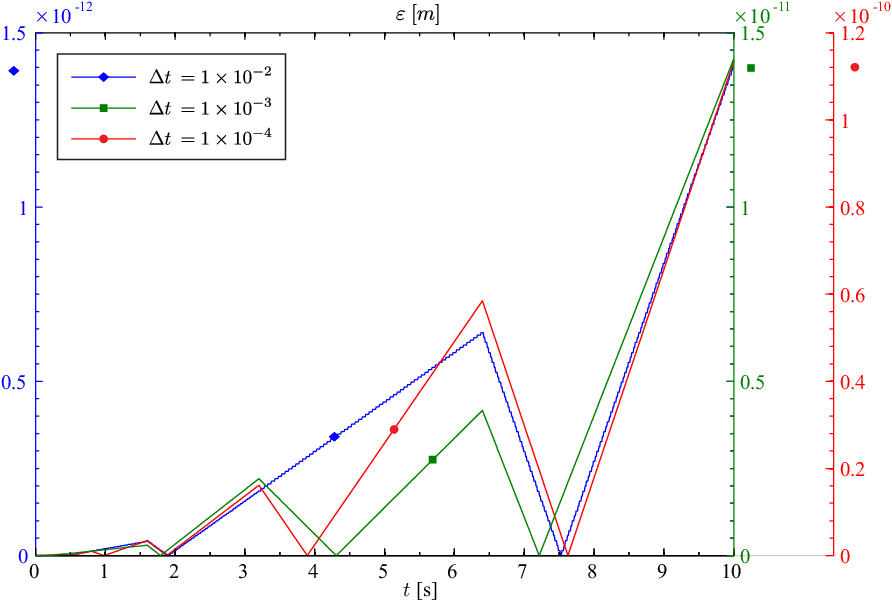,width=8cm}
\end{array}
$ }
\caption{Drift of COM from analytic solution when integrating a) the $SE(3)$%
, and b) the $SO(3)\times {\Bbb R}^{3}$ formulation.}
\label{ErrorFreeBody_COM_Rot3AxisLinVel}
\end{figure}
\begin{figure}[h]
\centerline{\ $
\begin{array}{c@{\hspace{6ex}}c}
a)\psfig{file=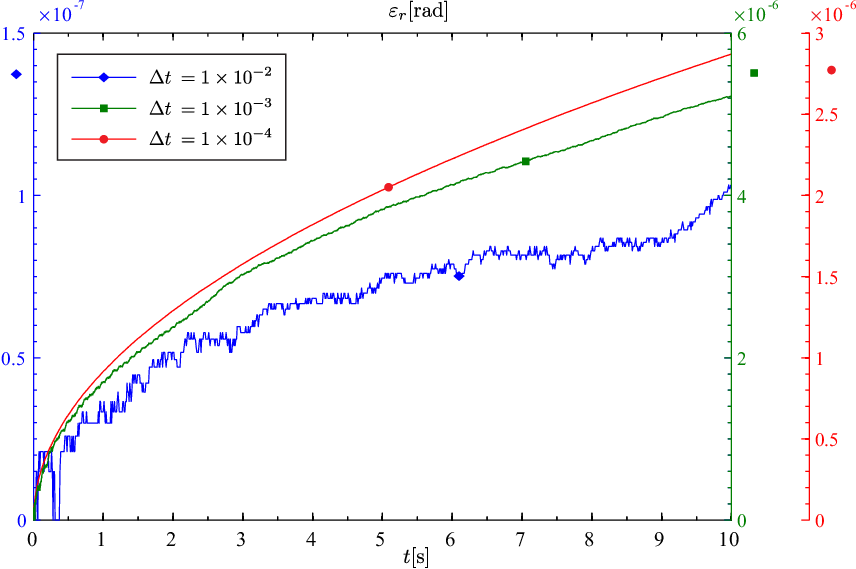,width=8cm}
& b)\psfig{file=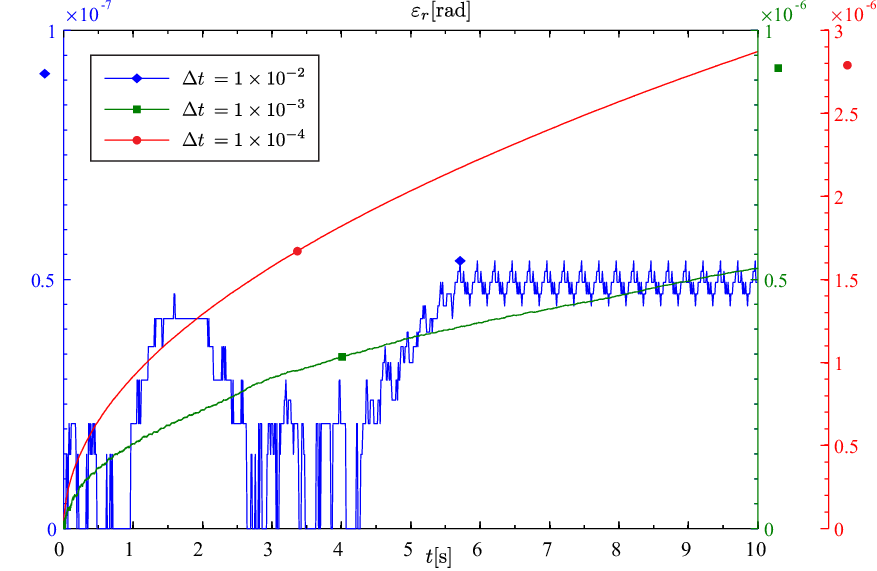,width=8cm}
\end{array}
$ }
\caption{Drift $\protect\varepsilon _{r}$ of the orientation from analytic
solution when integrating a) the $SE(3)$, and b) the $SO(3)\times {\Bbb R}%
^{3}$ formulation.}
\label{AngleErrorFreeBody_COM_Rot3AxisLinVel}
\end{figure}


\subsection{Unconstrained Rigid Body -- Arbitrary Body-Fixed Reference Frame}

Besides geometric constraints, another cause of coupling of rotations and
translations is the use of body-fixed reference frames not located at the
COM. Consider again the rigid body from the previous example, but now with a
reference frame parallel translated to an arbitrary body-fixed point $P$,
i.e. now $C=\left( {\bf R},{\bf r}^{\text{s}}\right) $ represents the
configuration of this reference frame. Denote with ${\bf r}_{0}=\left(
0.4,0,0\right) $\thinspace m the COM location expressed in the body-fixed
reference frame. Then, transforming the inertia tensor ${\bf \Theta }_{0}$
to the reference frame yields the inertia tensor ${\bf \Theta }_{P}={\rm %
diag~}\left( 1.224,18.504{,}19.584\right) \,$kg\thinspace m$^{2}$. The
motion of the body is indeed the same regardless of the used reference
frame, but the motion of the introduced reference frame varies with its
location. This has implications for the numerical configuration update.

\subsubsection{Body-Fixed Newton-Euler Equations on $SE(3)$}

The Newton-Euler equations (\ref{NE1}) are coordinate invariant, and using
body-fixed velocities the dynamics of the body is still governed by (\ref%
{NE1}), but now with 
\begin{equation}
{\bf J}_{P}=\left( 
\begin{array}{cc}
{\bf \Theta }_{P} & -m\widehat{{\bf r}}_{0}^{T} \\ 
-m\widehat{{\bf r}}_{0} & m{\bf I}%
\end{array}%
\right)
\end{equation}%
where ${\bf M},{\bf F}$ are expressed in the reference frame. The vector $%
{\bf V}=\left( {\bf \omega },{\bf v}\right) $ is the twist of the reference
frame at $P$.

\subsubsection{Newton-Euler Equations on $SO(3)\times {\Bbb R}^{3}$}

The Newton-Euler equations are no longer decoupled if the reference frame is
not at the COM, and (\ref{NE2}) must be replaced by 
\begin{equation}
\left( 
\begin{array}{cc}
{\bf \Theta }_{P} & -m({\bf R}\widehat{{\bf r}}_{0})^{T} \\ 
-m{\bf R}\widehat{{\bf r}}_{0} & m{\bf I}%
\end{array}%
\right) \left( 
\begin{array}{c}
\dot{{\bf \omega }} \\ 
\dot{{\bf v}}^{\text{s}}%
\end{array}%
\right) +\left( 
\begin{array}{c}
\widehat{{\bf \omega }}{\bf \Theta }_{P}{\bf \omega } \\ 
{\bf R}\widehat{{\bf \omega }}\widehat{{\bf \omega }}{\bf r}_{0}%
\end{array}%
\right) =\left( 
\begin{array}{c}
{\bf M} \\ 
{\bf F}^{\text{s}}%
\end{array}%
\right) .
\end{equation}

${\bf V}^{\text{m}}=\left( {\bf \omega },{\bf v}^{\text{s}}\right) $ is the
mixed velocity of the reference frame at $P$.

\subsubsection{Numerical Results}

As above the initial configuration of the body is set to $C_{0}=\left({\bf I}%
,{\bf 0}\right) $, and the dynamic equations are integrated with three
different time step sizes.

\subparagraph{Spatial Rotation}

The situation is the same as in the previous section, where the initial
angular velocity is ${\bf \omega }_{0}=(10\pi ,2\pi ,0)$\thinspace rad/s and
the COM is at rest. When the body is rotating about its COM the reference
frame exhibits the initial linear velocity ${\bf v}_{0}={\bf r}_{0}\times 
{\bf \omega }_{0}$. Again, the unconstrained body without gravity forces
should perform rotations about its COM. That is, the vector $\Delta {\bf r}=%
{\bf r}-{\bf r}_{0}$ should remain zero. Figure \ref%
{FigErrorFreeBody_RotSkewAxis} shows the drift $\varepsilon :=||\Delta {\bf r%
}||$ of the COM position. The position accuracy achieved by the $SE(3)$
update is within the computation accuracy of $10^{-15}$ for all step sizes,
while that achieved with the $SO\left( 3\right) \times {\Bbb R}^{3}$
formulation depends on the integration step size. Both formulation achieve
the same rotation update and exhibit identical drift of the kinetic energy $%
T_{0}=696.399$\thinspace J (not shown here). Note that the $SE\left(
3\right) $ update satisfies the constraints for arbitrary time step sizes. 
\begin{figure}[h]
\centerline{\ $
\begin{array}{c@{\hspace{4ex}}c}
a)\psfig{file=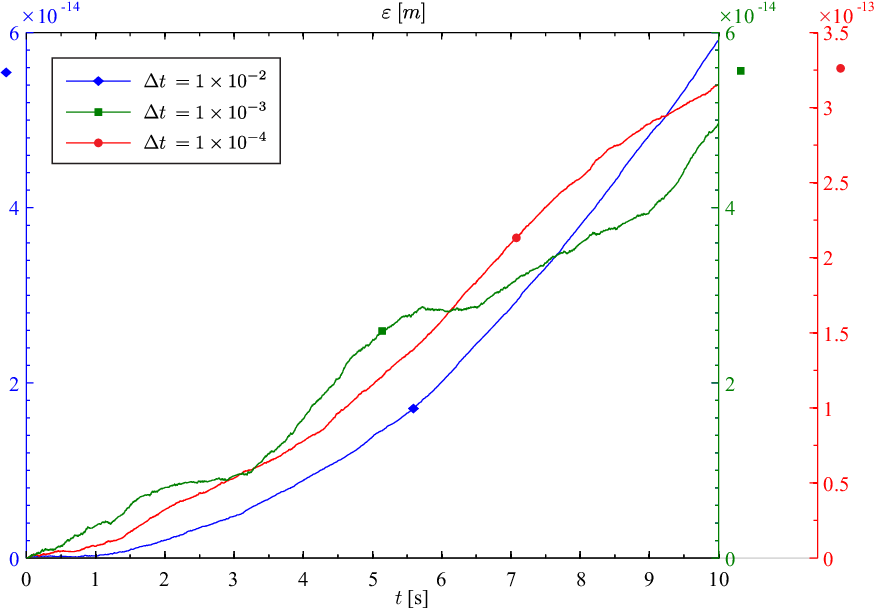,width=8cm} & b)\psfig{file=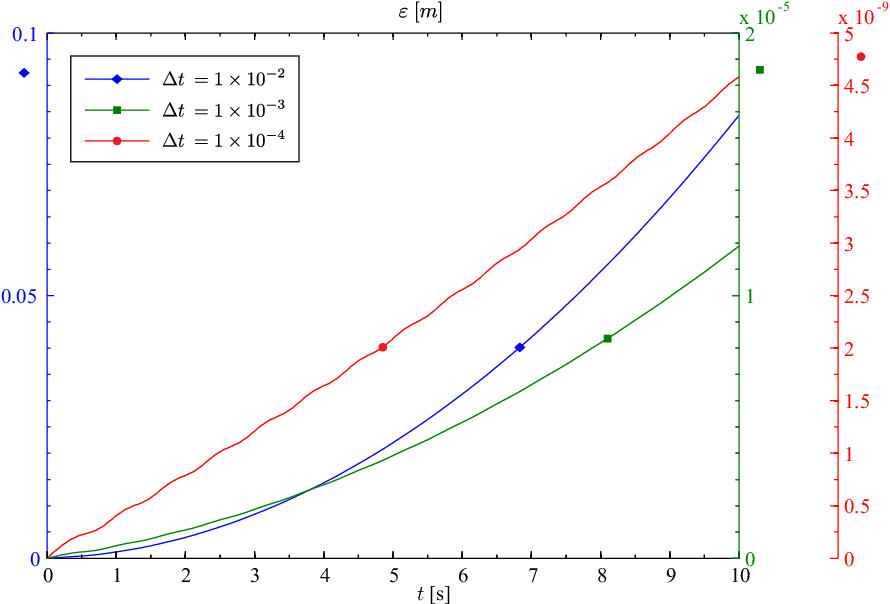,width=8cm}
\end{array}
$ }
\caption{Drift $\protect\varepsilon $ of the COM when integrating a) the $%
SE(3)$, and b) the $SO(3)\times {\Bbb R}^{3}$ formulation.}
\label{FigErrorFreeBody_RotSkewAxis}
\end{figure}

\subparagraph{Rotation about fixed axis plus linear translation}

If, in addition to the initial angular velocity ${\bf \omega }_{0}=(10\pi
,2\pi ,0)$\thinspace rad/s, the body is again given a constant initial
linear velocity of 10\thinspace m/s, the initial velocity of the reference
frame is ${\bf v}_{0}={\bf r}_{0}\times {\bf \omega }_{0}+(10,0,0)$%
\thinspace m/s. The body performs the same spatial rotation as above
together with the linear motion of its COM, i.e. a screw motion of the
reference frame.

\begin{figure}[h]
\centerline{\ $
\begin{array}{c@{\hspace{4ex}}c}
a)\psfig{file=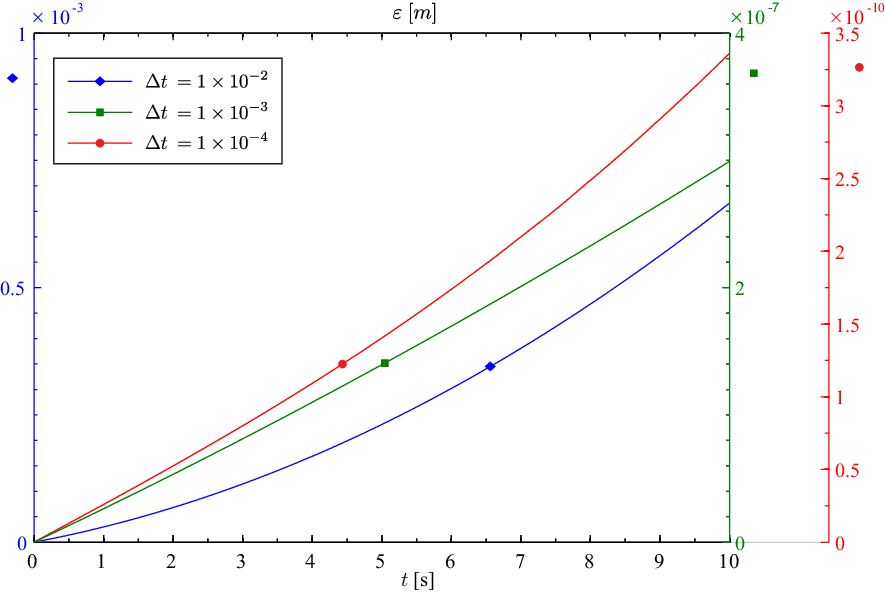,width=8cm} & b)\psfig{file=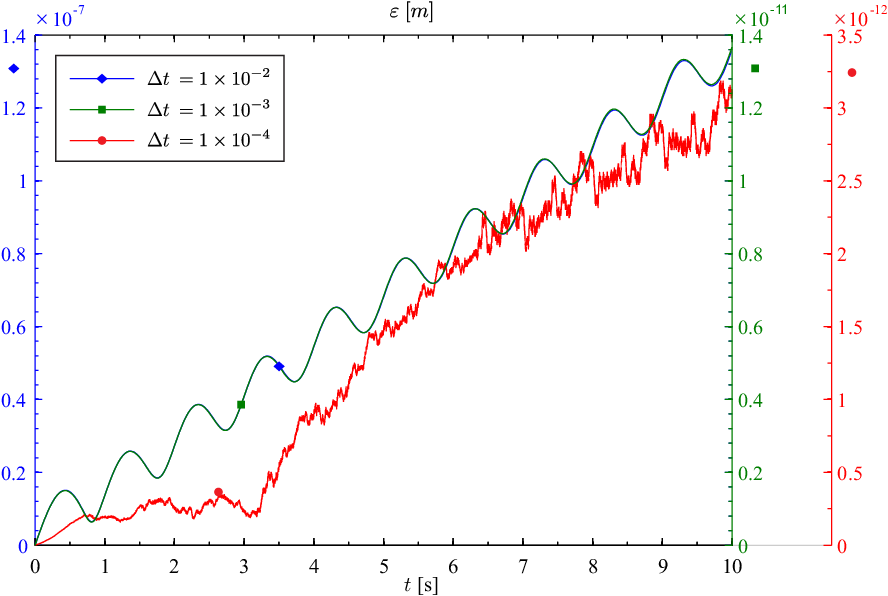,width=8cm}
\end{array}
$ }
\caption{Drift $\protect\varepsilon $ of the COM from the analytic solution
when integrating a) the $SE(3)$, and b) the $SO(3)\times {\Bbb R}^{3}$
formulation.}
\label{FigErrorFreeBody_Rot3AxisLinVel}
\end{figure}

Figure \ref{FigErrorFreeBody_Rot3AxisLinVel} reveals that the accuracy of
both formulations depends on the step size. The $SO(3)\times {\Bbb R}^{3}$
update yields the best accuracy as for the above model with the reference
frame at the COM. As for the above model, the lower accuracy of the $%
SE\left( 3\right) $ update is attributed to the fact that the motion does
not belong to a subgroup. The result is to be expected since the $SE\left(
3\right) $ update scheme is frame invariant. Also the $SO(3)\times {\Bbb R}%
^{3}$ update cannot exactly estimate the finite motion since with the chosen
reference frame this does not belong to a subgroup.

\subsection{Heavy Top}

\begin{figure}[h]
\centerline{\ \psfig{file=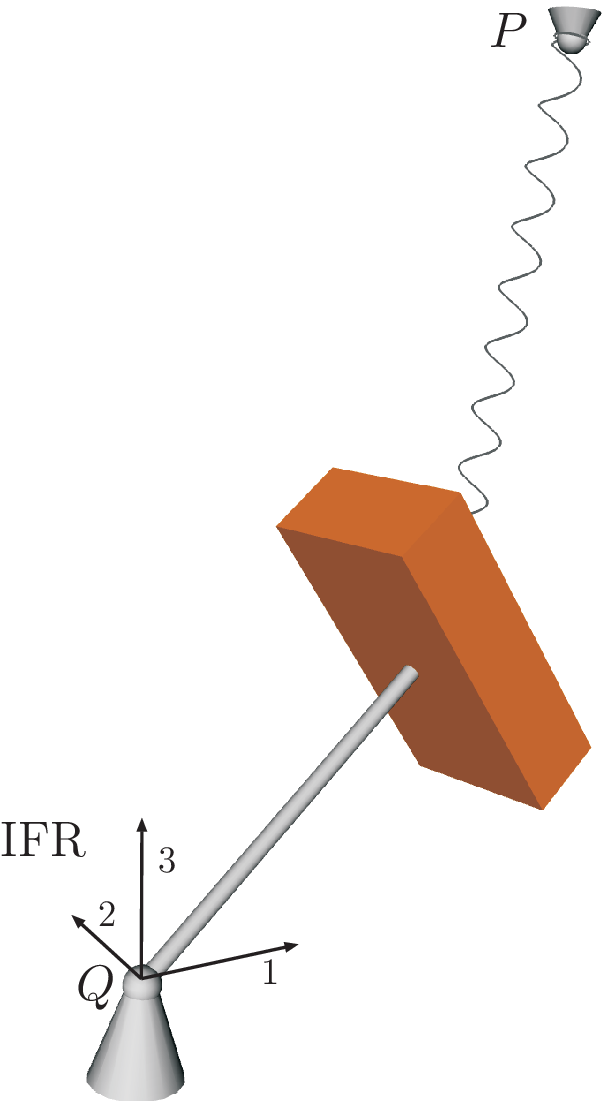,width=4.0cm} 
}
\caption{Model of a heavy top consisting of a rigid body, constrained by a
massless rod to rotate about a fixed point, suspended with a spring attached
to its COM.}
\label{figHeavyTopSpring}
\end{figure}

A heavy top is modeled as a rigid body pivoted to the ground at the point $Q$
as shown in figure \ref{figHeavyTopSpring}. The top is suspended by a spring
attached to its COM (the geometric center of the box) and a space-fixed
point $P$. The body-fixed reference frame is located at the COM, which, in
the reference configuration, has the same orientation as the shown inertial
frame. The pivot point is the center of rotation of the top moving in
gravity field. 

The spring is modeled by a linear force law with stiffness $c=10$\thinspace
N/mm. Then, in the spatial representation, the vector of applied forces in (%
\ref{TopEOMspatial}) is ${\bf F}^{\text{s}}=c\left( {\bf p}_{0}^{\text{s}}-%
{\bf r}^{\text{s}}\right) +m{\bf g}^{\text{s}}$, with the gravity vector $%
{\bf g}^{\text{s}}=\left( 0,0,-g\right) $, $g=9.81~$m/s$^{2}$. Here, ${\bf p}%
_{0}^{\text{s}}=\left( 1,0,0.5\right) $~m is the space-fixed position vector
of the suspension point $P$ of the spring. Further, in (\ref{TopEOMspatial}) 
${\bf M}={\bf 0}$, since no torques act on the COM. The force vector in the
body-fixed representation (\ref{TopEOMbody}) is ${\bf F}={\bf R}^{T}{\bf F}^{%
\text{s}}$.

The body is a solid box with side lengths $0.1\times 0.2\times 0.4$%
\thinspace m made of aluminium. The body has a mass of $m=21.6$\thinspace kg
and its inertia tensor w.r.t. the COM is ${\bf \Theta }_{0}={\rm diag}%
\,\left( 0.36,{0.306,0.09}\right) \,$kg\thinspace m$^{2}$. The top is
constrained to rotate about a fixed point. Denote with ${\bf r}%
_{0}=(-0.5,0,0)$\thinspace m the position vector of the pivot point measured
in the body-fixed reference frame. The configuration of the reference frame
is represented by $C=\left( {\bf R},{\bf r}^{\text{s}}\right) $, with
rotation matrix ${\bf R}$ and ${\bf r}^{\text{s}}$ denoting the position of
the COM expressed in the spatial inertial frame (IFR).

\subsubsection{Motion Equations in Body-Fixed Representation}

The body-fixed twist is denoted ${\bf V}=\left( {\bf \omega },{\bf v}\right) 
$. The geometric constraints imposed by the spherical joint (pivot) are 
\begin{equation}
h\left( C\right) ={\bf r}^{\text{s}}+{\bf Rr}_{0}={\bf 0.}
\label{geomContsTop}
\end{equation}%
The time differentiation, and assuming (\ref{geomContsTop}), yields the
velocity constraints (\ref{VelConstr}) 
\begin{equation}
\left( 
\begin{array}{cc}
\widehat{{\bf r}}_{0} & -{\bf I}%
\end{array}%
\right) \left( 
\begin{array}{c}
{\bf \omega } \\ 
{\bf v}%
\end{array}%
\right) ={\bf JV}={\bf 0}
\end{equation}%
and the acceleration constraints 
\begin{equation}
\left( 
\begin{array}{cc}
\widehat{{\bf r}}_{0} & -{\bf I}%
\end{array}%
\right) \left( 
\begin{array}{c}
\dot{{\bf \omega }} \\ 
\dot{{\bf v}}%
\end{array}%
\right) =\widehat{{\bf \omega }}\widehat{{\bf \omega }}{\bf r}_{0}+\widehat{%
{\bf \omega }}{\bf v.}  \label{accConstrTop}
\end{equation}%
The body-fixed Newton-Euler equations w.r.t. to the COM combined with (\ref%
{accConstrTop}) yield the overall index 1 DAE system

\begin{equation}
\left( 
\begin{array}{ccc}
{\bf \Theta }_{0} & {\bf 0} & \widehat{{\bf r}}_{0}^{T} \\ 
{\bf 0} & m{\bf I} & -{\bf I} \\ 
\widehat{{\bf r}}_{0} & -{\bf I} & {\bf 0}%
\end{array}%
\right) \left( 
\begin{array}{c}
\dot{{\bf \omega }} \\ 
\dot{{\bf v}} \\ 
{\bf \lambda }%
\end{array}%
\right) =\left( 
\begin{array}{c}
\widehat{{\bf \omega }}{\bf \Theta }_{0}{\bf \omega } \\ 
{\bf F}-m\widehat{{\bf \omega }}{\bf v} \\ 
\widehat{{\bf \omega }}\widehat{{\bf \omega }}{\bf r}_{0}+\widehat{{\bf %
\omega }}{\bf v}%
\end{array}%
\right) .  \label{TopEOMbody}
\end{equation}%
The vector ${\bf F}$ is the external force acting upon the COM represented
in the body-fixed frame.

\subsubsection{Motion Equations in Mixed Velocity Representation}

In the mixed velocity representation ${\bf V}^{\text{m}}=\left( {\bf \omega }%
,{\bf v}^{\text{s}}\right) $, the geometric constraints (\ref{geomContsTop})
gives rise to the following velocity and acceleration constraints,
respectively, 
\begin{eqnarray}
\left( 
\begin{array}{cc}
{\bf R}\widehat{{\bf r}}_{0} & -{\bf I}%
\end{array}%
\right) \left( 
\begin{array}{c}
{\bf \omega } \\ 
{\bf v}^{\text{s}}%
\end{array}%
\right) &=&{\bf 0} \\
\left( 
\begin{array}{cc}
{\bf R}\widehat{{\bf r}}_{0} & -{\bf I}%
\end{array}%
\right) \left( 
\begin{array}{c}
\dot{{\bf \omega }} \\ 
\dot{{\bf v}}^{\text{s}}%
\end{array}%
\right) &=&{\bf R}\widehat{{\bf \omega }}\widehat{{\bf \omega }}{\bf r}_{0}.
\end{eqnarray}%
This, together with the Newton-Euler equations w.r.t. to the COM in the
mixed representation, yields 
\begin{equation}
\left( 
\begin{array}{ccc}
{\bf \Theta } & {\bf 0} & -\widehat{{\bf r}}_{0}{\bf R}^{T} \\ 
{\bf 0} & m{\bf I} & -{\bf I} \\ 
{\bf R}\widehat{{\bf r}}_{0} & -{\bf I} & {\bf 0}%
\end{array}%
\right) \left( 
\begin{array}{c}
\dot{{\bf \omega }} \\ 
\dot{{\bf v}}^{\text{s}} \\ 
{\bf \lambda }%
\end{array}%
\right) =\left( 
\begin{array}{c}
-\widehat{{\bf \omega }}{\bf \Theta }_{0}{\bf \omega } \\ 
{\bf F}^{\text{s}} \\ 
{\bf R}\widehat{{\bf \omega }}\widehat{{\bf \omega }}{\bf r}_{0}%
\end{array}%
\right)  \label{TopEOMspatial}
\end{equation}%
where ${\bf F}^{\text{s}}$ is the external force acting upon the COM
represented in the spatial inertial frame.

\subsubsection{Numerical Results}

The systems (\ref{TopEOMspatial}) and (\ref{TopEOMbody}) were integrated for
8\thinspace s starting from the initial configuration $C_{0}=\left( {\bf I},%
{\bf r}_{0}\right) $ with the initial velocities ${\bf \omega }%
_{0}=(0,0,0.5) $\thinspace rad/s and ${\bf v}_{0}={\bf r}_{0}\times {\bf %
\omega }_{0}$.

\begin{figure}[h]
\centerline{\ $
\begin{array}{c@{\hspace{4ex}}c}
a)\psfig{file=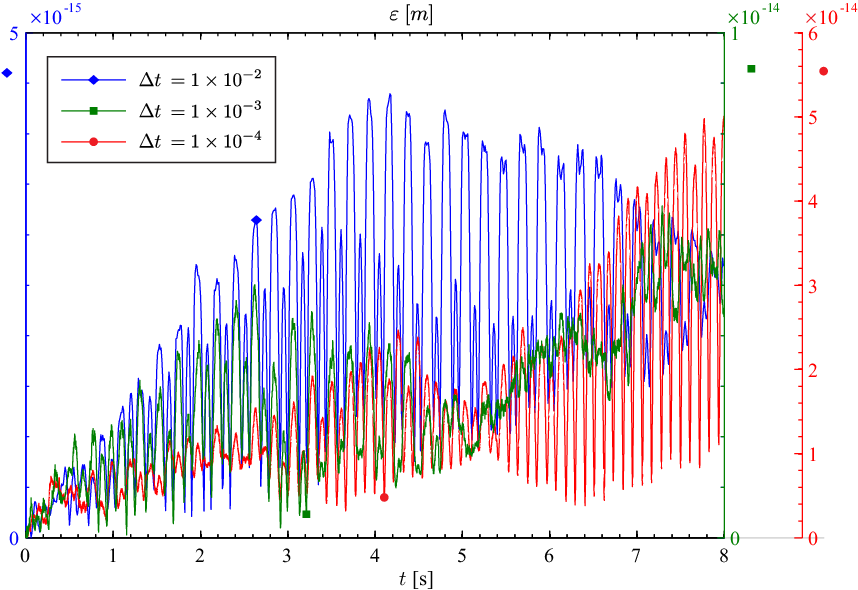,width=8cm} & b)\psfig{file=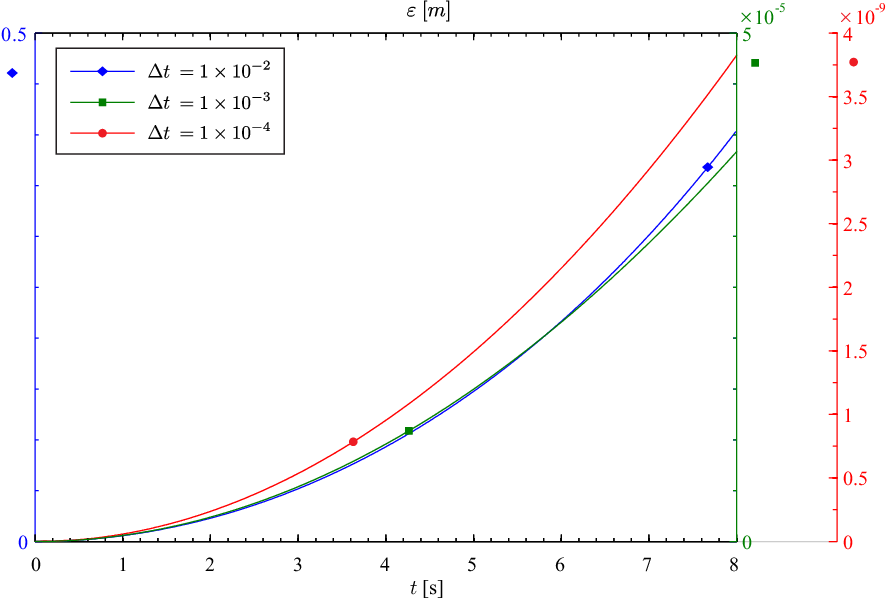,width=8cm}
\end{array}
$ }
\caption{Error $\protect\varepsilon $ in the position constraints of the top
when integrating a) the $SE(3)$, and b) the $SO(3)\times {\Bbb R}^{3}$
formulation.}
\label{FigErrorHeavyTop_GravitySpring}
\end{figure}

\begin{figure}[h]
\centerline{\ $
\begin{array}{c@{\hspace{4ex}}c}
a)\psfig{file=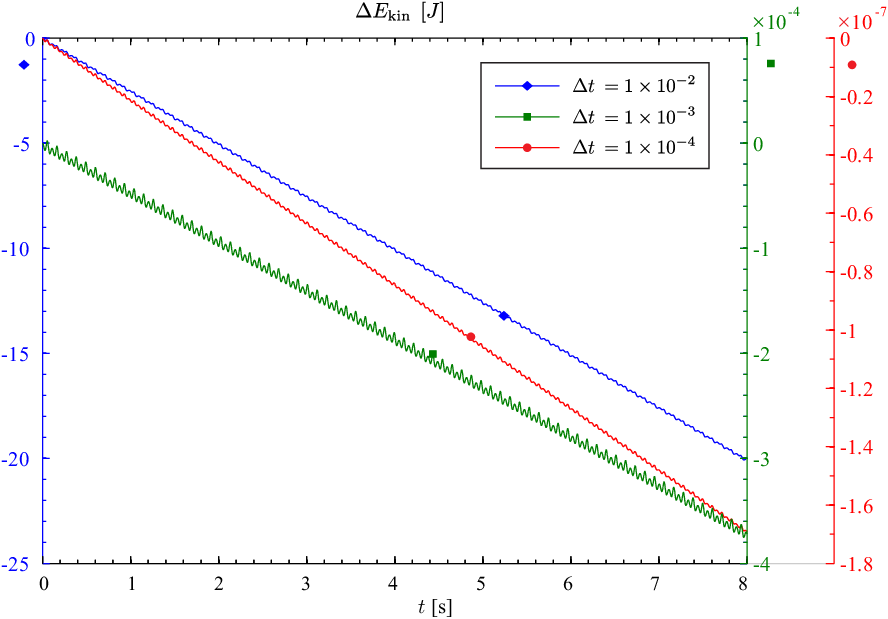,width=8cm} & b)\psfig{file=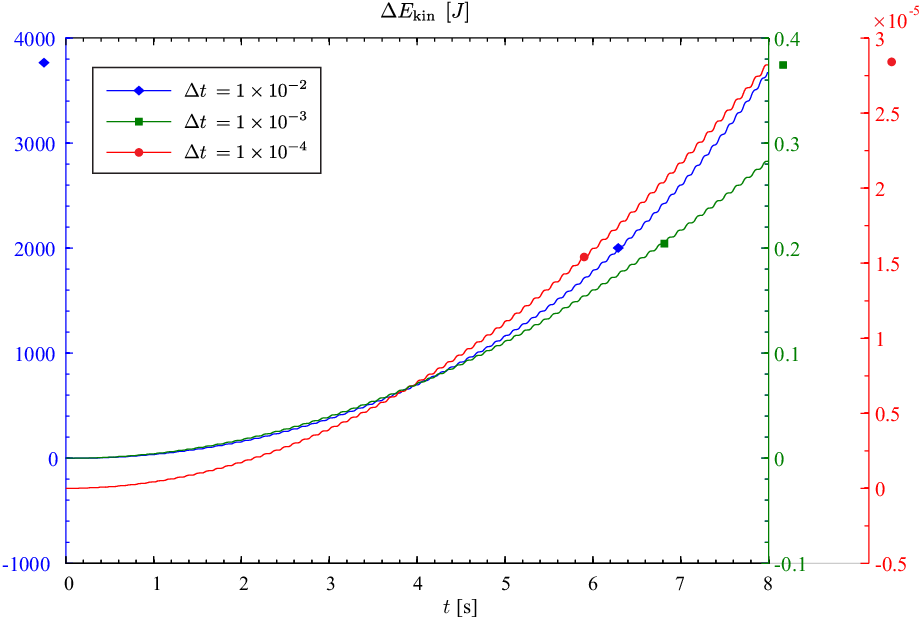,width=8cm}
\end{array}
$ }
\caption{Drift of the the total energy when integrating a) the $SE(3)$, and
b) the $SO(3)\times {\Bbb R}^{3}$ formulation.}
\label{FigEkinHeavyTop_GravitySpring}
\end{figure}

Figure \ref{FigErrorHeavyTop_GravitySpring} reveals that the $SE\left(
3\right) $ update preserves the center of rotation with accuracy closed to
the computation precision for all three step size values (shown is the error
in the translation constraints $\varepsilon \left( C\right) =\left\Vert
h\left( C\right) \right\Vert $). The $SO\left( 3\right) \times {\Bbb R}^{3}$
formulation on the other hand leads to a significant drift. The difference
is clearly visible in the error of the total energy (starting from $E\left(
0\right) =5000.69$~J) in figure \ref{FigEkinHeavyTop_GravitySpring}.
Although the configuration update does not affect the order of accuracy of
the overall integration scheme it is also beneficial for the energy
conservation.

The pure rotational motion of the top belongs to the rotation group as a $%
SE\left( 3\right) $ subgroup. But since it is not a rotation about the
origin of the reference frame (leading to coupled translations), it is not a
subgroup of $SO\left( 3\right) \times {\Bbb R}^{3}$.


\clearpage

\subsection{Spherical Double Pendulum in the Gravity Field}

\begin{figure}[h]
\centerline{\ \psfig{file=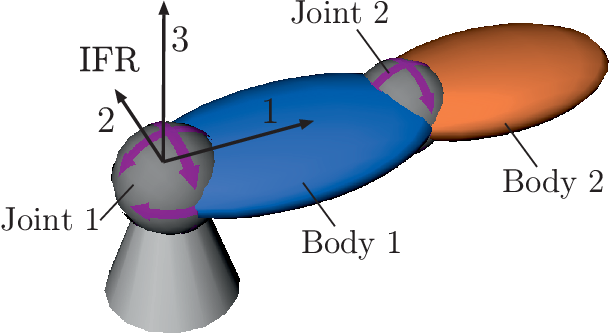,width=5cm}
} \vspace{-1ex}
\caption{Spherical double pendulum.}
\label{figSphericalPendulum}
\end{figure}

The double pendulum as shown in figure \ref{figSphericalPendulum} consists
of the two rigid bodies. The two bodies are interconnected and the pendulum
as a whole is connected to the ground by the spherical joints. The two links
are flat boxes with the side length $a,b,c$ along the axes of the COM
reference frame. Both have the same dimension with the lengths $a=0.2\,$m$%
,b=0.1\,$m$,c=0.05\,$m. Figure \ref{figSphericalPendulum} shows the inertia
ellipsoids. The configuration of the system is represented by $C_{1}=\left( 
{\bf R}_{1},{\bf r}_{1}^{\text{s}}\right) $ and $C_{2}=\left( {\bf R}_{2},%
{\bf r}_{2}^{\text{s}}\right) $. The two links are subject to the geometric
constraints \vspace{-5ex}

\begin{eqnarray}
h_{1}\left( C_{1}\right) &=&{\bf R}_{1}{\bf r}_{0}+{\bf r}_{1}^{\text{s}}=%
{\bf 0}  \nonumber \\
h_{2}\left( C_{1},C_{2}\right) &=&{\bf R}_{1}{\bf r}_{10}+{\bf r}_{1}^{\text{%
s}}-{\bf R}_{2}{\bf r}_{20}-{\bf r}_{2}^{\text{s}}={\bf 0}
\label{GeomConstrPend}
\end{eqnarray}%
%
%
%
%
%
where ${\bf r}_{i0},i=1,2$ is the position vector from the COM frame on body 
$i$ to the spherical joint connecting the two links, expressed in the COM
frame. The variable ${\bf r}_{0}$ is the position vector from the COM frame
on body 1 to the spherical joint connecting it to the ground expressed in
this COM frame. Denote with ${\Theta }_{i0}$ the inertia tensor of body $i$
w.r.t. the COM frame, and with $m_{i}$ its mass.

\paragraph{Motion Equations in the Body-Fixed Representation}

The velocity and acceleration constraints corresponding to (\ref%
{GeomConstrPend}), in the terms of the body-fixed twists ${\bf V}_{1},{\bf V}%
_{2}\in se\left( 3\right) $, are \vspace{-5ex}

\begin{eqnarray}
{\bf 0} &=&\left( 
\begin{array}{cccc}
\widehat{{\bf r}}_{0} & -{\bf I} & {\bf 0} & {\bf 0} \\ 
{\bf R}_{1}\widehat{{\bf r}}_{10} & -{\bf R}_{1} & -{\bf R}_{2}\widehat{{\bf %
r}}_{20} & {\bf R}_{2}%
\end{array}%
\right) \left( 
\begin{array}{c}
{\bf \omega }_{1} \\ 
{\bf v}_{1} \\ 
{\bf \omega }_{2} \\ 
{\bf v}_{2}%
\end{array}%
\right) ={\bf JV}  \label{velConstrPendulum} \\
{\bf J}\dot{{\bf V}} &=&-\dot{{\bf J}}{\bf V},\ {\rm with\ }-\dot{{\bf J}}%
{\bf V}=\left( 
\begin{array}{c}
\widehat{{\bf \omega }}_{1}{\bf v}_{1}{\bf +}\widehat{{\bf \omega }}_{1}%
\widehat{{\bf \omega }}_{1}{\bf r}_{0} \\ 
{\bf R}_{1}\widehat{{\bf \omega }}_{1}\widehat{{\bf \omega }}_{1}{\bf r}%
_{10}-{\bf R}_{2}\widehat{{\bf \omega }}_{2}\widehat{{\bf \omega }}_{2}{\bf r%
}_{20}+{\bf R}_{1}\widehat{{\bf \omega }}_{1}{\bf v}_{1}-{\bf R}_{2}\widehat{%
{\bf \omega }}_{2}{\bf v}_{2}%
\end{array}%
\right) .  \label{accConstrPendulum}
\end{eqnarray}%
\vspace{-3ex}

Using a reference frame at the COM yields \vspace{-4ex}

\begin{equation}
\left( 
\begin{array}{cccccc}
{\Theta }_{10} & {\bf 0} & {\bf 0} & {\bf 0} & -\widehat{{\bf r}}_{0} & -%
\widehat{{\bf r}}_{10}{\bf R}_{1}^{T} \\ 
{\bf 0} & m_{1}{\bf I} & {\bf 0} & {\bf 0} & -{\bf I} & -{\bf R}_{1}^{T} \\ 
{\bf 0} & {\bf 0} & {\Theta }_{20} & {\bf 0} & {\bf 0} & \widehat{{\bf r}}%
_{20}{\bf R}_{2}^{T} \\ 
{\bf 0} & {\bf 0} & {\bf 0} & m_{2}{\bf I} & {\bf 0} & {\bf R}_{2}^{T} \\ 
\widehat{{\bf r}}_{0} & -{\bf I} & {\bf 0} & {\bf 0} & {\bf 0} & {\bf 0} \\ 
{\bf R}_{1}\widehat{{\bf r}}_{10} & -{\bf R}_{1} & -{\bf R}_{2}\widehat{{\bf %
r}}_{20} & {\bf R}_{2} & {\bf 0} & {\bf 0}%
\end{array}%
\right) 
\hspace{-1.2ex}%
\left( 
\begin{array}{c}
\dot{{\bf \omega }}_{1} \\ 
\dot{{\bf v}}_{1} \\ 
\dot{{\bf \omega }}_{2} \\ 
\dot{{\bf v}}_{2} \\ 
{\bf \lambda }_{1} \\ 
{\bf \lambda }_{2}%
\end{array}%
\right) 
\hspace{-1.2ex}%
=%
\hspace{-1ex}%
\left( 
\hspace{-1ex}%
\begin{array}{c}
-\widehat{{\bf \omega }}_{1}{\Theta }_{10}{\bf \omega }_{1} \\ 
{\bf F}_{1}-m_{1}\widehat{{\bf \omega }}_{1}{\bf v}_{1} \\ 
-\widehat{{\bf \omega }}_{2}{\Theta }_{20}{\bf \omega }_{2} \\ 
{\bf F}_{2}-m_{2}\widehat{{\bf \omega }}_{2}{\bf v}_{2} \\ 
\ast \\ 
\ast \ast%
\end{array}%
\right)
\end{equation}%
where $\ast $ and $\ast \ast $ are the terms in the two rows of the right
hand side of (\ref{accConstrPendulum}). Since only the gravity forces act on
the system, the body-fixed forces are ${\bf F}_{i}={\bf R}_{i}^{T}{\bf g}^{%
\text{s}}$, where ${\bf g}^{\text{s}}=\left( 0,0,-g\right) $ is the gravity
vector w.r.t. space-fixed frame. The Lagrange multiplier ${\bf \lambda }%
_{i}\in {\Bbb R}^{3}$ is the reaction force in the joint $i$.

\paragraph{Motion Equations in the Mixed Velocity Representation}

With the mixed velocities ${\bf V}_{i}^{\text{m}}=\left( {\bf \omega }_{i},%
{\bf v}_{i}^{\text{s}}\right) $, the velocity and the acceleration
constraints are%
\vspace{-6ex}%

\begin{eqnarray}
{\bf 0} &=&\left( 
\begin{array}{cccc}
{\bf R}_{1}\widehat{{\bf r}}_{0} & -{\bf I} & {\bf 0} & {\bf 0} \\ 
{\bf R}_{1}\widehat{{\bf r}}_{10} & -{\bf I} & -{\bf R}_{2}\widehat{{\bf r}}%
_{20} & {\bf I}%
\end{array}%
\right) \left( 
\begin{array}{c}
{\bf \omega }_{1} \\ 
{\bf v}_{1} \\ 
{\bf \omega }_{2} \\ 
{\bf v}_{2}%
\end{array}%
\right) ={\bf JV}^{\text{m}}  \label{velConstrPendulumHyb} \\
{\bf J}\dot{{\bf V}}^{\text{m}} &=&-\dot{{\bf J}}{\bf V}^{\text{m}},\ {\rm %
with\ }-\dot{{\bf J}}{\bf V}^{\text{m}}=\left( 
\begin{array}{c}
{\bf R}_{1}\widehat{{\bf \omega }}_{1}\widehat{{\bf \omega }}_{1}{\bf r}_{0}
\\ 
{\bf R}_{1}\widehat{{\bf \omega }}_{1}\widehat{{\bf \omega }}_{1}{\bf r}%
_{10}-{\bf R}_{2}\widehat{{\bf \omega }}_{2}\widehat{{\bf \omega }}_{2}{\bf r%
}_{20}%
\end{array}%
\right) .  \label{accConstrPendulumHyb}
\end{eqnarray}%
\vspace{-4ex}

The DAE index 1 motion equations are, with the force vectors ${\bf F}_{i}^{%
\text{s}}={\bf g}^{\text{s}}$,

\[
\left( 
\begin{array}{cccccc}
{\Theta }_{10} & {\bf 0} & {\bf 0} & {\bf 0} & -\widehat{{\bf r}}_{0}{\bf R}%
_{1}^{T} & -\widehat{{\bf r}}_{10}{\bf R}_{1}^{T} \\ 
{\bf 0} & m_{1}{\bf I} & {\bf 0} & {\bf 0} & -{\bf I} & -{\bf I} \\ 
{\bf 0} & {\bf 0} & {\Theta }_{20} & {\bf 0} & {\bf 0} & \widehat{{\bf r}}%
_{20}{\bf R}_{2}^{T} \\ 
{\bf 0} & {\bf 0} & {\bf 0} & m_{2}{\bf I} & {\bf 0} & {\bf I} \\ 
{\bf R}_{1}\widehat{{\bf r}}_{0} & -{\bf I} & {\bf 0} & {\bf 0} & {\bf 0} & 
{\bf 0} \\ 
{\bf R}_{1}\widehat{{\bf r}}_{10} & -{\bf I} & -{\bf R}_{2}\widehat{{\bf r}}%
_{20} & {\bf I} & {\bf 0} & {\bf 0}%
\end{array}%
\right) 
\hspace{-1.2ex}%
\left( 
\begin{array}{c}
\dot{{\bf \omega }}_{1} \\ 
\dot{{\bf v}}_{1}^{\text{s}} \\ 
\dot{{\bf \omega }}_{2} \\ 
\dot{{\bf v}}_{2}^{\text{s}} \\ 
{\bf \lambda }_{1} \\ 
{\bf \lambda }_{2}%
\end{array}%
\right) 
\hspace{-1.2ex}%
=%
\hspace{-1ex}%
\left( 
\hspace{-1ex}%
\begin{array}{c}
-\widehat{{\bf \omega }}_{1}{\Theta }_{10}{\bf \omega }_{1} \\ 
{\bf F}_{1}^{\text{s}} \\ 
-\widehat{{\bf \omega }}_{2}{\Theta }_{20}{\bf \omega }_{2} \\ 
{\bf F}_{2}^{\text{s}} \\ 
\ast \\ 
\ast \ast%
\end{array}%
\right) 
\]%
\vspace{-3ex}

The terms $\ast $ and $\ast \ast $ are the two rows entries at the right
hand side of (\ref{accConstrPendulumHyb}).

\paragraph{Numerical Results}

In the initial configuration $C_{1}\left( 0\right) =\left( {\bf I},\left(
a_{1}/2,0,0\right) \right) $, $C_{2}\left( 0\right) =\left( {\bf I},\left(
a_{1}+a_{2}/2,0,0\right) \right) $ the pendulum is aligned with the
space-fixed x-axis as shown in figure \ref{figSphericalPendulum}. The
pendulum is moving in the gravity field with the initial velocities set to $%
{\bf \omega }_{10}=(10,0,0)$\thinspace rad/s and ${\bf \omega }%
_{20}=(10\,\pi ,10\,\pi ,20\,\pi )$\thinspace rad/s.

\begin{figure}[th]
\centerline{\ $
\begin{array}{c@{\hspace{4ex}}c}
a)\psfig{file=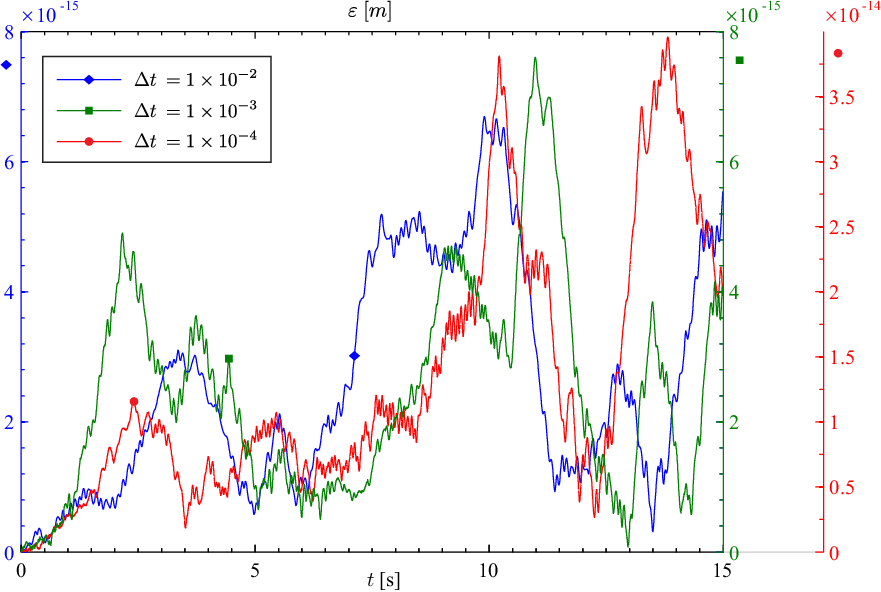,width=8cm} & b)\psfig{file=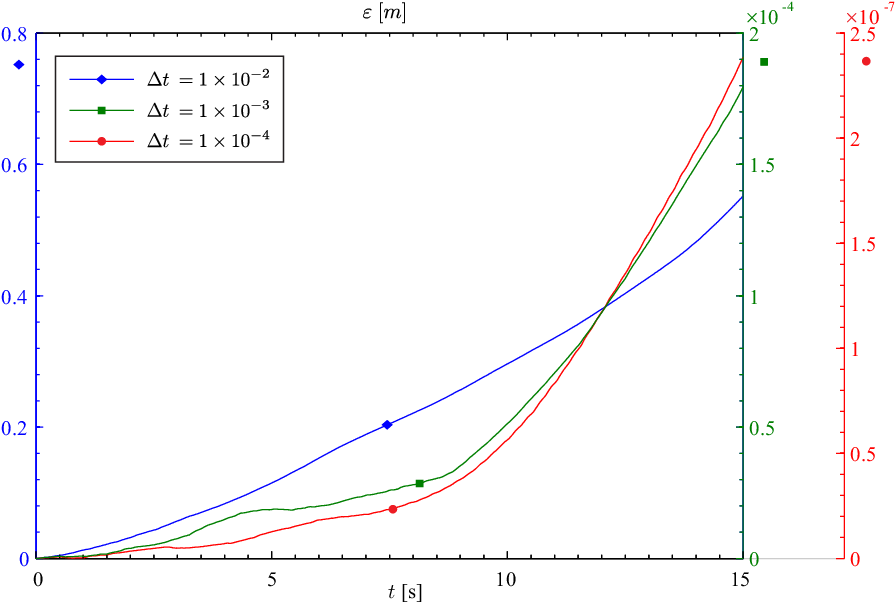,width=8cm}
\end{array}
$ }
\caption{Violation of geometric constraints of joint 1 when integrating the $%
SE(3)$ (a), and $SO(3)\times {\Bbb R}^{3}$ (b) formulation. \protect\vspace{%
-2ex}}
\label{DoublePendPosError1}
\end{figure}

\begin{figure}[b]
\centerline{\ $
\begin{array}{c@{\hspace{4ex}}c}
a)\psfig{file=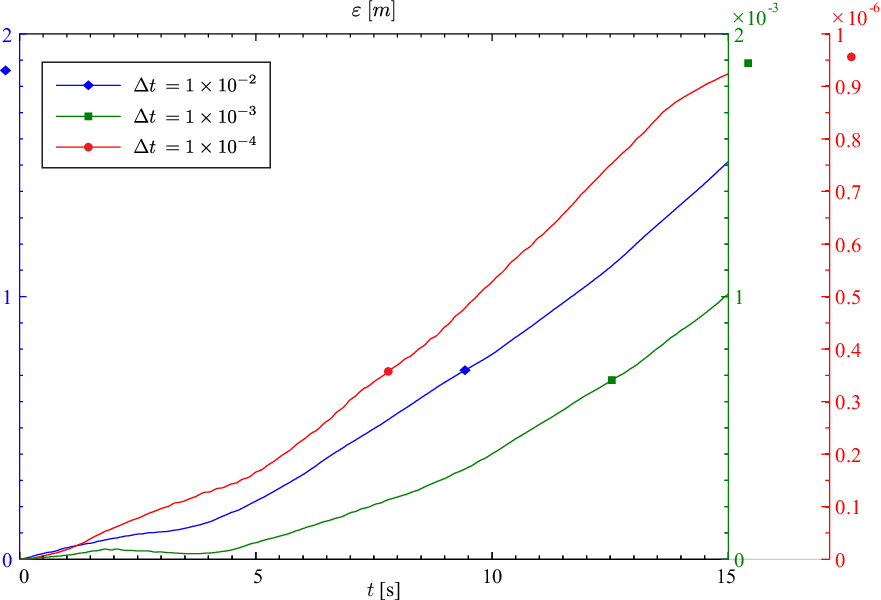,width=8cm} & b)\psfig{file=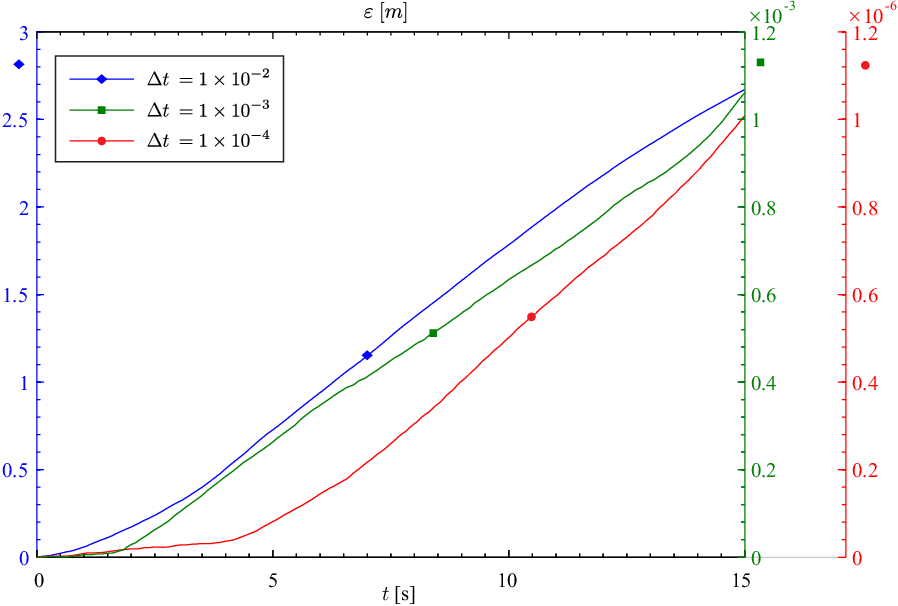,width=8cm}
\end{array}
$ }
\caption{Violation of geometric constraints of joint 2 when integrating a)
the $SE(3)$, and b) $SO(3)\times {\Bbb R}^{3}$ formulation. \protect\vspace{%
-2ex}}
\label{DoublePendPosError2}
\end{figure}

The figure \ref{DoublePendPosError1}a) shows the exact constraint
satisfaction up to the computation precision for joint 1 for $SE\left(
3\right) $. Both formulations achieve a similar constraint satisfaction for
joint 2 (figure \ref{DoublePendPosError2}). The joint 1 constrains the body
1 to perform the rotational motion, i.e. to a subgroup of $SE\left( 3\right) 
$. Thus, the condition of the corollary 1 is satisfied. Overall, the $%
SE\left( 3\right) $ update performs best as concluded in corollary \ref%
{corollary1} (for joint 1) and \ref{corollary3} (for joint 2).

\begin{figure}[h]
\centerline{\ $
\begin{array}{c@{\hspace{4ex}}c}
a)\psfig{file=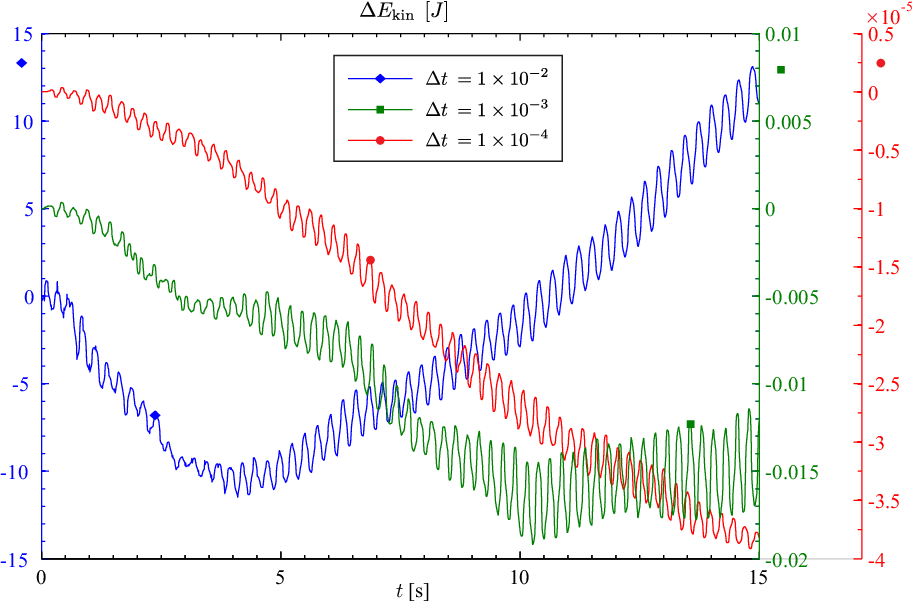,width=8cm} & b)\psfig{file=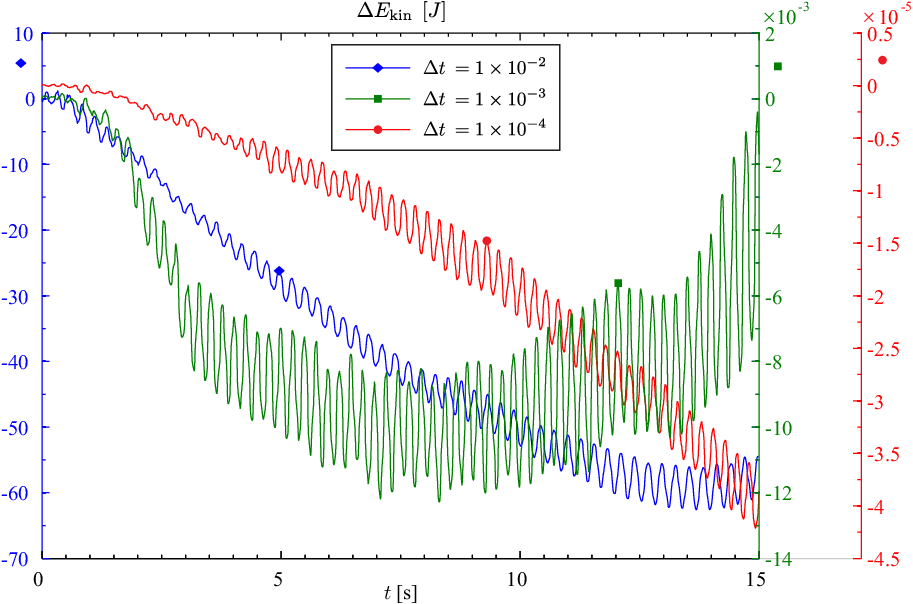,width=8cm}
\end{array}
$ }
\caption{Drift of total energy for a) the $SE(3)$, and b) $SO(3)\times {\Bbb %
R}^{3}$ update. \protect\vspace{-2ex}}
\label{DoublePendEnergyDrift}
\end{figure}

The constraint violation is also reflected in the drift of the total energy
in figure \ref{DoublePendEnergyDrift} that becomes significant for the $%
SO\left( 3\right) \times {\Bbb R}^{3}$ update with the large step size.
Similar results were obtained for a spherical three-bar and a four-bar
pendulum.

\vspace{-4ex}

\subsection{Planar 4-Bar Mechanism}

\begin{figure}[h]
\centerline{\
\psfig{file=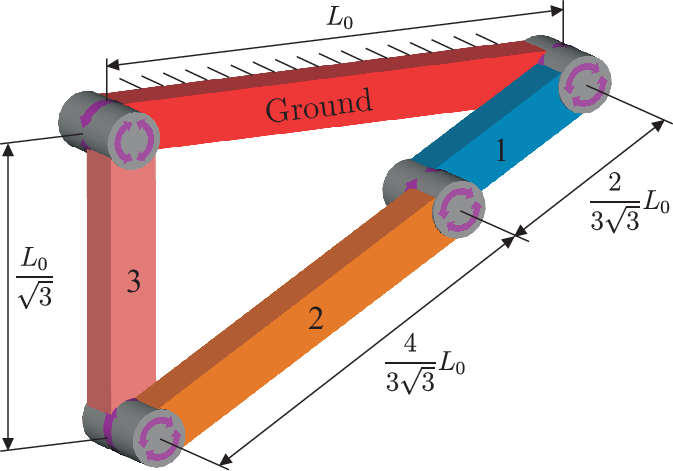,width=5cm} }
\caption{Planar 4-bar mechanism comprising four revolute joints with
parallel axes.}
\label{Fig4BarReference}
\end{figure}

The closed loop planar 4-bar mechanism in figure \ref{Fig4BarReference}
consists of the three rigid bodies linked to the ground. Body 1, 2, and 3
are mutually connected with the revolute joints. Furthermore, the bodies 1
and 3 are connected to the ground by the revolute joints. The motion
equations are omitted for the sake of the brevity. The geometric parameters
are indicated in figure \ref{Fig4BarReference}, and $L_{0}=0.5\,$m is used
in the simulations.

\begin{figure}[h]
\centerline{\ $
\begin{array}{c@{\hspace{4ex}}c}
a)\psfig{file=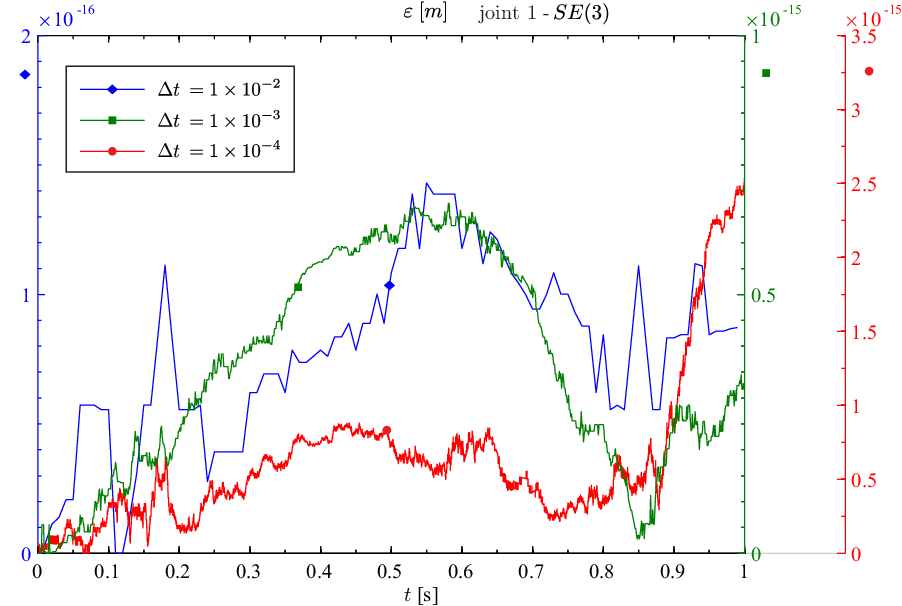,width=8cm} & b)\psfig{file=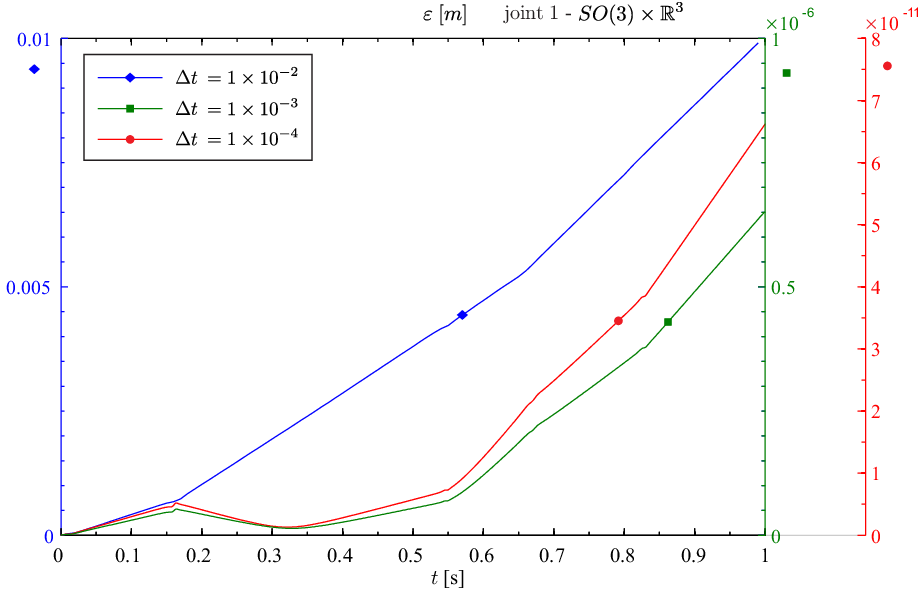,width=8cm} \\ 
c)\psfig{file=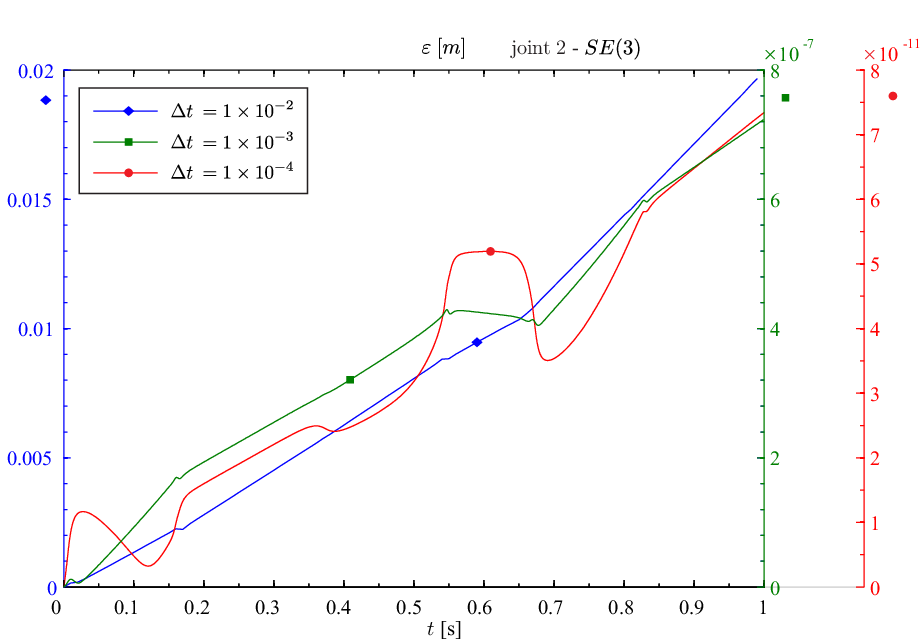,width=8cm} & d)\psfig{file=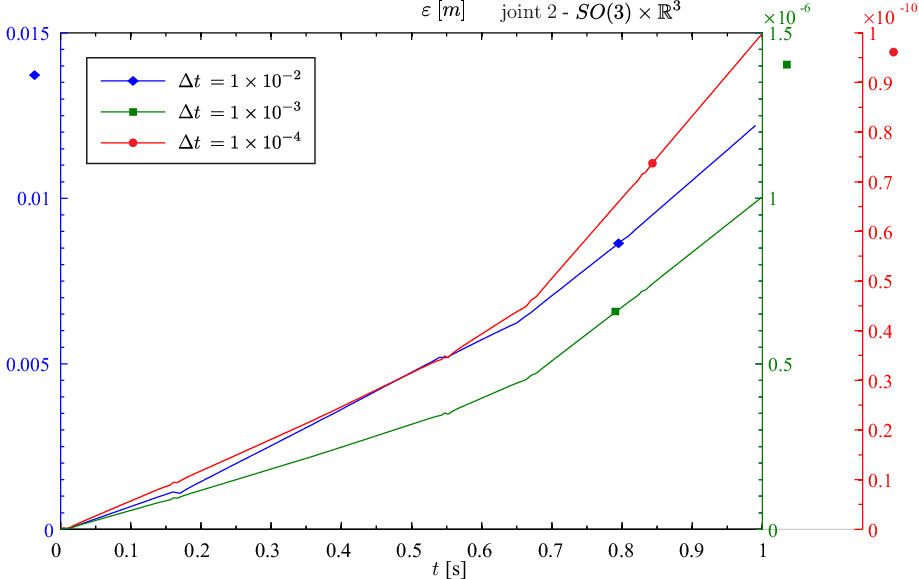,width=8cm}
\end{array}
$ }
\caption{Violation of geometric constraints of revolute joint 1 and 2 when
integrating a,c) the $SE(3)$, and b,d) the $SO(3)\times {\Bbb R}^{3}$
formulation. \protect\vspace{-2ex}}
\label{FigError_Planar4BarJoint2}
\end{figure}

\begin{figure}[h]
\centerline{\ $
\begin{array}{c@{\hspace{4ex}}c}
a)\psfig{file=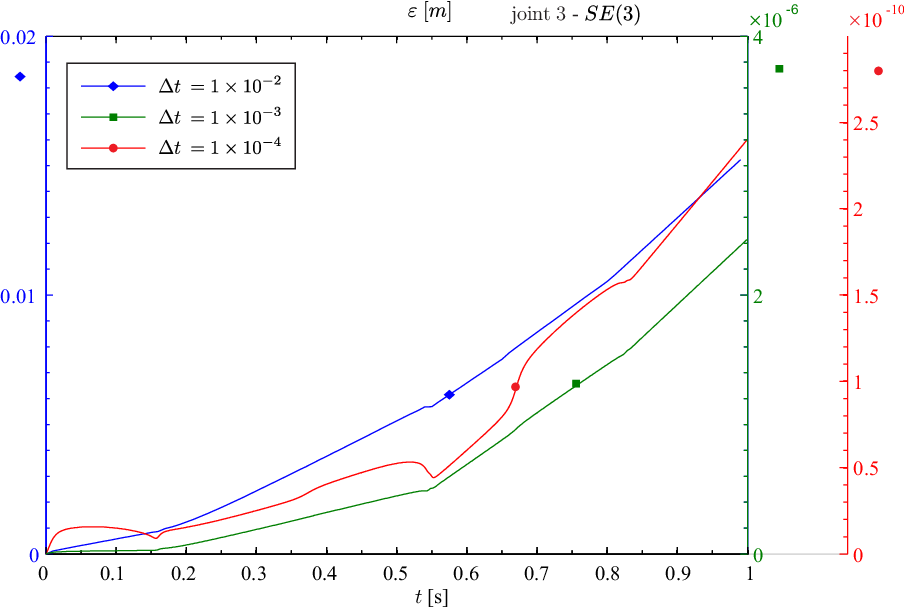,width=8cm} & b)\psfig{file=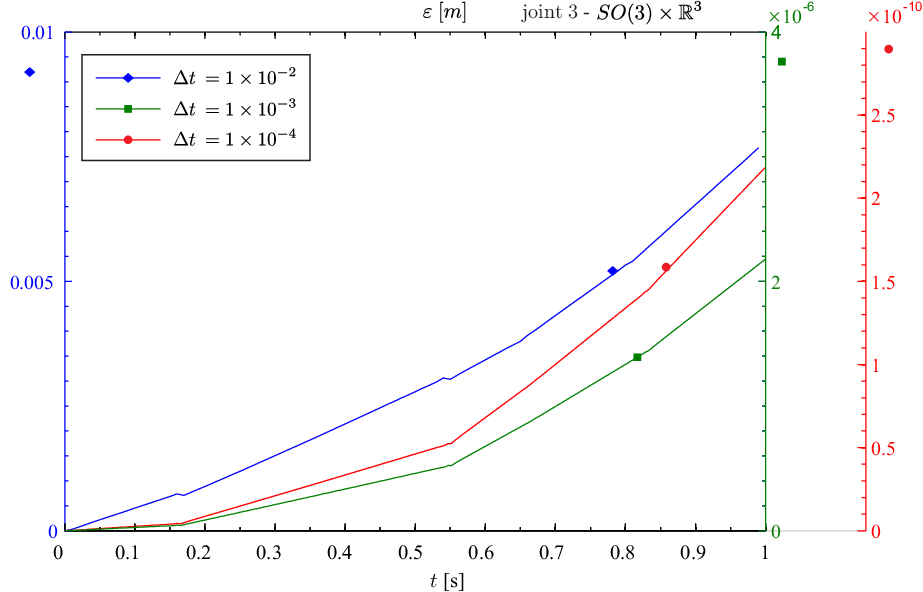,width=8cm} \\ 
c)\psfig{file=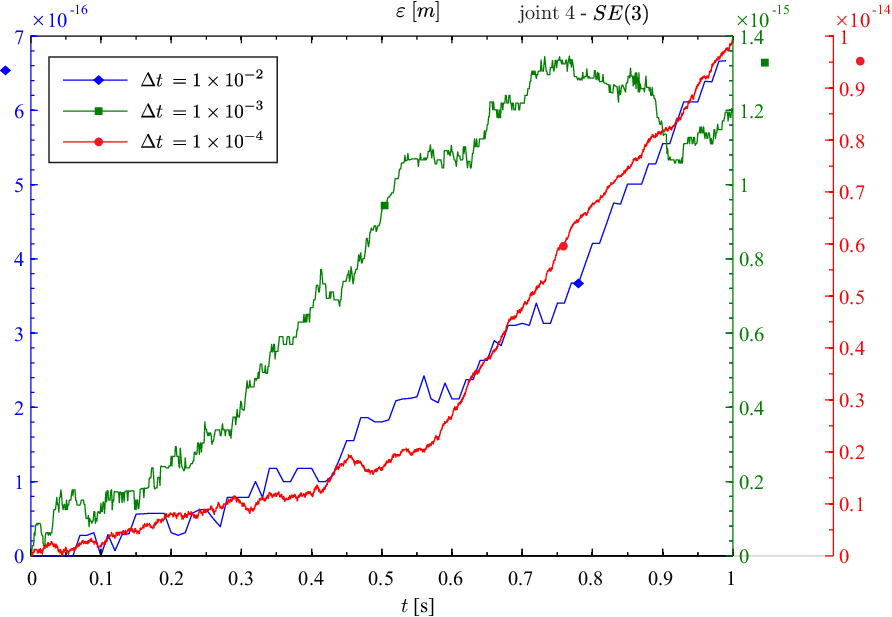,width=8cm} & d)\psfig{file=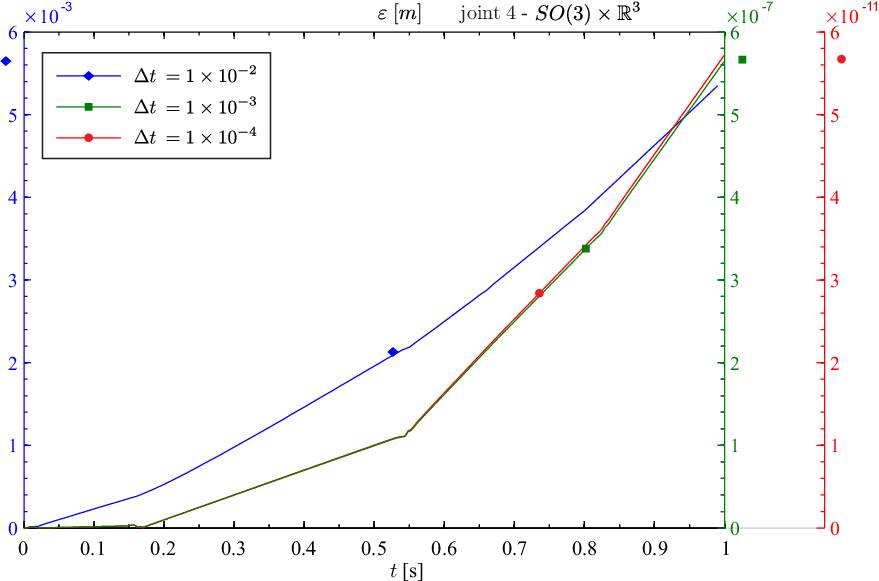,width=8cm}
\end{array}
$ }
\caption{Violation of geometric constraints of revolute joint 3 and
spherical joint 4 when integrating a,c) the $SE(3)$, and b,d) the $%
SO(3)\times {\Bbb R}^{3}$ formulation. \protect\vspace{-2ex}}
\label{FigError_Planar4BarJoint4}
\end{figure}

The initial configuration is shown in figure \ref{Fig4BarReference}. The
initial velocities are set in the way that the input crank (body 1) rotates
with the angular velocity $\omega _{0}=10\,\pi $\thinspace rad/s. Thus the
motion equations were integrated with the initial conditions ${\bf \omega }%
_{10}=\left( 0,0,\omega _{0}\right) ,{\bf \omega }_{20}=\left( 0,0,-\omega
_{0}/2\right) ,{\bf \omega }_{30}=\left( 0,0,0\right) $, ${\bf v}%
_{10}=\left( L_{0}\omega _{0}/(3\sqrt{3}),0,0\right) ,{\bf v}_{20}=\left(
0,2L_{0}\omega _{0}/(3\sqrt{3}),0,0\right) ,{\bf v}_{30}=\left( 0,0,0\right) 
$. The numerical results are shown in figures \ref{FigError_Planar4BarJoint2}
and \ref{FigError_Planar4BarJoint4}. The orientation constraints are exactly
satisfied by the both formulations. This example confirms again the
statement that the $SE\left( 3\right) $ update exactly preserves the lower
pair constraints restricting the body 1 and 3 to a subgroup (corollary 1),
and also that the joint constraints are exactly satisfied (lemma \ref{lemma1}%
). In this example, the bodies 1 and 3 are constrained so to perform the
rotational motion about a constant axis.


\clearpage

\subsection{RP Mechanisms%
\label{secRP}%
\label{secRP}}

\begin{figure}[h]
\centerline{\ \psfig{file=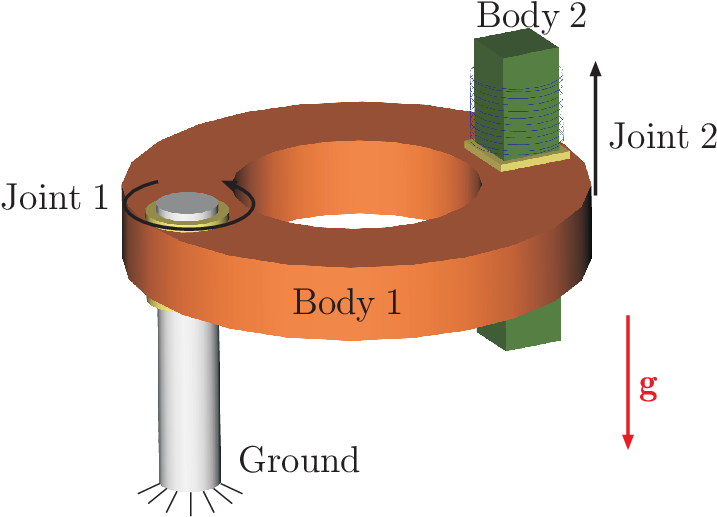,width=6cm} }
\caption{RP kinematic chain: Two bodies connected by a prismatic joint 2,
and connected to the ground by revolute joint 1. \protect\vspace{-2ex}}
\label{FigRP_reference}
\end{figure}

As next example, the mechanism in figure \ref{FigRP_reference} is considered
consisting of a cylindrical ring (body 1) connected to the ground by a
revolute joint 1 and to another body (body 2) by a prismatic joint 2. The
two bodies are elastically coupled by a longitudinal spring along the
prismatic joint. The system is moving in the gravity field as indicated in
figure \ref{FigRP_reference} with the gravity vector ${\bf g}^{\text{s}%
}=\left( 0,0,-9.81\right) \,$m/s$^{2}$. By assuming that material is
aluminium, the respective mass is $m_{1}=6.82825$\thinspace kg and $%
m_{2}=0.864$\thinspace kg. The inertia matrix of body 1 and 2, w.r.t. the
body-fixed reference frames at the COM, is ${\bf \Theta }_{1}={\rm diag}%
\,\left( 0.0507567,0.0507567,0.0986682\right) \,$kg\thinspace m$^{2}$ and $%
{\bf \Theta }_{2}={\rm diag}\,\left( 0.0002304,0.0029952,0.0029952\right) \,$%
kg\thinspace m$^{2}$ , respectively. The dynamics of the bodies is governed
by the Newton-Euler equations supplemented by the joint constraints (not
shown here for the sake of brevity). The spring stiffness is set to $10^{4}$%
\thinspace N/m.

The motion equations are integrated with the initial conditions ${\bf \omega 
}_{10}={\bf \omega }_{20}=(0,0,20\,\pi )\,$rad/s and ${\bf v}_{10}={\bf r}%
_{1}\times {\bf \omega }_{10}$, ${\bf v}_{20}={\bf r}_{21}\times {\bf \omega 
}_{20}+(0,0,1)\,$m/s, where ${\bf r}_{1i}$ is the position vector of the
axis of joint 1 w.r.t. to the COM frame on the body $i$. That is, the
prismatic joint 2 is given an initial velocity of 1\thinspace m/s.

\begin{figure}[h]
\centerline{\ $
\begin{array}{c@{\hspace{4ex}}c}
a)\psfig{file=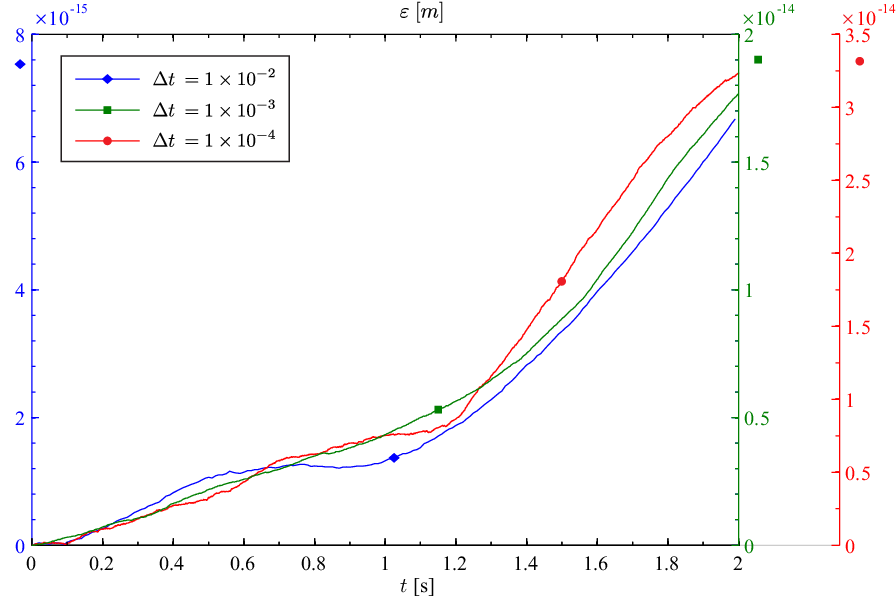,width=8cm} & b)\psfig{file=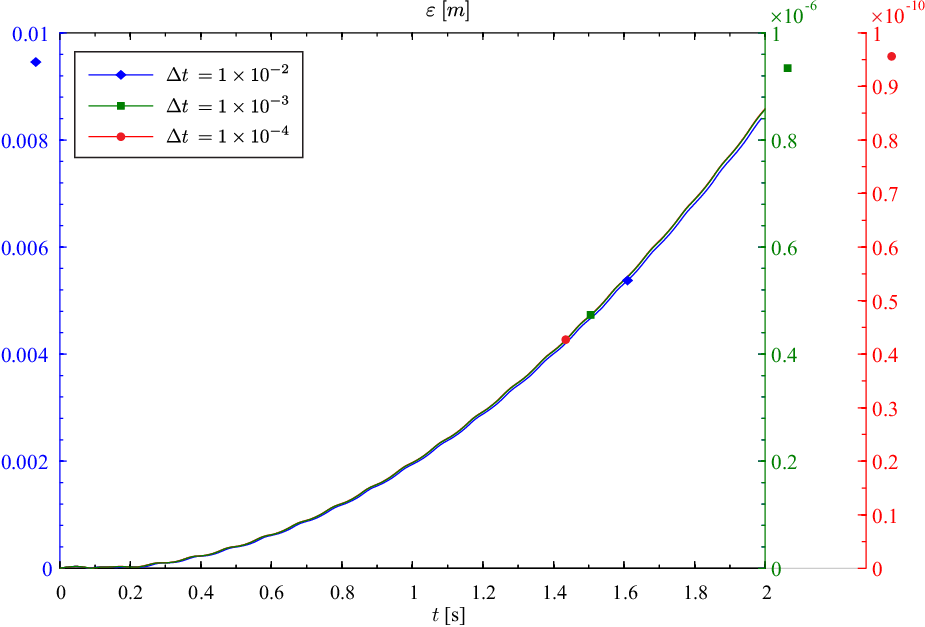,width=8cm}
\end{array}
$ }
\caption{Violation of geometric constraints of revolute joint 1 when
integrating a) the $SE(3)$, and b) $SO(3)\times {\Bbb R}^{3}$ formulation. 
\protect\vspace{-2ex}}
\label{FigErrorRP_Joint1}
\end{figure}

The numerical results in figure \ref{FigErrorRP_Joint1}a) and \ref%
{FigErrorRP_Joint2}a) show the constraint satisfaction, indicating an
excellent motion reconstruction, when the $SE(3)$ update is used. Using the $%
SO(3)\times {\Bbb R}^{3}$ formulation yields the accuracy depending on the
integration step size as shown in figures \ref{FigErrorRP_Joint1}b) and \ref%
{FigErrorRP_Joint2}b). This is an example where lemma \ref{lemma1} applies.
Both bodies are free to move in a $SE(3)$ subgroup and the actions of the
revolute and prismatic joint commute. The body 1 performs a pure rotational
motion (subgroup $SO\left( 2\right) $) and the body 2 a 'cylindrical'
motion, i.e. a rotation about a fixed axis plus a translation along this
axis (subgroup $SO\left( 2\right) \times {\Bbb R}$). The same accuracy is
achieved when the joint 1 is a cylindrical joint. Note again that the step
size can be arbitrarily large and the $SE\left( 3\right) $ update still
satisfies the kinematic joint constraints.

Using $SO(3)\times {\Bbb R}^{3}$ as c-space, the motion space of body 1 is
only a submanifold of the subgroup $SO(2)\times {\Bbb R}^{2}$ of the
c-space, and the motion space of the body 2 is a submanifold of the c-space
subgroup $SO(2)\times {\Bbb R}^{3}$, thus this update cannot satisfy the
joint constraints. 

\begin{figure}[h]
\centerline{\ $
\begin{array}{c@{\hspace{4ex}}c}
a)\psfig{file=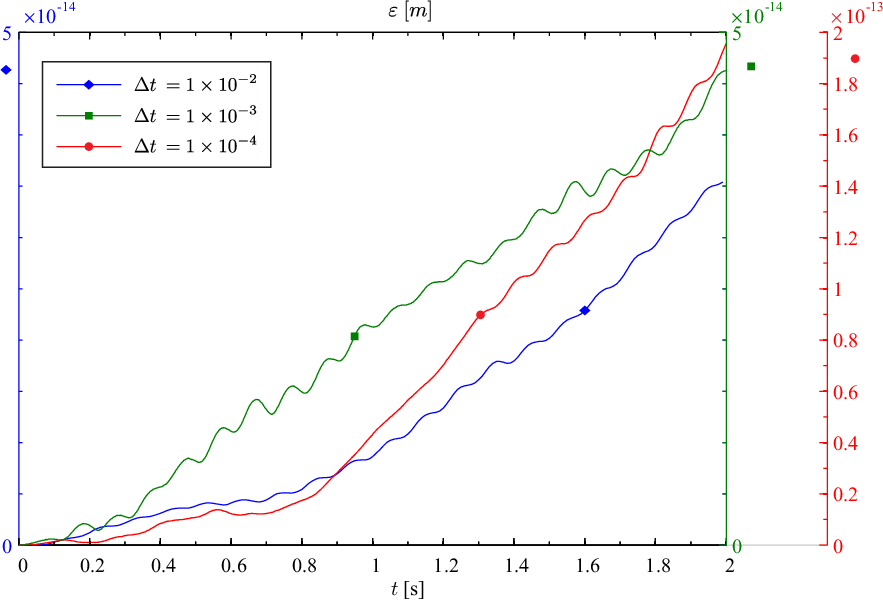,width=8cm} & b)\psfig{file=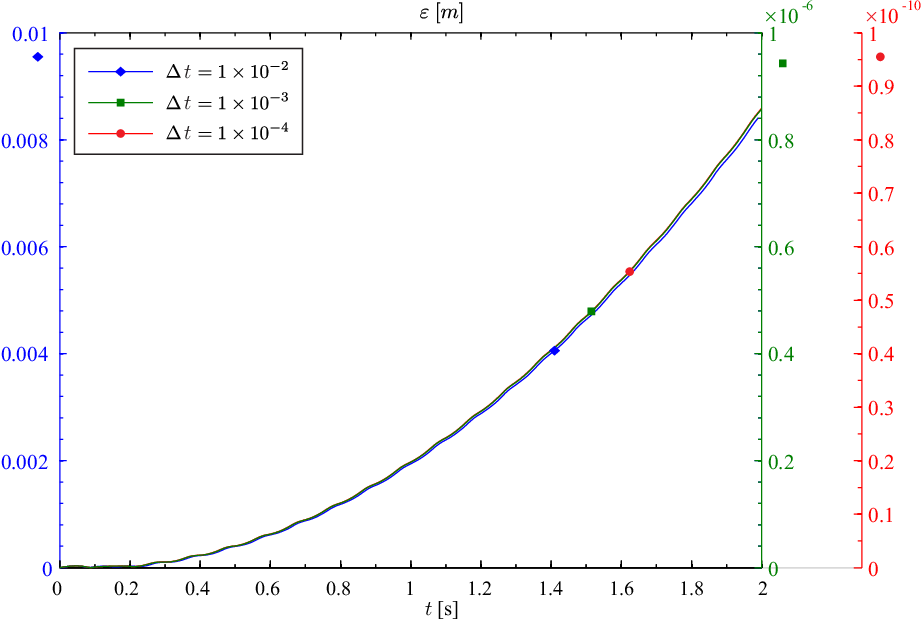,width=8cm}
\end{array}
$ }
\caption{Violation of geometric constraints of prismatic joint 2 when
integrating a) the $SE(3)$, and b) $SO(3)\times {\Bbb R}^{3}$ formulation. 
\protect\vspace{-2ex}}
\label{FigErrorRP_Joint2}
\end{figure}

\newpage

\subsection{Cardanic Transmission}

\begin{figure}[h]
\centerline{\psfig{file=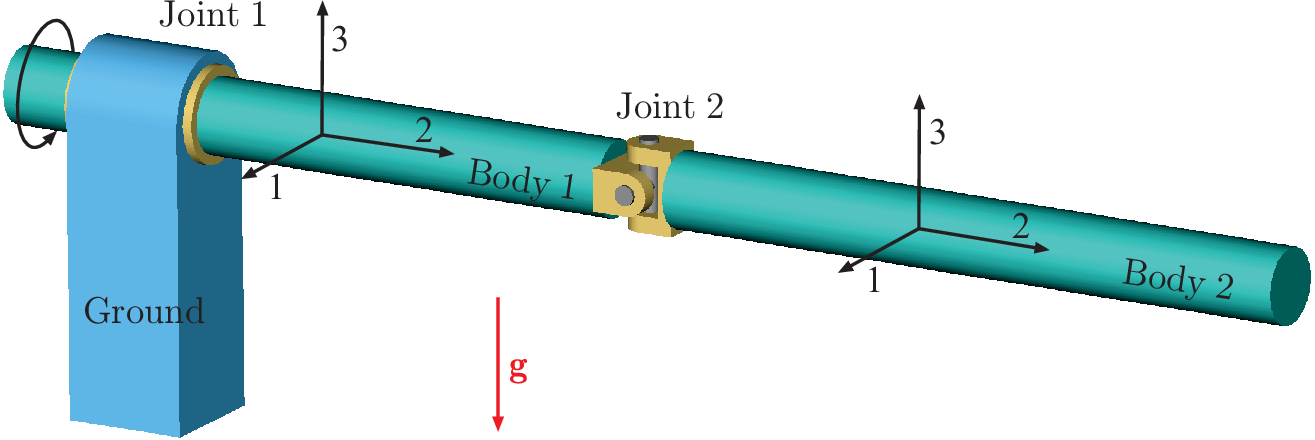,width=7.5cm} }
\caption{A Cardanic transmission consisting of an input shaft (body 1),
mounted at the ground by a revolute joint 1, connected to the drive shaft
(body 2) by a universal joint 2. \protect\vspace{-2ex}}
\label{FigHook}
\end{figure}

The last example is a Cardanic transmission presented in figure \ref{FigHook}%
. The input shaft (body 1), mounted at the ground by a revolute joint (joint
1), is connected to a drive shaft (body 2) by a universal joint (joint 2).
In this example, the body 1 is constrained to perform rotation about the
fixed axis of the revolute joint 1. Since the reference (COM) frame of the
body 1 is located at the rotation axis of the revolute joint, there is no
translation component. Hence, the motion observed at this point is a pure
rotation and the both formulations yield a perfect constraint satisfaction
(otherwise $SO(3)\times {\Bbb R}^{3}$ would not). The motion equations are
integrated with initial velocities ${\bf \omega }_{1}=\left( 0,\pi ,0\right) 
$ rad/s, and ${\bf \omega }_{2}=\left( \pi ,\pi ,0\right) $ rad/s. In the
initial configuration the two bodies are aligned along the 3-axis.

\begin{figure}[h]
\centerline{\ $
\begin{array}{c@{\hspace{4ex}}c}
a)\psfig{file=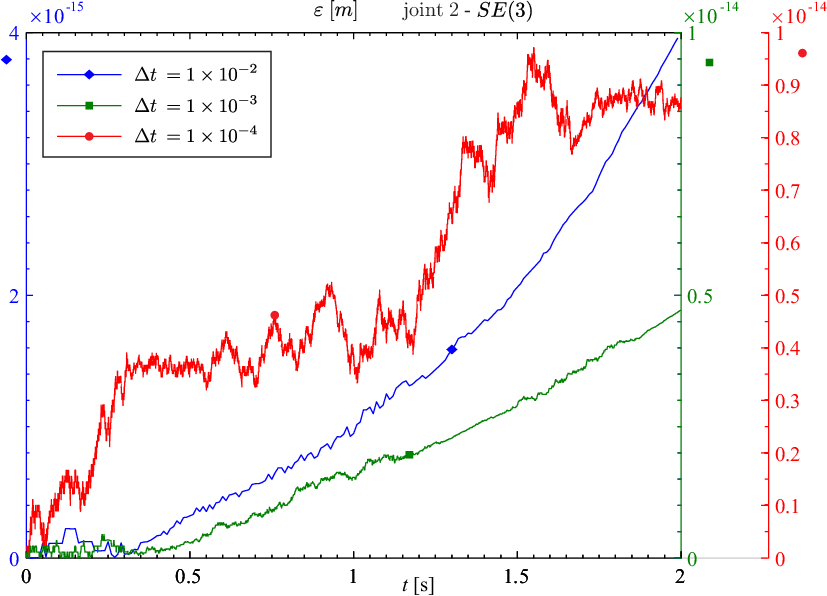,width=8cm} & b)\psfig{file=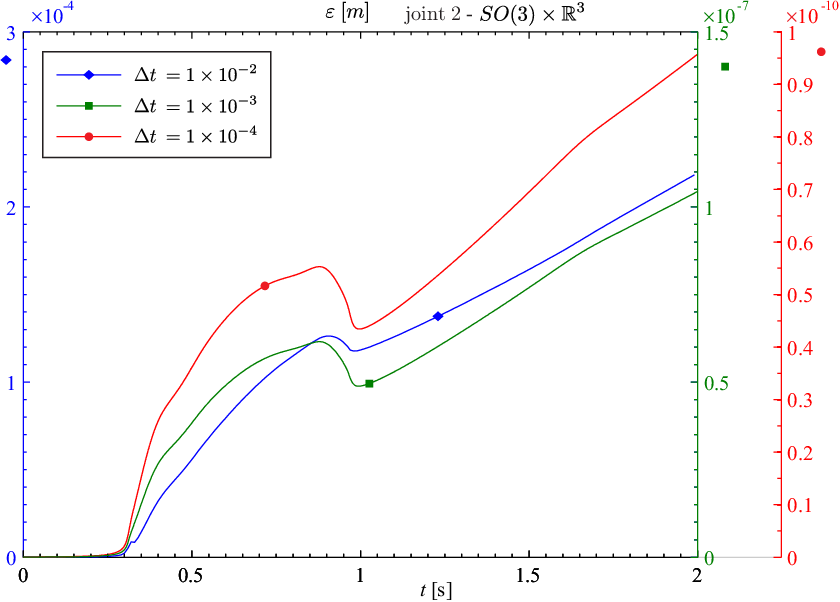,width=8cm}
\end{array}
$ }
\caption{Violation of position constraints of universal joint 2 when
integrating a) the $SE(3)$, and b) $SO(3)\times {\Bbb R}^{3}$ formulation. 
\protect\vspace{-2ex}}
\label{FigPosErrorHook_Joint2}
\end{figure}

The universal joint connecting body 1 and body 2 adds two rotational degrees
of freedom (DOF) so that body 2 is constrained to perform free spatial
rotations about the intersection point of the two hook joint axes. The
combination of the two joints is equivalent to a spherical joint
constraining the reference frame on body 2 to move on a sphere centered at
the hook joint. That is, body 2 is constrained to the subgroup $SO(3)$ of $%
SE\left( 3\right) $. Due to the translation components, the motion does,
however, not form a subgroup of $SO(3)\times {\Bbb R}^{3}$. The condition of
corollary \ref{corollary1} is thus fulfilled for the $SE\left( 3\right) $
update only. Consequently, the position constraints of the joint 2 should be
perfectly satisfied. This is confirmed in figure \ref{FigPosErrorHook_Joint2}%
a). This is not so for the $SO(3)\times {\Bbb R}^{3}$ update as shown in
figure \ref{FigPosErrorHook_Joint2}b).

\begin{figure}[h]
\centerline{\ $
\begin{array}{c@{\hspace{4ex}}c}
a)\psfig{file=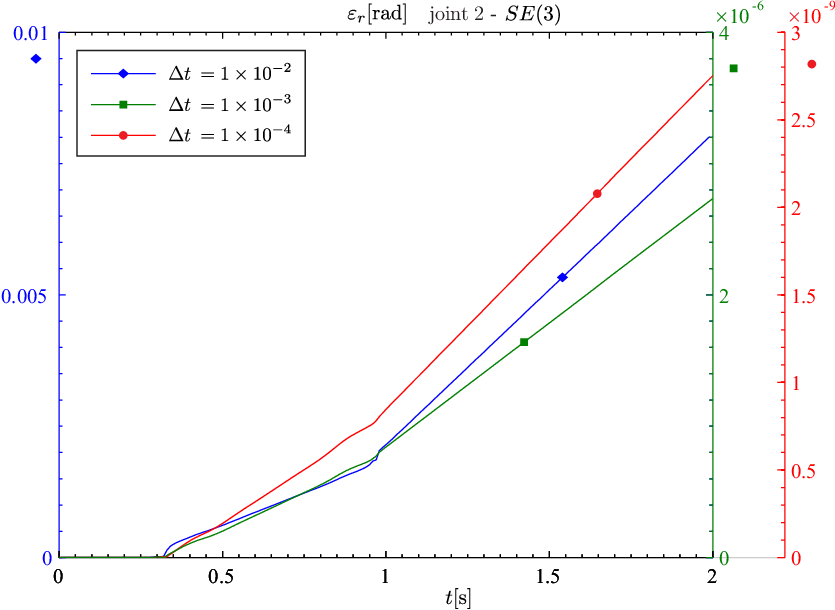,width=8cm} & b)\psfig{file=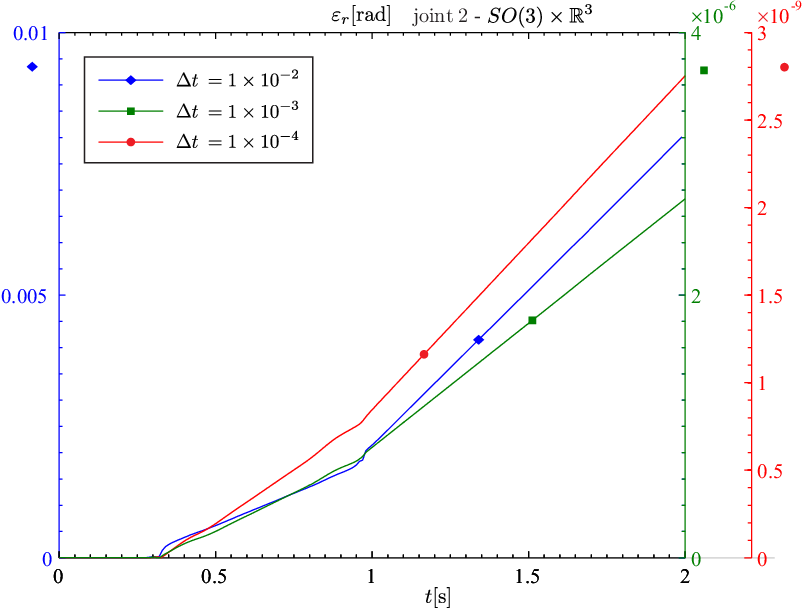,width=8cm}
\end{array}
$ }
\caption{Violation of orientation constraints of universal joint 2 when
integrating a) the $SE(3)$, and b) $SO(3)\times {\Bbb R}^{3}$ formulation. 
\protect\vspace{-2ex}}
\label{FigRotErrorHook_Joint2}
\end{figure}

The hook joint does not define a subgroup of $SE\left( 3\right) $ or $%
SO\left( 3\right) \times {\Bbb R}^{3}$ so that the condition of lemma \ref%
{lemma1} is not satisfied. Therefore, the orientation constraints of the
hook joint are not exactly satisfied, but only according to the time step
size and order of accuracy of the RK4 integration method. This is documented
in figure \ref{FigRotErrorHook_Joint2} revealing the same constraint
satisfaction of both configuration update variants. 

\newpage

\section{Discussion and Conclusion%
\label{secConclusion}%
\label{secConclusion}}

The dynamics of a constrained MBS comprising rigid bodies is governed by the
Newton-Euler equations (\ref{BH}a) together with the kinematic
reconstruction equations (\ref{BH}b). The latter reveal the actual motion of
the MBS corresponding to its velocity. Moreover, in a numerical integration
scheme these are used to estimate finite motion increments. Since (\ref{BH}%
b) respectively (\ref{ODEvec1}b), (\ref{ODEvec2}b), and (\ref{ODE2}b),
reflect the geometry of the c-space Lie group, a different c-space leads to
different numerical results. This has a consequence for the constraint
satisfaction when (\ref{BH}) is replaced by the index 1 systems of the form (%
\ref{ODEvec1}), (\ref{ODEvec2}), or (\ref{ODE2}) and is solved with the ODE
integrators. While the accuracy of the numerical solution of the ODE is
determined by the integration method, the geometric constraints should be
perfectly satisfied, i.e. the MBS as a whole shall remain assembled. This
poses the question, if, and under which conditions, an integration scheme,
solving the ODE, can inherently satisfy the constraints independently of its
order of accuracy.

Taking a kinematic perspective it was shown in this paper that this issue is
connected with the particular form of the kinematic reconstruction
equations, i.e. the c-space Lie group, rather than with the integration
scheme. The two relevant c-space Lie groups are $SE\left( 3\right) $ and $%
SO\left( 3\right) \times {\Bbb R}^{3}$. The latter is commonly used (at
least implicitly) in modeling of rigid body MBS, whereas the importance of
using $SE\left( 3\right) $ for the finite element integration schemes has
already been recognized \cite%
{Borri2001a,Borri2001b,Borri2002,Borri2003,Bottasso2002}. The analysis is
facilitated by a general formulation of the motion equations accounting for
a general c-space Lie group. For completeness, the Lie groups $SE\left(
3\right) $ and $SO\left( 3\right) \times {\Bbb R}^{3}$ are also related to
the dual quaternion formulation and the Euler parameter formulation that are
frequently used in the applications. The actual motion equations are
formulated as a vector space ODE as well as an ODE on the c-space Lie group.

It was shown in section \ref{secDiscussion} that, if a body is constrained
to a subgroup of its c-space Lie group, the corresponding configuration
update always yields a body configuration within that subgroup, thus
satisfies the {\it motion constraints imposed on the body}, independently of
the actual integration accuracy. On the contrary, when the motion of a rigid
body within the MBS is not constrained to a subgroup, but to an arbitrary
submanifold, the update step can lead out of that submanifold and thus
violate the constraints. This conclusion does not immediately account for
the satisfaction of the joint constraints.

It is know that $SE\left( 3\right) $ is the proper Lie group of rigid body
motions. But it is shown that this only implies better constraint
satisfaction for special MBS. Since subgroups of $SE\left( 3\right) $
account for lower pair joints (Reuleaux pairs), an immediate consequence is
that the $SE\left( 3\right) $ update achieves a perfect satisfaction of the 
{\it joint constraints} when a rigid body is connected to the ground by a
lower pair joint. Another result is that the lower pair joint constraints
are exactly satisfied if the joint connects two bodies that are respectively
constrained to a $SE\left( 3\right) $ subgroup. The subgroups of the direct
product group $SO\left( 3\right) \times {\Bbb R}^{3}$, which is assumed by
the standard form (\ref{VelRigidBody}), account for situations with a less
practical relevance. As such, a prominent example is the free floating body
(section \ref{secFreeBody}) with the reference frame located at its COM,
where the momentum conservation constrains the body to rotate while moving
on a straight line. This motion belongs to a subgroup of $SO\left( 3\right)
\times {\Bbb R}^{3}$ but not of $SE\left( 3\right) $.

Numerical results are presented in section \ref{secExamples} confirming the
above conclusion. The results show that in the general cases, where the
constrained motions do not form a subgroup of either c-space Lie group, both 
$SE\left( 3\right) $ and $SO\left( 3\right) \times {\Bbb R}^{3}$ lead to the
same order of accuracy. Hence, $SE\left( 3\right) $ is superior only for a
special class of MBS. Nevertheless, it allows for the construction of
tailored simulation models.\ Since $SE\left( 3\right) $ can achieve exact
constraint satisfaction for certain MBS comprising lower pair joints, and
since for a general MBS the $SE\left( 3\right) $ and $SO\left( 3\right)
\times {\Bbb R}^{3}$ update scheme yield the same order of accuracy, it can
be concluded that $SE\left( 3\right) $ is the best choice for a c-space Lie
group. From the computational point, however, the expressions (\ref%
{dexpInvSE3Park}) and (\ref{dexpInvSE3}) used in the $SE\left( 3\right) $
update are slightly more complex than (\ref{Vhdexp}). In order to minimize
the overall complexity while ensuring the exact constraint satisfaction
whenever possible, a c-space Lie group can be tailored to a specific model.

Notice that the c-space only affects the kinematic equations (\ref{BH}b) but
not the motion equations (\ref{BH}a) or the integration scheme. From the
practical point of view, the question about the c-space boils down to the
appropriate configuration update since the c-space is merely a conceptual
construct. All said is equally valid for the Euler parameters and dual
quaternions as well as for any 3-angle parameterization.

Finally, it should be remarked that the presented results are relevant in
general since any DAE formulation (index 1, 2, or 3) of MBS in terms of the
non-holonomic velocities involves the kinematic reconstruction equations (%
\ref{BH}b).

\appendix{}

\section*{Appendix}

\section{Geometric Background on Rigid Body Motions\label{appendixA}}

\subsection{Spatial Rotation}

The relative orientation of two frames can be represented by a rotation
matrix ${\bf R}\in SO\left( 3\right) $. The special orthogonal group $%
SO\left( 3\right) $ is the Lie group of \thinspace $3\times 3$ orthogonal
matrices with corresponding Lie algebra $so\left( 3\right) $ the vector
space of skew symmetric $3\times 3$ matrices. The latter is isomorphic to $%
{\Bbb R}^{3}$ via ${\bf x}\in {\Bbb R}^{3}\mapsto \widehat{{\bf x}}\in
so\left( 3\right) $, where $\widehat{{\bf x}}$ is the cross product matrix
associated to ${\bf x}$. A rotation can be associated an instantaneous
rotation vector ${\bf \xi }\in {\Bbb R}^{3}$. The finite rotation, i.e. the
rotation matrix, corresponding to the (generally not constant rotation
vector) is given by the exponential mapping 
\begin{equation}
\exp \widehat{{\bf \xi }}={\bf I}+\frac{\sin \left\Vert {\bf \xi }%
\right\Vert }{\left\Vert {\bf \xi }\right\Vert }\widehat{{\bf \xi }}+\frac{%
1-\cos \left\Vert {\bf \xi }\right\Vert }{\left\Vert {\bf \xi }\right\Vert }%
\widehat{{\bf \xi }}^{2}.  \label{expSO3}
\end{equation}%
This is the Euler-Rodriguez formula describing the rotation about the
rotation axis ${\bf \xi }/\left\Vert {\bf \xi }\right\Vert $ by the angle $%
\left\Vert {\bf \xi }\right\Vert $.

Since $so\left( 3\right) $ is isomorphic to ${\Bbb R}^{3}$ with the cross
product as Lie bracket, the Lie bracket $\left[ \widehat{{\bf x}},\widehat{%
{\bf y}}\right] =\widehat{{\bf x}}\widehat{{\bf y}}-\widehat{{\bf y}}%
\widehat{{\bf x}}=\widehat{{\bf x}\times {\bf y}}$ on $so\left( 3\right) $
is expressed in vector notation as $\left[ {\bf x},{\bf y}\right] =\widehat{%
{\bf x}}{\bf y}={\bf x}\times {\bf y}$.

Expressing the angular velocity in terms of the time derivative of the
rotation vector requires the differential of the exp mapping (\ref{expSO3}).
As for any Lie group this is introduced as ${\rm dexp}:so\left( 3\right)
\times so\left( 3\right) \rightarrow so\left( 3\right) $, via right
translation so that ${\rm dexp}_{\widehat{{\bf \xi }}}\dot{{\bf \xi }}=\dot{%
{\bf R}}{\bf R}^{-1}$, with ${\bf R}=\exp \widehat{{\bf \xi }}$. It admits
the series expansion 
\begin{equation}
{\rm dexp}_{\widehat{{\bf x}}}\left( \widehat{{\bf y}}\right) =\sum_{i\geq 0}%
\frac{1}{\left( i+1\right) !}{\rm ad}_{\widehat{{\bf x}}}^{i}\left( \widehat{%
{\bf y}}\right)  \label{dexp}
\end{equation}%
where ${\rm ad}_{\widehat{{\bf x}}}\left( \widehat{{\bf y}}\right) =\left[ 
\widehat{{\bf x}},\widehat{{\bf y}}\right] $. With vector notation ${\bf ad}%
_{\widehat{{\bf x}}}\widehat{{\bf y}}=\widehat{{\bf x}}{\bf y}$ this leads
to the closed form matrix expressions 
\begin{eqnarray}
{\bf dexp}_{{\bf \xi }} &=&{\bf I}+\frac{1-\cos \left\Vert {\bf \xi }%
\right\Vert }{\left\Vert {\bf \xi }\right\Vert ^{2}}\widehat{{\bf \xi }}+%
\frac{\left\Vert {\bf \xi }\right\Vert -\sin \left\Vert {\bf \xi }%
\right\Vert }{\left\Vert {\bf \xi }\right\Vert ^{3}}\widehat{{\bf \xi }}^{2}
\label{dexpSO3} \\
&=&\frac{1}{\left\Vert {\bf \xi }\right\Vert ^{2}}%
\big%
[({\bf I}-\exp \widehat{{\bf \xi }})\widehat{{\bf \xi }}+{\bf \xi \xi }^{T}%
\big%
].
\end{eqnarray}%
This matrix is frequently called tangent operator of rotations \cite%
{Simo1985}. Its significance is that it allows for expressing the body-fixed
angular velocity in terms of the time derivative of the rotation vector
since 
\begin{equation}
{\bf \omega }={\bf dexp}_{-{\bf \xi }}\dot{{\bf \xi }}
\end{equation}%
where the left Poisson equation $\widehat{{\bf \omega }}:={\bf R}^{T}\dot{%
{\bf R}}$ defines the body-fixed angular velocity vector ${\bf \omega }$.
Equally important is the inverse relation. The inverse of the dexp mapping
possesses the expansion 
\begin{equation}
{\rm dexp}_{\widehat{{\bf x}}}^{-1}\left( \widehat{{\bf y}}\right)
=\sum_{i\geq 0}\frac{B_{i}}{i!}{\rm ad}_{\widehat{{\bf x}}}^{i}\left( 
\widehat{{\bf y}}\right)  \label{dexpInvSeries}
\end{equation}%
with $B_{i}$ being the Bernoulli numbers. The inverse of (\ref{dexpSO3}) is
explicitly 
\begin{equation}
{\bf dexp}_{{\bf \xi }}^{-1}={\bf I}-\frac{1}{2}\widehat{{\bf \xi }}+\left(
1-\frac{\left\Vert {\bf \xi }\right\Vert }{2}\cot \frac{\left\Vert {\bf \xi }%
\right\Vert }{2}\right) \frac{\widehat{{\bf \xi }}^{2}}{\left\Vert {\bf \xi }%
\right\Vert ^{2}}.  \label{dexpInvSO3}
\end{equation}%
This expression seems to have first appeared in \cite{BulloMurray1995}. For
numerical implementations the second-order approximation of (\ref%
{dexpInvSeries}) 
\begin{equation}
{\rm dexp}_{\widehat{{\bf x}}}^{-1}\left( \widehat{{\bf y}}\right) \approx 
\widehat{{\bf y}}-\frac{1}{2}{\rm ad}_{\widehat{{\bf x}}}\left( \widehat{%
{\bf y}}\right) +\frac{1}{12}{\rm ad}_{\widehat{{\bf x}}}^{2}\left( \widehat{%
{\bf y}}\right)  \label{dexpInvApprox}
\end{equation}%
can be used. Interestingly, since $B_{3}=0$, this is also the third-order
approximation. Moreover, since $B_{i}=0$ for odd $i$, the approximation of
order $2i$ and $2i+1$ are identical. Application of higher order
approximations is, however, not advisable as the computational complexity
exceeds that of the closed form (\ref{dexpInvSO3}).

\subsection{Proper Rigid Body Motions and Body-Fixed Twists}

A rigid body is kinematically represented by a body-fixed reference frame.
Its configuration with respect to a world-fixed inertial reference frame is
described by the position vector ${\bf r}\in {\Bbb R}^{3}$ of its origin and
the relative rotation matrix ${\bf R}\in SO\left( 3\right) $. The {\em %
configuration} of a rigid body is thus represented by the pair $C=\left( 
{\bf R},{\bf r}\right) $, and a rigid body motion is a curve $C\left(
t\right) $. The set of such pairs may hence be regarded as the rigid body
configuration space. In order to capture the rigid body kinematics this set
must be equipped with slightly more structure. This structure has to do with
the transition from one configuration to another. The combination of two
successive rigid-body configurations is given by $C_{2}\cdot C_{1}=\left( 
{\bf R}_{2}{\bf R}_{1},{\bf r}_{2}+{\bf R}_{2}{\bf r}_{1}\right) $. The
rigid body configurations equipped with this multiplication describes frame
transformations, i.e. finite rigid body {\em motions} that constitute the
6-dimensional Lie group $SE\left( 3\right) =SO\left( 3\right) \ltimes {\Bbb R%
}^{3}$ --the special Euclidian group in 3 dimensions, i.e. the group of
isometric orientation preserving transformations of 3-dimensional Euclidian
spaces. It is crucial to note that this is the {\it semidirect} product of
the special orthogonal group $SO\left( 3\right) $ and the translation group,
represented by ${\Bbb R}^{3}$, hence the above multiplication law. This is
also encoded in the $4\times 4$ matrix representation of a rigid body motion 
\begin{equation}
{\bf C}=\left( 
\begin{array}{cc}
{\bf R} & {\bf r} \\ 
{\bf 0} & 1%
\end{array}%
\right) .  \label{C}
\end{equation}%
All matrices of the form (\ref{C}) constitute the representation of $%
SE\left( 3\right) $ as a matrix group, denoted as ${\bf C}\in SE\left(
3\right) $. The group multiplication, i.e. the concatenation of two rigid
body motions, is then given by the matrix product 
\begin{equation}
{\bf C}_{2}{\bf C}_{1}=\left( 
\begin{array}{cc}
{\bf R}_{2}{\bf R}_{1} & {\bf r}_{2}+{\bf R}_{2}{\bf r}_{1} \\ 
{\bf 0} & 1%
\end{array}%
\right) .  \label{MultSE3}
\end{equation}

To emphasize that $C\in SE\left( 3\right) $ represent frame transformations
they are occasionally termed 'proper rigid body motions'. A general motion
of a rigid body is a screw motion, i.e. an interconnected rotation and
translation along a screw axis. The velocity corresponding to the screw
motion of a rigid body is a twist characterized by the angular velocity $%
{\bf \omega }$ and the linear velocity vector ${\bf v}$, both expressed in
the body-fixed frame, summarized in the twist coordinate vector ${\bf V}%
=\left( {\bf \omega },{\bf v}\right) \in {\Bbb R}^{6}$. Given a rigid body
motion ${\bf C}\left( t\right) $, the twist of the body in body-fixed
representation is defined as 
\begin{equation}
\widehat{{\bf V}}:={\bf C}^{-1}\dot{{\bf C}}\text{ \ \ with \ \ }\widehat{%
{\bf V}}=\left( 
\begin{array}{cc}
\widehat{{\bf \omega }} & {\bf v} \\ 
{\bf 0} & 0%
\end{array}%
\right) \in se\left( 3\right)  \label{Vhat}
\end{equation}%
where ${\bf v}={\bf R}^{T}\dot{{\bf r}}$, and the matrix is given in terms
of the twist coordinate vector. Twists are interpreted as instantaneous
screws. A general screw coordinate vector is denoted as ${\bf X}=\left( {\bf %
\xi },{\bf \eta }\right) \in {\Bbb R}^{6}$. The vector space of matrices of
the form (\ref{Vhat}) forms the matrix representation of the Lie algebra $%
se\left( 3\right) $, that is the algebra of screws. The hat operator $%
\widehat{\cdot }$ is used interchangeably to produce either $so\left(
3\right) $ or $se\left( 3\right) $ matrices. Hence, any $se\left( 3\right) $%
-matrix is given in terms of some twist coordinates.

The Lie bracket on $se\left( 3\right) $ is the matrix commutator that
assumes the explicit form 
\begin{equation}
\lbrack \widehat{{\bf X}}_{1},\widehat{{\bf X}}_{2}]=\widehat{{\bf X}}_{1}%
\widehat{{\bf X}}_{2}-\widehat{{\bf X}}_{2}\widehat{{\bf X}}_{1}=\left( 
\begin{array}{cc}
\widehat{{\bf \xi }_{1}\times {\bf \xi }_{2}} & {\bf \xi }_{1}\times {\bf %
\eta }_{2}-{\bf \xi }_{2}\times {\bf \eta }_{1} \\ 
{\bf 0} & 0%
\end{array}%
\right) .  \label{se3bracket}
\end{equation}%
In other words, the Lie bracket of two twists is given by the screw product $%
\left[ {\bf X}_{1},{\bf X}_{2}\right] =\left( {\bf \xi }_{1}\times {\bf \xi }%
_{2},{\bf \xi }_{1}\times {\bf \eta }_{2}-{\bf \xi }_{2}\times {\bf \eta }%
_{1}\right) $ \cite{SeligBook}. This reflects the fact that $se\left(
3\right) =so\left( 3\right) \oplus _{s}{\Bbb R}^{3}$ is the semidirect sum
of $so\left( 3\right) $ and ${\Bbb R}^{3}$. The Lie bracket also embodies
the adjoint action of $se\left( 3\right) $ on itself, defined as ${\rm ad}_{%
\widehat{{\bf X}}_{1}}\widehat{{\bf X}}_{2}=[\widehat{{\bf X}}_{1},\widehat{%
{\bf X}}_{2}]$. Being a linear operator it can be expressed as a matrix
acting on screw coordinate vector ${\bf ad}_{{\bf X}_{1}}{\bf X}_{2}$. For a
screw ${\bf X}\in se\left( 3\right) $ this matrix is 
\begin{equation}
{\bf ad}_{{\bf X}}=\left( 
\begin{array}{cc}
\widehat{{\bf \xi }} & \ \ {\bf 0} \\ 
\widehat{{\bf \eta }} & \ \ \widehat{{\bf \xi }}%
\end{array}%
\right) .  \label{adse3}
\end{equation}%
This is the matrix representation of the Lie bracket of two twists.

$se\left( 3\right) $ is isomorphic to the algebra of screws via the
identification (\ref{Vhat}). A screw coordinate vector ${\bf X}=\left( {\bf %
\xi },{\bf \eta }\right) $ describes an instantaneous screw motion, i.e. a
rotation about the axis ${\bf \xi }$ together with a translation along this
axis. The latter is given in terms of the position vector ${\bf r}$ of a
point on the screw as ${\bf \eta }={\bf r}\times {\bf \xi }+h_{{\bf X}}{\bf %
\xi }$, where the scalar $h_{{\bf X}}$ is the pitch of the screw. If $h_{%
{\bf X}}=0$, then ${\bf X}$ are simply the Pl\"{u}cker coordinates of a line
along the screw axis. Since any rigid body motion is a finite screw motion
it can be expressed in terms of some (time dependent) instantaneous screw
coordinates via the exponential mapping.

The exponential mapping, generating $SE\left( 3\right) $ from its Lie
algebra $se\left( 3\right) $, i.e. the matrix exponential that gives rigid
body configuration (\ref{C}) in terms of screw coordinates, admits the
explicit form 
\begin{eqnarray}
{\bf X}=\left( {\bf \xi },{\bf \eta }\right) \longmapsto \exp \widehat{{\bf X%
}}{} &{=}&\left( 
\begin{array}{cc}
\exp \widehat{{\bf \xi }} & \ \;{\bf dexp}_{{\bf \xi }}{\bf \eta }%
\vspace{2mm}
\\ 
{\bf 0} & 1%
\end{array}%
\right)  \label{expX} \\
&=&\left( 
\begin{array}{cc}
\exp \widehat{{\bf \xi }} & \;\;\;\frac{1}{\left\Vert {\bf \xi }\right\Vert
^{2}}\left( {\bf I}-\exp \widehat{{\bf \xi }}\right) \left( {\bf \xi }\times 
{\bf \eta }\right) +h_{{\bf X}}\,{\bf \xi }%
\vspace{2mm}
\\ 
{\bf 0} & 1%
\end{array}%
\right)  \label{expX2}
\end{eqnarray}%
with pitch $h_{{\bf X}}:={\bf \xi }\cdot {\bf \eta }/\left\Vert {\bf \xi }%
\right\Vert ^{2}$, and ${\bf dexp}_{{\bf \xi }}$ from (\ref{dexpSO3}). The
exponential mapping (\ref{expX}) determines the finite screw motion as a
result of the time evolution of the instantaneous screw ${\bf X}\left(
t\right) $.

The corresponding twist is determined by the (right) dexp mapping ${\rm dexp}%
:se\left( 3\right) \times se\left( 3\right) \rightarrow se\left( 3\right) $
as 
\begin{equation}
\widehat{{\bf V}}={\rm dexp}_{-\widehat{{\bf X}}}(\dot{\widehat{{\bf X}}}).
\label{VdexpSE3}
\end{equation}%
The expression (\ref{dexp}), being applicable to any Lie group, gives rise
to several closed forms. A straightforward application of (\ref{dexp}) using
(\ref{adse3}) yields the matrix form \cite%
{Bottasso1998,Borri2001a,ParkChung2005}%
\begin{equation}
{\bf dexp}_{{\bf X}}=\left( 
\begin{array}{cc}
{\bf dexp}_{{\bf \xi }} & {\bf 0}%
\vspace{2mm}
\\ 
{\bf P} & {\bf dexp}_{{\bf \xi }}%
\end{array}%
\right)  \label{dexpSE3}
\end{equation}%
for ${\bf X}=\left( {\bf \xi },{\bf \eta }\right) $, where 
\begin{equation}
{\bf P}\left( {\bf X}\right) =\frac{\beta }{2}\widehat{{\bf \eta }}+\frac{%
1-\alpha }{\left\Vert {\bf \xi }\right\Vert ^{2}}\left( \widehat{{\bf \eta }}%
\widehat{{\bf \xi }}+\widehat{{\bf \xi }}\widehat{{\bf \eta }}\right) +h_{%
{\bf X}}\frac{\alpha -\beta }{\left\Vert {\bf \xi }\right\Vert }\widehat{%
{\bf \xi }}+\frac{h_{{\bf X}}}{\left\Vert {\bf \xi }\right\Vert ^{2}}\left( 
\frac{\beta }{2}-\frac{3(1-\alpha )}{\left\Vert {\bf \xi }\right\Vert }%
\right) \widehat{{\bf \xi }}^{2}  \label{P1}
\end{equation}%
with $\alpha :=\frac{2}{\left\Vert {\bf \xi }\right\Vert }\sin \frac{%
\left\Vert {\bf \xi }\right\Vert }{2}\cos \frac{\left\Vert {\bf \xi }%
\right\Vert }{2},\beta :=\frac{4}{\left\Vert {\bf \xi }\right\Vert ^{2}}\sin
^{2}\frac{\left\Vert {\bf \xi }\right\Vert }{2}$, and the pitch $h={\bf \xi }%
\cdot {\bf \eta }/\left\Vert {\bf \xi }\right\Vert ^{2}$. For pure rotation,
i.e. $h_{{\bf X}}=0$, (\ref{P1}) simplifies. Therewith the vector of the
body-fixed twist is given in terms of the time derivative of ${\bf X}$ as $%
{\bf V}={\bf dexp}_{-{\bf X}}\dot{{\bf X}}$.

The inverse of the dexp mapping matrix is readily found as 
\begin{equation}
{\bf dexp}_{{\bf X}}^{-1}=\left( 
\begin{array}{cc}
{\bf dexp}_{{\bf \xi }}^{-1} & {\bf 0}%
\vspace{2mm}
\\ 
{\bf U} & {\bf dexp}_{{\bf \xi }}^{-1}%
\end{array}%
\right)  \label{dexpInvSE3Park}
\end{equation}%
with 
\begin{equation}
{\bf U}\left( {\bf X}\right) =\frac{1-\gamma }{\left\Vert {\bf \xi }%
\right\Vert ^{2}}\left( \widehat{{\bf \eta }}\widehat{{\bf \xi }}+\widehat{%
{\bf \xi }}\widehat{{\bf \eta }}\right) +\frac{h_{{\bf X}}}{\left\Vert {\bf %
\xi }\right\Vert ^{3}}\left( \frac{1}{\beta }+\gamma -2\right) \widehat{{\bf %
\xi }}^{2}-\frac{1}{2}\widehat{{\bf \eta }}
\end{equation}%
and $\gamma :=\frac{2}{\left\Vert {\bf \xi }\right\Vert }\cot \frac{%
\left\Vert {\bf \xi }\right\Vert }{2}$. Another closed form expression was
reported in \cite{SeligBook} 
\begin{equation}
{\bf dexp}_{{\bf X}}^{-1}={\bf I}-\frac{1}{2}{\bf ad}_{{\bf X}}+\left( \frac{%
2}{\left\Vert {\bf \xi }\right\Vert ^{2}}+\frac{\left\Vert {\bf \xi }%
\right\Vert +3\sin \left\Vert {\bf \xi }\right\Vert }{4\left\Vert {\bf \xi }%
\right\Vert \left( \cos \left\Vert {\bf \xi }\right\Vert -1\right) }\right) 
{\bf ad}_{{\bf X}}^{2}+\left( \frac{1}{\left\Vert {\bf \xi }\right\Vert ^{4}}%
+\frac{\left\Vert {\bf \xi }\right\Vert +\sin \left\Vert {\bf \xi }%
\right\Vert }{4\left\Vert {\bf \xi }\right\Vert ^{3}\left( \cos \left\Vert 
{\bf \xi }\right\Vert -1\right) }\right) {\bf ad}_{{\bf X}}^{4}.
\label{dexpInvSE3}
\end{equation}%
The formulations (\ref{dexpSE3}) and (\ref{dexpInvSE3}) admit reformulating
the kinematic relation as explicit ODE $\dot{{\bf X}}={\bf dexp}_{-{\bf X}%
}^{-1}{\bf V}$.

\begin{remark}
It is instructive to show to show how the 3-parameter description of
rotations fits into the Lie group formulation. If rotations are described by
successive rotations about intermediate axes the rotation matrix is given as 
${\bf R}=\exp (\widehat{{\bf e}}_{i}\theta _{i})\exp (\widehat{{\bf e}}%
_{j}\theta _{j})\exp (\widehat{{\bf e}}_{k}\theta _{k})$, where ${\bf e}_{i},%
{\bf e}_{j},{\bf e}_{k}$ are constant unit vectors along the intermediate
instantaneous rotation axes, and $\theta _{i},\theta _{j},\theta _{k}$ are
the corresponding rotation angles. The only condition is ${\bf e}_{i}\neq 
{\bf e}_{j},{\bf e}_{j}\neq {\bf e}_{k}$. In particular, using $i=k=3,j=1$
corresponds to the Euler angle description, which is called 'degenerate'
since ${\bf e}_{i}={\bf e}_{k}$. Setting $i=1,j=2,k=3$ yields the
Bryant/Cardan angle parameterization. Noting that for constant ${\bf e}_{i}$
in ${\bf R}_{i}=\exp (\widehat{{\bf e}}_{i}\theta _{i})$, ${\bf R}_{i}^{-1}%
\dot{{\bf R}}_{i}={\rm dexp}_{{\bf e}_{i}\theta _{i}}(\widehat{{\bf e}}_{i}%
\dot{\theta}_{i})=\widehat{{\bf e}}_{i}\dot{\theta}_{i}$, etc., the
corresponding angular velocity is 
\begin{equation}
{\bf \omega }=\left( 
\begin{array}{ccc}
&  &  \\ 
{\bf R}_{k}^{T}{\bf R}_{j}^{T}{\bf e}_{i}\ \ \vdots & {\bf R}_{k}^{T}{\bf e}%
_{j}\ \ \vdots & {\bf e}_{k} \\ 
&  & 
\end{array}%
\right) \left( 
\begin{array}{c}
\dot{\theta}_{i} \\ 
\dot{\theta}_{j} \\ 
\dot{\theta}_{k}%
\end{array}%
\right) ={\bf B\dot{\theta}}  \label{3param}
\end{equation}%
as in (\ref{VelRigidBody}). E.g. for $i=k=3,j=1$, these are the kinematic
Euler equations, and the $3\times 3$ matrix in (\ref{3param}) is the
corresponding coefficient matrix \cite{McCauley}. This is easily extended to 
$SE\left( 3\right) $ by using a unit basis including translations.
\end{remark}

\subsection{Direct Product Representation of Rigid Body Configurations and
Mixed Velocity Representation%
\label{secHybrid}%
\label{secHybrid}}

If the concatenation of frame transformations is ignored, a rigid body
configuration can simply be regarded as $C=\left( {\bf R},{\bf r}\right) \in
SO\left( 3\right) \times {\Bbb R}^{3}$. That is, the configuration space of
a rigid body is the direct product of $SO\left( 3\right) $ and ${\Bbb R}^{3}$
with group multiplication 
\begin{equation}
C_{1}\cdot C_{2}=({\bf R}_{1}{\bf R}_{2},{\bf r}_{1}+{\bf r}_{2}).
\label{MultDirectProduct}
\end{equation}%
This does clearly not represent a frame transformation, i.e. a rigid body
motion. Nevertheless it is commonly used as geometric model for MBS
kinematics giving rise to the equations (\ref{VelRigidBody}), and for Lie
group integration methods of rigid body dynamics \cite%
{BruelsCardona2010,BruelsCardonaArnold2012,CelledoniOwren1999,Krysl,TerzeIMSD2012}%
. The drawback becomes critical if the multiplication (\ref%
{MultDirectProduct}) is employed in the position update of integration
schemes as will be shown later.

The inverse element is $({\bf R},{\bf r})^{-1}=({\bf R}^{T},-{\bf r})$. If
desired, $SO\left( 3\right) \times {\Bbb R}^{3}$ can be represented by the
group of $7\times 7$ matrices of the form%
\begin{equation}
{\bf C}=\left( 
\begin{array}{ccc}
{\bf R} & {\bf 0} & {\bf 0} \\ 
{\bf 0} & {\bf I} & {\bf r} \\ 
{\bf 0} & {\bf 0} & 1%
\end{array}%
\right) ,\text{with\thinspace\ }{\bf R}\in SO(3),{\bf r}\in {\Bbb R}^{3}.
\label{MatrixSO3xR3}
\end{equation}%
The product (\ref{MultDirectProduct}) is then represented by ${\bf C}%
_{1}\cdot {\bf C}_{2}$. The Lie algebra of the direct product $SO\left(
3\right) \times {\Bbb R}^{3}$ is the direct sum $so\left( 3\right) \oplus 
{\Bbb R}^{3}$ whose elements can also be represented as vectors ${\bf X}%
=\left( {\bf \xi },{\bf r}\right) \in {\Bbb R}^{6}$. With ${\Bbb R}^{3}$
being a commutative algebra, the Lie bracket on this algebra is 
\begin{equation}
\left[ {\bf X}_{1},{\bf X}_{2}\right] =\left( {\bf \xi }_{1}\times {\bf \xi }%
_{2},{\bf 0}\right) .  \label{bracketSO3xR3}
\end{equation}%
This can be expressed as ${\bf ad}_{{\bf X}_{2}}{\bf X}_{1}$ with matrix%
\begin{equation}
{\bf ad}_{{\bf X}}=\left( 
\begin{array}{cc}
\widehat{{\bf \xi }} & \ \ {\bf 0} \\ 
{\bf 0} & \ \ {\bf 0}%
\end{array}%
\right) .  \label{adso3xr3}
\end{equation}%
If the matrix representation (\ref{MatrixSO3xR3}) is used, the corresponding
Lie algebra consists of matrices%
\begin{equation}
\widehat{{\bf X}}=\left( 
\begin{array}{ccc}
\widehat{{\bf \xi }} & {\bf 0} & {\bf 0} \\ 
{\bf 0} & {\bf 0} & {\bf r} \\ 
{\bf 0} & {\bf 0} & 0%
\end{array}%
\right) ,\text{with\thinspace\ }\widehat{{\bf \xi }}\in so(3),{\bf r}\in 
{\Bbb R}^{3}  \label{Matrixso3xR3}
\end{equation}%
with matrix commutator as Lie bracket. The exponential mapping on the direct
product group is 
\begin{equation}
{\bf X}=\left( {\bf \xi },{\bf r}\right) \longmapsto \exp {\bf X}=(\exp 
\widehat{{\bf \xi }},{\bf r})  \label{expSO3xR3}
\end{equation}%
with the exponential mapping (\ref{expSO3}) on $SO\left( 3\right) $.
Apparently ${\bf X}$ is not a screw coordinate vector, but rather consists
of a rotation and a translation vector. The exponential of (\ref%
{Matrixso3xR3}) yields the matrices (\ref{MatrixSO3xR3}).

The rigid body velocity for $C\left( t\right) =\left( {\bf R},{\bf r}\right)
\in SO\left( 3\right) \times {\Bbb R}^{3}$ is introduced in vector notation
as ${\bf V}^{\text{m}}=\left( {\bf \omega },{\bf v}^{\text{s}}\right) $ with 
${\bf v}^{\text{s}}:=\dot{{\bf r}}$. This is not a proper twist. It contains
an apparent mix of body-fixed angular velocity ${\bf \omega }$ and spatial
linear velocity ${\bf v}^{\text{s}}\equiv \dot{{\bf r}}$, and is therefore
commonly referred to as {\it mixed representation} of rigid body velocities.
With (\ref{Matrixso3xR3}) this is, in analogy to (\ref{Vhat}), defined by 
\begin{equation}
\widehat{{\bf V}}^{\text{m}}:={\bf C}^{-1}\dot{{\bf C}}\text{ \ \ with \ \ }%
\widehat{{\bf V}}^{\text{m}}=\left( 
\begin{array}{ccc}
\widehat{{\bf \omega }} & {\bf 0} & {\bf 0} \\ 
{\bf 0} & {\bf 0} & {\bf v}^{\text{s}} \\ 
{\bf 0} & {\bf 0} & 0%
\end{array}%
\right) \in so\left( 3\right) \times {\Bbb R}^{3}.  \label{hybridvel}
\end{equation}%
The dexp mapping corresponding to (\ref{expSO3xR3}) is 
\begin{equation}
{\rm dexp}_{{\bf X}_{1}}({\bf X}_{2})=({\rm dexp}_{\widehat{{\bf \xi }}_{1}}(%
\widehat{{\bf \xi }}_{2}),{\bf r}_{2})  \label{dexpSO3xR3}
\end{equation}%
with the dexp mapping on $SO\left( 3\right) $ in (\ref{dexpSO3}). In vector
form of the mixed velocity, the relation%
\begin{equation}
{\bf V}^{\text{m}}={\bf dexp}_{-{\bf X}}\dot{{\bf X}}=\left( 
\begin{array}{cc}
{\bf dexp}_{-\widehat{{\bf \xi }}} & \ \ {\bf 0}%
\vspace{2mm}
\\ 
{\bf 0} & \ \ {\bf I}%
\end{array}%
\right) \left( 
\begin{array}{c}
\dot{{\bf \xi }} \\ 
\dot{{\bf r}}%
\end{array}%
\right)  \label{Vhdexp}
\end{equation}%
for ${\bf X}=\left( {\bf \xi },{\bf r}\right) $, resembles (\ref%
{VelRigidBody}). In fact the angular and linear velocities are considered as
decoupled. The matrix ${\bf dexp}_{-\widehat{{\bf \xi }}}$ in (\ref{Vhdexp})
is the matrix ${\bf B}$ in (\ref{VelRigidBody}) when the rotation vector
parameterization is used. Moreover, (\ref{Vhdexp}) covers any general
3-angle parameterization by concatenation of three successive mappings as in
(\ref{3param}).

\begin{remark}
In applications frequently the mixed velocity ${\bf V}^{\text{m}}=\left( 
{\bf \omega },\dot{{\bf r}}\right) $ is used (section \ref{secHybrid}). If
desired, the body-fixed twist can be transformed to the mixed velocity.
Since $\dot{{\bf r}}={\bf Rv}$ this is achieved by premultiplication of the
second column in (\ref{dexpSE3}) with ${\bf R}=\exp \widehat{{\bf \xi }}$.
Noting that $\exp \widehat{{\bf \xi }}\cdot {\bf dexp}_{-\widehat{{\bf \xi }}%
}={\bf dexp}_{\widehat{{\bf \xi }}}$ this leads to the expression ${\bf V}^{%
\text{m}}={\bf A}\left( {\bf X}\right) \dot{{\bf X}}$ with%
\begin{equation}
{\bf A}\left( {\bf X}\right) :=\left( 
\begin{array}{cc}
{\bf I} & \ \ {\bf 0}%
\vspace{2mm}
\\ 
{\bf 0} & \ \ {\bf R}%
\end{array}%
\right) \cdot {\bf dexp}_{-{\bf X}}=\left( 
\begin{array}{cc}
{\bf dexp}_{-{\bf \xi }} & \ \ {\bf 0}%
\vspace{2mm}
\\ 
{\bf R}\cdot {\bf P}\left( -{\bf X}\right) & \ \ {\bf dexp}_{{\bf \xi }}%
\end{array}%
\right) .  \label{hybridtwist}
\end{equation}
\end{remark}

\begin{remark}
The matrix representation (\ref{hybridvel}) is not used in implementations
but merely introduced to emphasize the formally identical definition of
velocity (\ref{Vhat}) and (\ref{hybridvel}) for both groups, used in the
general form of the motion equations in section \ref{secEOMLieGroup}. The
kinematic reconstruction equations (\ref{Videxp}) are generally valid and
only the c-space Lie group is to be replaced. The velocities are
left-invariant vector fields. For $SE(3)$ this means that they are invariant
w.r.t. changes of global reference frame (described by left multiplication
of $C$).
\end{remark}

In summary the only difference of the two c-space representations is the
concatenation of configurations, but this is the critical aspect for
numerical reconstruction of finite motions from velocities where the
solution is advanced from one time step to the next by an incremental
motion. In other words, both, $SE\left( 3\right) $ and $SO\left( 3\right)
\times {\Bbb R}^{3}$, can be used to represent the {\it configuration} of a
rigid body but only $SE\left( 3\right) $ allows for representing rigid body 
{\it motions}. In particular successive frame transformations, i.e.
concatenation of successive configurations, (\ref{MultSE3}) and (\ref%
{MultDirectProduct}), respectively, reveal the fundamental difference, which
is accordingly reflected in the exp mapping.

\subsection{Proper Rigid Body Motions in Terms of Dual Quaternions%
\label{secDualQuat}%
\label{secDualQuat}}

The fact that there is no 3-parametric global parameterization of rotations
leads to the well-known problem of parameterization singularities on $%
SO\left( 3\right) $ and any product group. This is apparent from the
axis-angle description, i.e. using exponential coordinates ${\bf \xi }$,
since (\ref{dexpSO3}), respectively (\ref{dexpSE3}), is singular for $%
\left\Vert {\bf \xi }\right\Vert =\pm k\pi ,k\in {\Bbb N}$. This can only be
avoided by using redundant parameters like Euler parameters (instead of
three independent in ${\bf \xi }$), or by avoiding the introduction of local
coordinates at all (see appendix B). The use of dependent global coordinates
corresponds to replace $SE\left( 3\right) $ and $SO\left( 3\right) \times 
{\Bbb R}^{3}$ by their covering Lie groups. This is outlined in the this and
next section.

The rigid body motion group $SE\left( 3\right) $ is homomorphic to the group
of dual quaternions, ${\Bbb H}_{\varepsilon }$ \cite{MccarthyBook,Yang1964}.
An ordinary quaternion is written as 4-vector ${\bf Q}=\left( q_{0},{\bf q}%
\right) \in {\Bbb H}$ where the vector part is ${\bf q}=\left(
q_{1},q_{2},q_{3}\right) $. Euler parameters are unit quaternions, i.e. $%
\left\Vert {\bf Q}\right\Vert =1$, used to parameterize spatial rotations.
Dual quaternions are generalizations of ordinary quaternions. Avoiding using
dual numbers, a dual quaternion can be expressed as an 8-dimensional vector $%
\hat{{\bf Q}}=\left( {\bf Q},{\bf Q}_{\varepsilon }\right) $, where ${\bf Q}%
\in {\Bbb H}$ is a unit quaternion and ${\bf Q}_{\varepsilon }$ is another
quaternion satisfying the Pl\"{u}cker condition ${\bf Q}\cdot {\bf Q}%
_{\varepsilon }=0$. The Euler parameters ${\bf Q}$ describe the rotation and 
${\bf Q}_{\varepsilon }$ the translation of a frame transformation. That is,
a frame transformation $C=\left( {\bf R},{\bf r}\right) \in SE\left(
3\right) $, and thus its screw coordinates ${\bf X}=\left( {\bf \xi },{\bf %
\eta }\right) $, can be mapped to a dual quaternion. This is called the
kinematic image space transformation or kinematic mapping \cite%
{Blaschke1942,BottemaRoth1979,RavaniRoth1984,Study1891}.

Multiplication of dual quaternions is expressible in matrix form as \cite%
{DooleyMcCarthy1991} 
\begin{equation}
\hat{{\bf Q}}\hat{{\bf P}}=\hat{{\bf Q}}^{+}\hat{{\bf P}}=\hat{{\bf P}}^{-}%
\hat{{\bf Q}}
\end{equation}%
with 
\begin{equation}
\hat{{\bf Q}}^{+}=\left( 
\begin{array}{cc}
{\bf Q}^{+} & {\bf 0} \\ 
{\bf Q}_{\varepsilon }^{+} & {\bf Q}^{+}%
\end{array}%
\right) ,\ \ \ \hat{{\bf Q}}^{-}=\left( 
\begin{array}{cc}
{\bf Q}^{-} & {\bf 0} \\ 
{\bf Q}_{\varepsilon }^{-} & {\bf Q}^{-}%
\end{array}%
\right)
\end{equation}%
and the Hamilton matrices 
\begin{equation}
{\bf Q}^{+}=\left( 
\begin{array}{cc}
q_{0} & -{\bf q}^{T} \\ 
{\bf q} & q_{0}{\bf I}_{3}+\widehat{{\bf q}}%
\end{array}%
\right) ,\ \ \ {\bf Q}^{-}=\left( 
\begin{array}{cc}
q_{0} & -{\bf q}^{T} \\ 
{\bf q} & q_{0}{\bf I}_{3}-\widehat{{\bf q}}%
\end{array}%
\right) .
\end{equation}%
The rotation matrix corresponding to ${\bf Q}$ is given as ${\bf R}={\bf DE}%
^{T}$ with 
\begin{equation}
{\bf D}\left( {\bf Q}\right) =\left( 
\begin{array}{cc}
-{\bf q\ \ } & q_{0}{\bf I}_{3}+\widehat{{\bf q}}%
\end{array}%
\right) ,\ \ \ {\bf E}\left( {\bf Q}\right) =\left( 
\begin{array}{cc}
-{\bf q\ \ } & q_{0}{\bf I}_{3}-\widehat{{\bf q}}%
\end{array}%
\right) .  \label{DE}
\end{equation}%
The corresponding Cartesian position vector ${\bf r}$ is given by the vector
part of 
\begin{equation}
\left( 
\begin{array}{c}
0 \\ 
{\bf r}%
\end{array}%
\right) =2{\bf Q}^{-T}{\bf Q}_{\varepsilon }.  \label{Qpos}
\end{equation}%
The mixed velocity ${\bf V}^{\text{m}}=\left( {\bf \omega },{\bf v}^{\text{s}%
}\right) $ is then given as 
\begin{equation}
{\bf V}^{\text{m}}={\bf H}^{\text{m}}(\hat{{\bf Q}})\dot{\hat{{\bf Q}}}
\label{hybridtwistDualQuat}
\end{equation}%
with 
\begin{equation}
{\bf H}^{\text{m}}(\hat{{\bf Q}}):=2\left( 
\begin{array}{cc}
{\bf D}\left( {\bf Q}\right) & {\bf 0}%
\vspace{2mm}
\\ 
-{\bf D}\left( {\bf Q}_{\varepsilon }\right) & \ \ \ {\bf D}\left( {\bf Q}%
\right)%
\end{array}%
\right) .  \label{Hhdual}
\end{equation}%
Using ${\bf v}={\bf ED}^{T}{\bf v}^{\text{s}}$ and ${\bf ED}^{T}{\bf D}%
=\left\Vert {\bf Q}\right\Vert {\bf E}$ with $\left\Vert {\bf Q}\right\Vert
=1$, the body-fixed twist is 
\begin{equation}
{\bf V}={\bf H}(\hat{{\bf Q}})\dot{\hat{{\bf Q}}}  \label{twistDualQuat}
\end{equation}%
with 
\begin{equation}
{\bf H}(\hat{{\bf Q}}):=\left( 
\begin{array}{cc}
{\bf D}\left( {\bf Q}\right) & \ \ \ {\bf 0}%
\vspace{2mm}
\\ 
-{\bf E}\left( {\bf Q}\right) {\bf D}^{T}\left( {\bf Q}\right) {\bf D}\left( 
{\bf Q}_{\varepsilon }\right) & \ \ \ {\bf E}\left( {\bf Q}\right)%
\end{array}%
\right) .  \label{Hdual}
\end{equation}

In summary, a frame transformation (a proper rigid body motion) is
represented by a dual quaternion, and the screw coordinates ${\bf X}=\left( 
{\bf \xi },{\bf \eta }\right) \in {\Bbb R}^{6}$ are replaced by the dual
quaternion ${\bf X}=\left( {\bf Q},{\bf Q}_{\varepsilon }\right) \in {\Bbb H}%
_{\varepsilon }$ that serve as dependent global coordinates. The two $\left( 
{\bf Q},{\bf Q}_{\varepsilon }\right) $ and $\left( {\bf Q}^{\ast },{\bf Q}%
_{\varepsilon }\right) $ correspond to the same frame transformation $C\in
SE\left( 3\right) $, where ${\bf Q}^{\ast }$ is the conjugate of ${\bf Q}$.
Relation (\ref{twistDualQuat}) and (\ref{hybridtwistDualQuat}) is the
counterpart of (\ref{dexpSE3}) and (\ref{hybridtwist}), respectively.

\begin{remark}
It may seem unnecessary to replace both, the parameterization of rotation
and translation by quaternions. But this is indeed necessary in order to
retain the semidirect product structure of $SE\left( 3\right) $, i.e. the
coupling of rotation and translation. This is apparent from (\ref{dexpSE3})
since even if the upper left dexp mapping is replaced by ${\bf D}$ the one
in the lower right corner cannot be replaced.
\end{remark}

\subsection{Direct Product Representation in Terms of Euler Parameters%
\label{secQuat}%
\label{secQuat}}

The substitution of the rotation parameters in the direct product group is
straightforward. $SO\left( 3\right) \times {\Bbb R}^{3}$ is homomorphic to $%
{\Bbb H}\times {\Bbb R}^{3}$, and Euler parameters can be introduced
immediately as global coordinates on $SO\left( 3\right) $. Denote ${\bf X}%
=\left( {\bf Q},{\bf r}\right) \in {\Bbb H}\times {\Bbb R}^{3}$, then%
\begin{equation}
{\bf V}^{\text{m}}={\bf H}^{\text{m}}({\bf X})\dot{{\bf X}}
\end{equation}%
with 
\begin{equation}
{\bf H}^{\text{m}}({\bf Q},{\bf r}):=2\left( 
\begin{array}{cc}
{\bf D}\left( {\bf Q}\right) & {\bf 0}%
\vspace{2mm}
\\ 
{\bf 0} & {\bf I}_{3}%
\end{array}%
\right)  \label{Hh}
\end{equation}%
replacing (\ref{Vhdexp}). The direct product $SO\left( 3\right) \times {\Bbb %
R}^{3}$ is thus replaced by the direct product ${\Bbb H}\times {\Bbb R}^{3}$.

\section{Munthe-Kaas Integration Method for Kinematic Reconstruction on
C-Space Lie Groups using Local Coordinates}

The vector space ODEs (\ref{ODEvec1}) and (\ref{ODEvec2}) can be solved with
any established numerical integration scheme. Integration of (\ref{ODE2})
requires application of Lie group integration schemes. Several numerical
schemes were developed for solving system of the form (\ref{ODE2}b).
Munthe-Kaas (MK) and Crouch-Grossman methods are the best-known \cite%
{Iserles2000,chrouchgrossman1993}. The basic idea behind the MK integration
scheme is to replace the equation (\ref{ODE2}b) on $G$ by a vector space ODE
on ${\frak g}$ and to solve this by a standard vector space integration
method. The MK method was originally introduced for right translated
equations, i.e. of the form $\dot{g}=\widehat{{\bf V}}g$ \cite%
{Marthinsen1997,muntekaas1998,muntekaas1999,MuntheKaasOwren1999,OwrenMarthinsen1999}%
. It is adapted to the left translated form (\ref{ODE2}b) in the following.
Starting from an initial value $g_{0}\in G$ a solution of (\ref{ODE2}b) can
be expressed in the form $g\left( t\right) =g_{0}\exp \widehat{{\bf \Phi }}%
\left( t\right) $, with ${\bf \Phi }\left( 0\right) =0$, where ${\bf \Phi }%
\left( t\right) ${\bf \ }is a curve in ${\frak g}$. It is know that ${\bf %
\Phi }$ must satisfy the linear ODE\vspace{-4ex}

\begin{equation}
\dot{{\bf \Phi }}\left( t\right) ={\rm dexp}_{-{\bf \Phi }\left( t\right)
}^{-1}{\bf V}\left( t,g\left( t\right) \right) ,\text{ with }{\bf \Phi }%
\left( 0\right) =0  \label{ODE2left}
\end{equation}
which is commonly attributed to Hausdorff \cite{Hairer2006,magnus54}. Hence,
although it is not involved in the formulation of the reconstruction
equations (\ref{ODE2}b), the dexp mapping appears again in the integration
scheme. In the MK method (\ref{ODE2left}) is solved with an explicit
Runge-Kutta (RK) method.

At the integration step $i$, i.e. in the transition from time $t_{i-1}$ to $%
t_{i}$, the system 
\begin{equation}
\dot{{\bf \Phi }}={\rm dexp}_{-\Phi }^{-1}{\bf V}(t,g_{i-1}\exp {\bf \Phi }%
),\ t\in \lbrack t_{i-1},t_{i}],\text{ with }{\bf \Phi }\left(
t_{i-1}\right) =0  \label{subst1}
\end{equation}
is solved with a RK scheme to obtain a solution ${\bf \Phi }^{\left(
i\right) }:={\bf \Phi }\left( t_{i}\right) $. The numerical solution of (\ref%
{ODE2}b) is then advanced as $g_{i}:=g_{i-1}\exp {\bf \Phi }^{\left(
i\right) }$. $g_{0}$ is the initial configuration of the MBS (satisfying the
assembly constraints). Resorting to the RK method an explicit $s$-stage MK
scheme at time step $i$ with step size $\Delta t$ follows immediately as 
\begin{eqnarray}
g_{i}&{:=} &g_{i-1}\exp \widehat{{\bf \Phi }}^{\left( i\right) }  \nonumber
\\
{\bf \Phi }^{\left( i\right) }& {:=} &\Delta t\sum_{j=1}^{s}b_{j}{\bf k}_{j}
\nonumber \\
{\bf k}_{j}&{:=} &{\rm dexp}_{-{\bf \Psi }_{j}}^{-1}{\bf V}%
(t_{i-1}+c_{j}\Delta t,g_{i-1}\exp \widehat{{\bf \Psi }}_{j})  \label{MK1} \\
{\bf \Psi }_{j} &{:=}&\Delta t\sum_{l=1}^{j-1}a_{jl}{\bf k}_{l},\ {\bf \Psi }%
_{1}=0  \nonumber
\end{eqnarray}
where $a_{jl},b_{j}$, and $c_{j}$ are the Butcher coefficients of the $s$%
-stage RK method, and ${\bf k}_{j},{\bf \Psi }_{j}\in {\frak g}$.

The equations (\ref{ODE2}b) are indeed not independent and dynamics
simulation requires integrating the dynamics equations (\ref{ODE2}a) with
the same RK scheme, so that the overall system (\ref{ODE2}) is integrated
with one RK scheme. Moreover, since the MK method reduces to the vector
space RK scheme, the motion equations (\ref{ODE2}) can be treated as a
system on a state space Lie group and integrated with the MK method. This
will be reported in a forthcoming publication giving special attention to
constraint stabilization. Furthermore, an appropriate integration method can
be used to integrate (\ref{ODE2}) according to their numerical conditioning,
and the RK method be replaced accordingly in the MK scheme.

\begin{remark}
Although the ODE (\ref{ODE2}) is coordinate-free the ${\bf \Phi }^{\left(
i\right) }$ serve as local coordinates on the c-space. Thus a coordinate
patch on $G$ is defined during the integration. This parameterization can be
assumed singularity-free since each ${\bf \Phi }^{\left( i\right) }$ need
only cover a finite neighborhood of $g_{i-1}$, depending on the step size
and the actual dynamics.
\end{remark}

\begin{remark}
Notice that integration of the motion equations (\ref{BH2}) with a RK method
is equivalent to the MK scheme (\ref{MK1}) where global exponential
coordinates ${\bf q}$ are used, instead of local coordinates ${\bf \Phi }$,
and the configuration update is $g_{i}{:=}g_{0}\exp {\bf q}^{\left( i\right)
}$. This underlines again that the Lie group concept is just another way of
looking at the rigid body kinematics.
\end{remark}

\end{document}